\documentclass{article}
\usepackage{amsmath,enumerate,amsfonts,amssymb,color,graphicx}

\def\RR{{\mathbb R}}

\def\SS{{\mathbb S}}

\def\FS{{\mathfrak S}}

\newtheorem{theorem}{Theorem}
\newtheorem{thm}{Theorem}
\newtheorem{proposition}{Proposition}
\newtheorem{corollary}{Corollary}
\newtheorem{lemma}{Lemma}
\newtheorem{lemm}{Lemma}
\newtheorem{definition}{Definition}
\newtheorem{remark}{Remark}

\newtheorem{example}{Example}
\newtheorem{question}{Question}

\def\bproof{\noindent{\bf Proof.\;}}
\def\eproof{\hfill$\square$\medskip}
\def\dist{{\rm dist}\,}
\def\sign{{\rm sign}\,}
\def\deg{{\rm deg}\,}

\newcounter{marnote}

\begin{document}
\title{A fully nonlinear version of the Yamabe problem on locally conformally flat manifolds with umbilic boundary}
\author{YanYan Li \footnote{Department of Mathematics, Rutgers University, 110 Frelinghuysen Road, Piscataway, NJ 08901, USA}~~\footnote{Partially supported by NSF grant DMS-0654267 and a Rutgers University Research Council grant}\, and Luc Nguyen \footnote{OxPDE, Mathematical Institute, University of Oxford, 24-29 St Giles', Oxford, OX1 3LB, UK}}
\maketitle
\section{Introduction} 

The Yamabe problem asks, in dimension $n \ge 3$,
 to find on a given smooth closed Riemannian manifold
$(M^n, g)$
 a conformal metric which has constant scalar curvature.
This is the same as finding a solution of
\begin{equation}
-L_g u= \lambda_1(-L_g) u^{ \frac {n+2}{n-2} },
\ \ u>0,\qquad\text{on}\
M,
\label{W1-1}
\end{equation}
where
$$
L_g:= \Delta _g -\frac{ n-2}{ 4(n-1) } R_g
$$
is the conformal Laplacian of $g$, $\Delta_g$
is the Laplace-Beltrami operator, $R_g$ is the scalar curvature of $g$, and
$\lambda_1(-L_g)$ is the first eigenvalue of $-L_g$.
The answer was proved affirmative through the works of Yamabe himself \cite{Y}, Trudinger \cite{T}, Aubin \cite{Aubin} and Schoen \cite{S}. Different solutions to the Yamabe problem in the case $n$ $\leq$ $5$ and in the case $(M,g)$ is locally conformally flat were later given by Bahri and Brezis \cite{BB} and Bahri \cite{B}.

Let ${\cal M}(M,g)$ denote the set of smooth solutions of
\eqref{W1-1}.  It is not difficult to see that
 ${\cal M}(M,g)$ consists of one element if
$ \lambda_1(-L_g)<0$, and is equal to
$\{a \bar u\ |\ a > 0\}$ for some positive function
$\bar u$  on
$M$ if $ \lambda_1(-L_g)=0$.  When $ \lambda_1(-L_g)>0$, things are 
much more complex.  Schoen proved in the pioneering work \cite{Schoen-91} that
 if  $(M^n, g)$ is locally conformally flat with positive
$\lambda_1(-L_g)$, $n\ge 3$, then
\begin{equation}
\sup\ \Big\{\ 
\|u\|_{ C^m(M)} +\|u^{-1}\|_{ C^m(M)}\
\Big|\ u\in {\cal M}(M,g)\Big\}<\infty,\qquad
\forall\ m=1,2,3,\ldots
\label{V2-1}
\end{equation}
It was conjectured in \cite{Schoen-91} and \cite{Schoen-report},
 with a suggested 
strategy, that \eqref{V2-1} holds for general Riemannian manifolds
$(M^n, g)$, $n\ge 3$, and $\lambda_1(-L_g)>0$. 
Following the strategy, Li and Zhu proved in \cite{LiZhu}
the compactness result \eqref{V2-1} in dimension $n=3$.
Subsequent independent works were carried out by three groups as follows.
Druet proved in \cite{Druet03} that \eqref{V2-1}
holds for $n=4$ if
$\sup\ \{\
\|u\|_{ H^1(M) }\
|\ u\in {\cal M}(M,g)\}<\infty
$; and in  \cite{Druet04} that  \eqref{V2-1}
holds for $n=4, 5$.
Li and Zhang proved in \cite{LiZhang04} that
$\sup\ \{\
\|u\|_{ H^1(M) }\
|\ u\in {\cal M}(M,g)\}<\infty,
$ for $n=4$; in \cite{LiZhang05} that
 \eqref{V2-1}
holds for $4\le n\le 7$, or if $|W_g|+|\nabla W_g|>0$ on $M$,
or for $n=8, 9$ provided that the positive mass theorem holds in 
these dimensions;  and in
\cite{LiZhang07}
that  \eqref{V2-1}
holds for $n=10, 11$  provided that the positive mass theorem holds in
these dimensions.
Marques proved in \cite{Marques} that  \eqref{V2-1}
holds for $4\le n\le 7$, or if  $|W_g|>0$ on $M$.
Surprisingly, Brendle gave in 
\cite{brendle} 
 a counterexample to  \eqref{V2-1}
in dimension $n\ge 51$.  In subsequent papers, Khuri, Marques and Schoen in
\cite{KMS} extended \eqref{V2-1}
to  $12\le n \le 24$ provided  that the positive mass theorem holds in
these dimensions, while Brendle and Marques in
\cite{Brendle-Marques} improved the dimensions for counterexamples to 
$n\ge 25$.
In the above works to establish the compactness
property \eqref{V2-1}
in dimension $n\le 24$, the Liouville type theorem of Caffarelli, Gidas and 
Spruck  in \cite{CGS}
has played a very important role.
 
For compact manifolds with boundary, analogues of the Yamabe problem has been studied by many authors,
 see e.g. Cherrier \cite{C}, Escobar \cite{E1}, \cite{E2}, Han and Li \cite{H-L-Duke}, \cite{H-L-CAM},
Ambrosetti, Li and Malchiodi
 \cite{ALM}, Brendle \cite{Brendle1}, \cite{Brendle2},
Felli and Ould Ahmedou
 \cite{FO1},
Marques  \cite{M1},  \cite{M2}, Almaraz  \cite{A}, Brendle and Chen \cite{Brendle-Chen}
 and the references therein.
Since the works of Viaclovsky \cite{V}, \cite{V-CAG} and of Chang,
Gursky and Yang \cite{C-G-Y-JAM}, \cite{C-G-Y-AnnM}, there has been
much activity on fully nonlinear versions of the Yamabe problem (see
e.g. \cite{C-Y}, \cite{Li-ICM}, \cite{Tru-ICM}, \cite{V-SSC}, and
the references therein) and of the Yamabe problem on manifolds with
boundary (\cite{L-L-JEMS}, \cite{Chen}, \cite{J-L-L}, \cite{J},
\cite{ChenArxiv}); see also \cite{Sch} and \cite{G} for some other
boundary conditions.
In particular, a very general fully nonlinear version was proposed, and solved when the manifold is closed and locally conformally flat, by Li and Li \cite{L-L-CPAM}, \cite{L-L-Acta} (see also \cite{G-W}). Along this line, we consider analogues for manifolds with boundary.

Let $\Gamma$ be an open cone in $\RR^n$ and $f$ be a function defined on $\Gamma$ such that 
\begin{align}
&\Gamma \subset \RR^n \text{ is an open convex symmetric cone with vertex at the origin},\label{ConvexCone}\\
&\Gamma_n := \{\lambda \in \RR^n| \lambda_i > 0\} \subset \Gamma \subset \Gamma_1 := \{\lambda \in \RR^n| \lambda_1 + \ldots + \lambda_n > 0\}.\label{PositiveCone}\\
&f \in C^\infty(\Gamma) \cap C^0(\bar\Gamma) \text{  is
a non-negative symmetric function of $\lambda$},\label{f-smoothness}\\
&f > 0, \frac{\partial f}{\partial \lambda_k} > 0 \text{ in } \Gamma
\ (1 \leq k \leq n), \text{ and } f|_{\partial\Gamma} =
0,\label{f-ellipticity}\\
&\sum_{k=1}^n \frac{\partial f}{\partial
\lambda_k} \ge \delta\ \text{in}\ \Gamma \ \text{for some constant}\
\delta>0. \label{delta}
\end{align}
In some cases, we will also assume that
\begin{equation}
f \text{ is homogeneous of degree one on } \Gamma,\label{f-homogeneity}
\end{equation}
and/or
\begin{equation}
f \text{ is concave in } \Gamma. \label{f-concavity}
\end{equation}
Note that if $(f, \Gamma)$ satisfies
(\ref{ConvexCone})-(\ref{f-ellipticity}), \eqref{f-homogeneity} and \eqref{f-concavity}, then (\ref{delta}) is automatically satisfied; see
\cite{Urbas}. Important examples of such $(f,\Gamma)$ are
$(\sigma_k^{\frac{1}{k}},\Gamma_k)$ where $\sigma_k$ is the $k$-th
elementary symmetric function, i.e. $\sigma_k(\lambda)$ $=$
$\sum_{i_1 < \ldots < i_k} \lambda_{i_1} \ldots \lambda_{i_k}$, and
$\Gamma_k$ $=$ $\{\lambda \in \RR^n: \sigma_l(\lambda) > 0, 1 \leq l
\leq k\}$. It is well known that $(f, \Gamma)=(\sigma_k^{\frac{1}{k}},
\Gamma_k)$ satisfies (\ref{ConvexCone})-(\ref{f-concavity}) and
$\Gamma_k$ is the connected component of $\{\sigma_k>0\}$ containing $\Gamma_n\equiv \{\lambda_i>0\}$. See e.g. \cite{C-N-S-Acta}
where fully nonlinear elliptic equations with eigenvalues of the
Hessian $\nabla^2 u$ in such cones were first studied. In Section \ref{GammaType}, we give a construction for $f$ satisfying \eqref{f-smoothness}-\eqref{f-concavity} for a given cone $\Gamma$ that satisfies \eqref{ConvexCone}, \eqref{PositiveCone} and admits a smooth concave defining function, i.e. there exists some function $h$ $\in$ $C^\infty(\Gamma) \cap C^0(\bar\Gamma)$ satisfying $\nabla^2 h$ $\leq$ $0$ in $\Gamma$, $h$ $>$ $0$ in $\Gamma$ and $h$ $=$ $0$ on $\partial\Gamma$.

Let $(M^n,g)$ be a smooth compact
Riemannian manifold with boundary ($n$ $\geq$ $3$), $A_g$ denote the Schouten tensor of $M$, i.e.
\[
A_g = \frac{1}{n-2}\Big(Ric_g - \frac{R_g}{2(n-1)}g\Big),
\]
and $\lambda(A_g)$ denote its eigenvalues with respect to the metric $g$. Here $Ric_g$ denotes the Ricci curvature of $g$. Let $h_g$ denote the mean curvature of $\partial M$ with respect to the inner normal (so that the mean curvature of a Euclidean ball is positive). Let $N_1$, \ldots, $N_m$ be the components of $\partial M$ and $c_1$, \ldots $c_m$ be real numbers. Consider the problem
\begin{equation}
\left\{\begin{array}{l}
f(\lambda(A_{u^\frac{4}{n-2}g})) = 1 \text{ in }M^\circ,\\
\lambda(A_{u^\frac{4}{n-2}g}) \in \Gamma \text{ and } u > 0 \text{ in } M,\\
h_{u^{\frac{4}{n-2}}g} = c_k \text{ on } N_k, k = 1 \ldots m.
\end{array}\right.
\label{Main-Eqn}
\end{equation}
Note that under a conformal change of the metric, the Schouten tensor and the mean curvature change according to
\[
A_{u^{\frac{4}{n-2}}g}
	= -\frac{2}{n-2}u^{-1}\,\nabla_g^2 u + \frac{2n}{(n-2)^2}\,u^{-2}\nabla_g u \otimes \nabla_g u - \frac{2}{(n-2)^2}\,u^{-2}|\nabla_g u|^2\,I + A_g,
\]
and
\[
h_{u^{\frac{4}{n-2}}g}
	= u^{-\frac{n}{n-2}}\Big[\frac{\partial u}{\partial\nu} + \frac{n-2}{2}h_g\,u\Big].
\]
Here $\nu$ is the outer unit normal to $\partial M$.

On closed manifolds, the primitive of problem \eqref{Main-Eqn},
i.e. the first two lines in \eqref{Main-Eqn} under
the assumption that $\lambda(A_g)\in \Gamma$,  has
been examined extensively in the literature.
The problem for  $(f, \Gamma)=(\sigma_k^{\frac 1k}, \Gamma_k)$ was 
first proposed  by  Viaclovsky in \cite{V}. 
Chang, Gursky and Yang 
\cite{C-G-Y-JAM} proved, for 
 $(f, \Gamma)=(\sigma_2^{ \frac 12}, \Gamma_2)$ in dimension $n=4$,
the existence and compactness of solutions.
Li and Li  \cite{L-L-Acta} proved the existence and compactness of solutions
on locally conformally flat Riemannian manifolds
for $(f, \Gamma)$ satisfying (\ref{ConvexCone})-(\ref{f-ellipticity})
and \eqref{f-concavity}; see \cite{G-W},  \cite{L-L-CPAM} and 
\cite{G-V-W} for earlier results on
 $(f, \Gamma)=(\sigma_k^{\frac 1k}, \Gamma_k)$. These existence and compactness results have been extended as follows. 
Gursky and Viaclovsky  \cite{G-V-AnnM}
 proved the existence and compactness
of solutions for 
 $(f, \Gamma)=(\sigma_k^{\frac 1k}, \Gamma_k)$, $n<2k$;
see \cite{V-CAG} and 
\cite{G-V-AdvM} for some earlier results, as well as a later related paper
\cite{T-W}. Trudinger and Wang \cite{TW-new}  proved the existence and compactness of solutions for 
$(f, \Gamma)=(\sigma_k^{\frac 1k}, \Gamma_k)$
in dimension $n=2k$ for all
$k\ge 2$.
Sheng, Trudinger and Wang \cite{S-T-W}
proved the existence result for
$(f, \Gamma)=(\sigma_2^{ \frac 12}, \Gamma_2)$
in dimensions $n$ $\geq$ $3$; while
an independent proof was given by Ge and Wang \cite{Ge-Wang}
in dimensions $n> 8$.
The previously mentioned Liouville type theorem in \cite{CGS} was extended by Li and Li \cite{L-L-Acta}
to $(f, \Gamma)$ satisfying \eqref{ConvexCone}-\eqref{f-ellipticity}; 
see
\cite{C-G-Y-JAM},
\cite{CGY-2003} and \cite{L-L-CPAM} for some earlier results.

The case where $(f,\Gamma)$ $=$ $(\sigma_1,\Gamma_1)$ boils down to the Yamabe problem on manifolds with boundary in the ``positive case''. This particular case was studied in  \cite{C},  \cite{E1}, \cite{E2},  \cite{H-L-Duke}, \cite{H-L-CAM},  \cite{ALM} and \cite{Brendle2}. In these works, it is important to consider the first eigenvalue $\lambda_1$ of the eigenvalue problem
\begin{equation}
\left\{\begin{array}{l}
-\Delta_g\varphi + \frac{n-2}{4(n-1)}R_g\varphi
    = \lambda\varphi~ \text{ in } M^\circ,\\
\frac{\partial\varphi}{\partial\nu} + \frac{n-2}{2}\,h_g\varphi
    = 0~ \text{ on } \partial M.
\end{array}\right.
\label{7-1}
\end{equation}
In the
literature, $(M,g)$ is sometimes classified as of positive, negative
or zero type according to whether $\lambda_1$ is positive, negative
or zero. The signs of $\lambda_1$ is invariant under a conformal
change of metrics, i.e.
$\text{sign}(\lambda_1(M,g))=\text{sign}(\lambda_1(M, \xi^{ \frac
4{n-2}}g))$ for any positive smooth function $\xi$ on $M$.
We
note that in order for \eqref{Main-Eqn} to have a solution when
$(f,\Gamma)$ $=$ $(\sigma_1,\Gamma_1)$ and $c_k$ $=$ $0$, it is
necessary that $(M,g)$ is of positive type, i.e. $(M,g)$ admits a
conformal metric of positive scalar curvature and minimal boundary.
Conversely, it was proved, when $c_1=\cdots =c_m=c\in \RR$,
 in \cite{H-L-Duke} and \cite{H-L-CAM}
that the problem \eqref{Main-Eqn} is solvable for $(f,\Gamma)$ $=$
$(\sigma_1,\Gamma_1)$ under one of the following hypotheses:
\begin{enumerate}[(i)]
\item $n$ $\geq$ $3$, $(M,g)$ is locally conformally flat, $\partial M$ is umbilic, and $\lambda_1(M,g)$ $>$ $0$;
\item $n$ $\geq$ $5$, $\partial M$ is not umbilic,
and $\lambda_1(M,g)$ $>$ $0$.
\end{enumerate}
Note that
the case $c_1=\cdots=c_m=0$ was proved earlier
in \cite{E1} and \cite{E2}, where other cases were also studied. Recall that a hypersurface is umbilic if its second fundamental form is a multiple of the metric.

We introduce the following definition.

\begin{definition}\label{Type}
Let $(M,g)$ be a smooth compact Riemannian manifold with boundary
$\partial M$ and  $\Gamma$ satisfy \eqref{ConvexCone}
and \eqref{PositiveCone}.
\begin{enumerate}[(a)]
\item We say that $(M,g)$ is of $\Gamma$-positive type if there is a $C^2$ positive function $u$ on $M$ such that the conformal metric $\tilde g$ $=$ $u^{\frac{4}{n-2}}\,g$ satisfies $\lambda(A_{\tilde g})$ $\in$ $\Gamma$ in $M$ and $h_{\tilde g}$ $\geq$ $0$ on $\partial M$.
\item We say that $(M,g)$ is of $\Gamma$-nonpositive type if there is a $C^{0,1}$ positive function $u$ on $M$ such that the conformal metric $\tilde g$ $=$ $u^{\frac{4}{n-2}}\,g$ satisfies $\lambda(A_{\tilde g})$ $\in$ $\partial\Gamma$ in $M^\circ$ and $h_{\tilde g}$ $\leq$ $0$ on $\partial M$
in the viscosity sense 
(see Definition \ref{ViscoSol} in Section \ref{GammaType}). 
\end{enumerate}
\end{definition}

Clearly, a necessary condition for
the solvability of \eqref{Main-Eqn} for $c_k$ $=$ $0$ is that
$(M,g)$ is of $\Gamma$-positive type. 
Also, if $(M,g)$ is of $\Gamma$-positive type,
then $\lambda_1(M,g)$ $>$ $0$.
We note that the two types are mutually exclusive,
 see Lemmas \ref{MutuallyExclusive}
in Section \ref{GammaType}.
  Moreover, if
$(M,g)$ is locally conformally flat, $\partial M$ is umbilic, $\lambda_1(M,g)$ $>$
$0$, and $\Gamma$ satisfies  
\eqref{ConvexCone}-\eqref{PositiveCone} and admits a smooth defining concave function,
then $(M,g)$ must be of $\Gamma$-positive or
$\Gamma$-nonpositive type; see
Lemma \ref{Exhaustive} in Section
\ref{GammaType}.

\begin{example}\label{ExampleOfGammaType}
Let $M$ $=$ $\SS^n\setminus (B_1 \cup B_2)$ where $\SS^n$ is the
standard sphere, $B_1$ and $B_2$ are two disjoint non-touching
geodesic balls in $\SS^n$. Then $M$ is of $\Gamma_k$-positive type
for $1$ $\leq$ $k$ $<$ $\frac{n}{2}$, and of
$\Gamma_k$-nonpositive type for $\frac{n}{2}$ $\leq$ $k$ $\leq$ $n$. For a proof, see Section \ref{GammaType}.
\end{example}

Concerning the existence of solutions to \eqref{Main-Eqn}, Jin, Li and Li showed in \cite{J-L-L} that if $M$ has umbilic boundary and is locally conformally flat near its boundary, then \eqref{Main-Eqn} is always solvable when $\Gamma$ $\subset$ $\Gamma_j$ for some $j$ $>$ $\frac{n}{2}$, $c_k$ $\geq$ $0$, and $(M,g)$ is of $\Gamma$-positive type. In the same paper, they showed that the requirement $\Gamma$ $\subset$ $\Gamma_j$ for some $j$ $>$ $\frac{n}{2}$ can be relaxed if $c_1$ $=$ \ldots $=$ $c_m$ $=$ $0$. A similar statement was proved independently by Chen \cite{Chen}.

Concerning estimates for solutions of \eqref{Main-Eqn}, local first
and second derivative estimates are fairly well established. Under
the assumption that $0$ $<$ $u$ $<$ $b$ for some positive constant
$b$, local interior gradient estimates were established in
\cite{Li-CPAM} for $(f, \Gamma)$ satisfying
(\ref{ConvexCone})-(\ref{f-homogeneity}); while under an additional
concavity assumption of $f$ in $\Gamma$, different proofs were given
in \cite{Chen05}, \cite{Li-CPAM} and \cite{Wang}.  The estimates for
$(f, \Gamma)=(\sigma_k^{\frac 1k}, \Gamma_k)$ were proved earlier in
\cite{GW03}. Local boundary gradient estimates were established in
\cite{Li-CPAM} for $(f, \Gamma)$ satisfying
(\ref{ConvexCone})-(\ref{f-homogeneity}); while on locally conformally flat
manifolds with umbilic boundary, a different proof was given in
\cite{Chen} under the additional assumptions that $f$ is concave in
$\Gamma$, $c_k\ge 0$, and $\Gamma\subset \Gamma_2$ if $c_k>0$.

Under a much stronger assumption that $a$ $<$ $u$ $<$ $b$ for some
positive constants $a$ and $b$, local interior and boundary
estimates for $(f, \Gamma)$ satisfying
(\ref{ConvexCone})-(\ref{f-homogeneity}) were established in 
\cite{L-L-PrivComm} and
\cite{J-L-L} respectively. Local interior second derivative
estimates were proved in \cite{GW03} for $(f, \Gamma)=
(\sigma_k^{\frac 1k}, \Gamma_k)$ and were extended in \cite{L-L-CPAM} to $(f, \Gamma)$
satisfying (\ref{ConvexCone})-(\ref{f-ellipticity}), (\ref{f-homogeneity}) and
(\ref{f-concavity}). Local boundary second derivative estimates for
$(f, \Gamma)$ satisfying (\ref{ConvexCone})-(\ref{f-ellipticity}), (\ref{f-homogeneity})
and (\ref{f-concavity}) on locally conformally flat $M$ with umbilic
boundary $\partial M$ were established in
 \cite{J-L-L} for $c_k\ge 0$. Similar results were obtained independently by different methods in \cite{Chen}.
The locally conformally flat assumption is unnecessary --- see
\cite{J-L-L} for $(f, \Gamma)=(\sigma_n^{\frac 1n}, \Gamma_n)$ and
\cite{J} for general $(f, \Gamma)$.

The main purpose of this paper is to establish $C^0$ bounds for solutions of \eqref{Main-Eqn} for
the case where $M$ is a smooth locally conformally flat compact manifold with umbilic boundary.
 Such estimates enable us to establish, in view of the aforementioned first and second derivative estimates,
 the degree theory for fully nonlinear elliptic operators of second order in 
\cite{L-CPDE}
 and the degree counting formula in \cite{H-L-Duke},
 the existence of solutions of \eqref{Main-Eqn} for non-negative $c_k$'s.
\begin{theorem}\label{Main-Thm}
Let $(M^n,g)$, $n$ $\geq$ $3$, be a smooth compact locally
conformally flat Riemannian manifold with umbilic boundary $\partial
M$ and $N_1$, \ldots, $N_m$ be the components of $\partial M$.
Assume that (M,g) is not conformally equivalent to the standard
half-sphere $\SS^n_+$. Let $(f,\Gamma)$ satisfy
\eqref{ConvexCone}-\eqref{f-homogeneity} and assume that $(M,g)$ is of
$\Gamma$-positive type. For any given $\beta$ $>$ $0$ there exists a
constant $C$ $=$ $C(n,(M,g),(f,\Gamma),\beta)$ such that if $u$ $\in$
$C^2(M)$ is a positive solution to \eqref{Main-Eqn}
for some constants $c_1$, \ldots, $c_m$ satisfying $|c_k|$ $\leq$
$\beta$, then
\begin{equation}
\|u\|_{C^0(M,g)} + \|u^{-1}\|_{C^0(M,g)} \leq C.
\label{Main-Est-Thm}
\end{equation}
\end{theorem}

\begin{remark}
We note that in Theorem \ref{Main-Thm}, we do not assume that $f$ be
 concave. Also, regarding the upper bound of $u$, it suffices to 
assume that $\lambda_1(M,g)$ $>$ $0$ instead of $(M,g)$ is of $\Gamma$-positive type (see Proposition \ref{UpperBound} in Section 3).
\end{remark}

\begin{theorem}\label{Existence}
Let $(M^n,g)$, $n$ $\geq$ $3$, be a smooth compact locally conformally flat Riemannian manifold with umbilic boundary $\partial M$ and $N_1$, \ldots, $N_m$ be the components of $\partial M$. Let $(f,\Gamma)$ satisfy 
\eqref{ConvexCone}-\eqref{f-ellipticity}
 and assume that $(M,g)$ is of $\Gamma$-positive type. Assume in addition condition \eqref{f-concavity}, i.e. $f$ is concave. Then, for any collection of non-negative numbers $c_k$ $\geq$ $0$, there exists a positive solution $u$ $\in$ $C^\infty(M)$ of \eqref{Main-Eqn}.

Moreover if $(M,g)$ is not conformally equivalent to the standard half-sphere, then any such $u$ satisfies, for any $l$ $=$ $0$, $1$, \ldots, 
\begin{equation}
\|u\|_{C^l(M)} + \|u^{-1}\|_{C^l(M)} \leq C,
\label{C4-Est}
\end{equation}
where $C$ depends on $n$, $(M,g)$, $(f,\Gamma)$, $\{c_k\}$ and $l$.
\end{theorem}

It should be noted that if one removes the non-negativity assumption on $\{c_k\}$ in Theorem \ref{Existence}, estimate \eqref{C4-Est} may fail. In fact, for any $\epsilon$ $>$ $0$ and any $2$ $\leq$ $k$ $\leq$ $n$, there exist two non-touching geodesic balls of the standard sphere $\SS^n$, denoted by $B_1$ and $B_2$, and a sequence of smooth positive solutions $\{u_l\}$ of \eqref{Main-Eqn} with $M$ $=$ $\SS^n \setminus (B_1 \cup B_2)$, $N_1$ $=$ $\partial B_1$, $N_2$ $=$ $\partial B_2$, $(f,\Gamma)$ $=$ $(\sigma_k^{\frac{1}{k}},\Gamma_k)$, $c_1$ $=$ $0$, $-\epsilon$ $<$ $c_2(l)$ $<$ $0$ such that $\|u_l\|_{C^2(M)}$ $\rightarrow$ $\infty$. See Lemma \ref{NegEx-Global} for details. As pointed out earlier, when $2$ $\leq$ $k$ $<$ $\frac{n}{2}$, $M$ is of $\Gamma_k$-positive type and therefore satisfies the hypotheses of Theorem \ref{Existence}. We thus ask the following question.

\begin{question}\label{question-negativeC}
Does the existence part of Theorem \ref{Existence} hold without the non-negativity assumption on $\{c_k\}$?
\end{question}
Another relevant question is:
\begin{question}\label{question-umbilic}
Does the existence part of Theorem \ref{Existence} hold for general boundary $\partial M$?
\end{question}

In the special case where $(f,\Gamma)$ $=$ $(\sigma_1,\Gamma_1)$, Theorems \ref{Main-Thm} and \ref{Existence} were proved by Han and Li in \cite{H-L-Duke}. Our proof of Theorem \ref{Main-Thm} is very different.
A proof along the line of \cite{H-L-Duke} would require more development in the analysis of fully nonlinear conformally invariant equations.

As mentioned earlier, 
an analogue of Theorem \ref{Main-Thm} for closed locally conformally flat Riemannian manifolds was established by Li and Li \cite{L-L-Acta}.
An important ingredient in their approach was the use of the positive mass theorem of Schoen and Yau \cite{S-Y} for locally conformally flat manifolds to reduce their analysis to that on a Euclidean domain. This approach is also useful in our setting. As a tool for the passage from locally conformally flat manifolds with umbilic boundary to Euclidean domains, we establish, based on the positive mass theorem,
\begin{theorem}\label{Structure}
Let $(M^n,g)$, $n$ $\geq$ $3$, be a smooth compact locally conformally flat Riemannian manifold with (non-empty) umbilic boundary $\partial M$. Assume in addition that $\lambda_1(M,g)$ $>$ $0$. Then there exist a non-empty family of non-overlapping geodesic open balls $\{B_\alpha\}$ in the standard sphere $(\SS^n, g_{\SS^n})$ and a closed subset $\Lambda$ of $\SS^n$ of Hausdorff dimension at most $\frac{n-2}{2}$ such that the following conclusions hold.
\begin{enumerate}[(i)]
\item There exists a smooth conformal covering map $\Psi:$ $G$ $:=$ $\SS^n \setminus (\cup B_\alpha \cup \Lambda)$ $\rightarrow$ $(M,g)$, where $G$  is equipped with the metric inherited from the round $\SS^n$.
\item If $\bar B_\alpha \cap \bar B_\beta$ $=$ $\{p\}$, then $p$ $\in$ $\Lambda$.
\item If $\{B_{\alpha_j}\}$ $\subset$ $\{B_\alpha\}$ is a sequence of distinct balls ``converging'' to a point $p$, i.e. their centers converge to $p$ in $\SS^n$ and their radii tend to $0$, then $p$ $\in$ $\Lambda$.
\item If we write the pull-back metric of $g$ to $G$ by $\Psi$ as $w^{\frac{4}{n-2}}\,g_{\SS^n}$, then
\[
w(p) \rightarrow \infty \text{ as } \dist_{\SS^n}(p,\Lambda) \rightarrow 0.
\]
\end{enumerate}
\end{theorem}

By virtue of Theorem \ref{Structure}, the analysis of \eqref{Main-Eqn} can be ``reduced'' to that of an analogue on a Euclidean domain with a possibly singular boundary. Let $B$ $=$ $B_r(x)$ $\subset$ $\RR^n$ be a ball and $\{B_\alpha = B_{r_\alpha}(x_\alpha)\}$ be a family of (possibly empty, at most countably many) open and mutually non-overlapping balls contained in
 $B$. Let $\Lambda$ be a (possibly empty) closed subset of $\bar B\setminus (\cup B_\alpha)$ satisfying
\begin{enumerate}[(i)]
\item $\big(\cup B_\alpha\big) \cup \Lambda \neq \emptyset$,
\item if $\bar B_\alpha \cap \bar B_\beta$ $=$ $\{p\}$,
or $\bar B_\alpha\cap \partial B=\{p\}$,
then $p$ $\in$ $\Lambda$,
\item if $\{B_{\alpha_j}\}$ $\subset$ $\{B_\alpha\}$ is a sequence of distinct balls ``converging'' to a point $p$ in $\bar B$ in the sense that $x_{\alpha_j}$ $\rightarrow$ $p$ and $r_{\alpha_j}$ $\rightarrow$ $0$, then $p$ $\in$ $\Lambda$.
\end{enumerate}
Note that by these assumptions, the set
\[
\Omega := B \setminus (\cup \bar B_\alpha \cup \Lambda)
\]
is and open set. Consider on $\Omega$ the following equation
\begin{equation}
\left\{\begin{array}{l}
f(\lambda(A_{u^{\frac{4}{n-2}}g_{flat}})) = 1 \text{ in } \Omega,\\
\lambda(A_{u^{\frac{4}{n-2}}g_{flat}}) \in \Gamma \text{ and } u > 0 \text{ in } \bar\Omega\setminus\Lambda,\\
u(x) \rightarrow \infty \text{ as } x \rightarrow \Lambda,\\
\frac{\partial u}{\partial\nu} + \frac{n-2}{2r}\,u = c(B)\,u^{\frac{n}{n-2}} \text{ on }\partial B \setminus \Lambda,\\
\frac{\partial u}{\partial\nu} - \frac{n-2}{2r_\alpha}\,u = c(B_\alpha)\,u^{\frac{n}{n-2}} \text{ on }\partial B_\alpha \setminus \Lambda.
\end{array}\right.
\label{Main-Eqn-Sing-Gen}
\end{equation}
Here $g_{flat}$ is the Euclidean flat metric and $\nu$ is the outer unit normal.
\begin{theorem}\label{By-product}
Let $(\Omega,\Lambda)$ be as above and $(f,\Gamma)$ satisfy \eqref{ConvexCone}-\eqref{delta}.
For any given $\beta$ $>$ $0$ and $\epsilon$ $>$ $0$ there exists a constant $C$ $=$ $C(n,f,\Gamma,(\Omega,\Lambda),\beta,\epsilon)$ such that if $u$ $\in$ $C^2(\bar \Omega\setminus\Lambda)$ satisfies \eqref{Main-Eqn-Sing-Gen} for some constants $c(B)$, $c(B_\alpha)$ $>$ $-\beta$, then
\[
\sup\{u(x): x \in \bar\Omega, \dist(x,\Lambda) \geq \epsilon\} \leq C.
\]
Here when $\Lambda$ $=$ $\emptyset$, the set $\{x \in \bar\Omega: \dist(x,\Lambda) \geq \epsilon\}$ is understood as $\bar\Omega$.
\end{theorem}
It is readily seen that Theorems \ref{Structure} and \ref{By-product} imply the ``upper bound'' in Theorem \ref{Main-Thm}.

\begin{remark}
(a) In fact, in the case where the family $\{B_\alpha\}$ is empty or finite, the estimate in Theorem \ref{By-product} takes a better form:
\[
\sup\{u(x): x \in \bar\Omega, \dist(x,\Lambda) \geq \epsilon\} \leq C(n,f,\Gamma,(\Omega,\Lambda),\beta)\,\epsilon^{-\frac{n-2}{2}}.
\]
\noindent(b) As far as interior estimate is concerned in Theorem \ref{By-product}, the assumption \eqref{delta} can be dropped. The conclusion one gets is
\[
u(x) \leq C(n,f,\Gamma,\beta)\,\dist(x,\partial\Omega)^{-\frac{n-2}{2}} \text{ in } \Omega.
\]
Note that the constant above does not depend on $(\Omega,\Lambda)$. Also, $\partial\Omega$ contains $\Lambda$ by definition. See Proposition \ref{DistanceEst-Sing} in Section \ref{C0Est}.
\end{remark}

One might have noticed that between Theorem \ref{Main-Thm} and Theorem \ref{Existence}, there is some additional restriction on the sign of the constants $c_k$'s. This is due to the limit of known local $C^2$ estimates which were used in the proof of Theorem \ref{Existence}. Since this is of special relevance to our result, we state it as a reference.

\begin{thm}[\cite{J-L-L}]\label{C2BdryEst}
Let $(f,\Gamma)$ satisfy \eqref{ConvexCone}-\eqref{f-ellipticity}, \eqref{f-homogeneity}
and \eqref{f-concavity}. 
Let $(M,g)$ be a smooth locally conformally flat Riemannian manifold with umbilic boundary, $\mathcal{O}$ be an open subset of $M$, and $\eta$ and $\psi$ $>$ $0$ be smooth functions satisfying one of the following conditions:
\begin{enumerate}[(i)]
\item $\eta$ $\equiv$ $0$ and $\frac{\partial \psi}{\partial \nu}$ $\equiv$ $0$ on $(\mathcal{O} \cap \partial M) \times \RR$,
\item or $\eta$ $>$ $0$ on $(\mathcal{O} \cap \partial M) \times \RR$.
\end{enumerate}
Then for any $\mathcal{O}_1$ $\Subset$ $\mathcal{O}$ and $C_0$ $>$ $0$, there exists a constant $C$ depends only on $(M,g)$, $(f,\Gamma)$, $(\eta,\psi)$, $\mathcal{O}_1$, $\mathcal{O}$ and $C_0$ such that for any $u$ $\in$ $C^2(\mathcal{O})$ satisfying
\begin{equation}
\left\{\begin{array}{l}
f(\lambda(A_{u^\frac{4}{n-2}g})) = \psi(x,u) \text{ in } \mathcal{O},\\
\lambda(A_{u^\frac{4}{n-2}g}) \in \Gamma \text{ and } u > 0 \text{ in } \mathcal{O},\\
\frac{\partial u}{\partial\nu} + \frac{n-2}{2}h_g\,u = \eta(x,u) \text{ on } \mathcal{O} \cap \partial M,
\end{array}\right.
\label{C2Eqn}
\end{equation}
and
\[
|u| + |u^{-1}| + |\nabla u| \leq C_0 \text{ in } \mathcal{O},
\]
there holds
\[
|\nabla^2 u| \leq C \text{ in } \mathcal{O}_1.
\]
\end{thm}
Note that in \cite{J-L-L}, $u$ $\in$ $C^4(\mathcal{O})$ was assumed. On the other hand, by standard elliptic theories, a $C^2(\mathcal{O})$ solution of \eqref{C2Eqn} is in $C^\infty(\mathcal{O})$; see e.g. \cite[Lemma 17.16]{GT}.

As mentioned earlier, similar results were proved independently by a different method
in \cite{Chen},  and  the locally conformally flat
assumption on $M$ is unnecessary (\cite{J-L-L}, \cite{J}).  See
also \cite{ChenArxiv} where related equations were treated.
Second derivative estimates for fully nonlinear
elliptic equations with Neumann type boundary
conditions can be found in \cite{LT} and the references in.

Naturally one asks if the estimate in Theorem \ref{C2BdryEst} holds for arbitrary choices of smooth $\eta$ and $\psi$ $>$ $0$. In particular, does it hold for $\eta$ $\equiv$ $c$ and $\psi$ $\equiv$ $1$ where $c$ is a negative constant? Evidently, such estimate holds for $(f,\Gamma)$ $=$ $(\sigma_1,\Gamma_1)$. It turns out that the answer is negative, in general, for any $(f,\Gamma)$ $=$ $(\sigma_k^{\frac{1}{k}},\Gamma_k)$ with $2$ $\leq$ $k$ $\leq$ $n$ (see Lemma \ref{NegEx} in Section \ref{Counterex}).

Concerning Question \ref{question-negativeC}, simple examples of manifolds that give a positive answer are ones that are conformally covered by the standard half sphere.
 
\begin{proposition}\label{PositiveView-High-k} 
\begin{enumerate}[(a)] 
\item Assume that $(M,g)$ is conformally covered by the standard half sphere $\SS^n_+$. Then for any $c$ $\in$ $\RR$, there exists a metric $\hat g$ which is conformal to $g$, has constant sectional curvature $1$, and has constant boundary mean curvature $c$.

\item Let $M$ be a smooth compact locally conformally flat Riemannian manifold with umbilic boundary. The following statements are equivalent.
\begin{enumerate}[(i)]
\item $M$ is conformally covered by the standard half sphere.
\item $M$ is of $\Gamma_n$-positive type.
\item $M$ is of $\tilde\Gamma$-positive type for some $\tilde \Gamma$ $\subset$ $\Gamma_k$, $\frac{n}{2}$ $\leq$ $k$ $\leq$ $n$. 
\end{enumerate} 
\end{enumerate} 
\end{proposition} 

To shed more light on the failure of $C^2$ estimates and the existence of solutions for the $\sigma_k$ equation with $2$ $\leq$ $k$ $<$ $\frac{n}{2}$, we state here a result for annuli which we will present in another paper. Consider
\begin{equation} 
\left\{\begin{array}{l}
\sigma_k(\lambda(A_{u^{\frac{4}{n-2}}g_{flat}})) = 1 \text{ in } B_{R} \setminus B_1,\\
\lambda(A_{u^{\frac{4}{n-2}}g_{flat}}) \in \Gamma_k \text{ and } u > 0 \text{ in } \bar B_{R} \setminus B_1,\\
\frac{\partial u}{\partial r} + \frac{n-2}{2}\,u = -c_1\frac{n-2}{2}\,u^{\frac{n}{n-2}} \text{ on } \partial B_1,\\
\frac{\partial u}{\partial r} + \frac{n-2}{2R}\,u = c_2\frac{n-2}{2R}\,u^{\frac{n}{n-2}} \text{ on } \partial B_R,
\end{array}\right. 
\label{CCCircEqn}
\end{equation} 
where $c_1$ and $c_2$ are given constants.
 
\begin{theorem}[\cite{LiNg}]\label{PositiveView-Low-k} 
Assume that $1$ $\leq$ $k$ $<$ $\frac{n}{2}$. 
\begin{enumerate}[(a)] 
\item If $c_1 + c_2$ $\geq$ $0$, then \eqref{CCCircEqn} has a radial solution for any $R$ $>$ $1$. 
 
\item If $c_1 + c_2$ $<$ $0$, there exists $R_*$ $=$ $R_*(c_1,c_2,n,k)$ $>$ $1$ such that \eqref{CCCircEqn} has a radial solution for $R$ $\geq$ $R_*$. Moreover, if $k$ $\geq$ $2$, $R_*$ can be chosen so that \eqref{CCCircEqn} has no radial solution for $1$ $<$ $R$ $<$ $R_*$. 
 
\item If $c_1$ $=$ $c_2$ $=$ $0$, \eqref{CCCircEqn} always has a radial solution, which is a constant multiple of the cylindrical metric on the punctured space. If, in addition, $R$ $>$ $\exp\frac{\pi}{\sqrt{n-2k}}$, then \eqref{CCCircEqn} has an additional radial solution. 
\end{enumerate} 
\end{theorem}

In connection with the failure of boundary $C^2$ estimates for negative $c_k$'s, it is natural to ask the following question.

\begin{question}\label{Quest-Regularity}
Even though local boundary $C^2$ estimate may fail, what can one say about $C^{1,\alpha}$ estimates for conformally invariant elliptic
equations? More generally, if $w$ $>$ $0$ satisfies a fully nonlinear
 elliptic partial differential equation of the form
\[
\left\{\begin{array}{l}
f(\lambda(\nabla^2 w +b\,\frac{|\nabla w|^2}{w}I)) = \psi(x,w) \text{ in } B_1^+(0),\\
\lambda(\nabla^2 w +b\,\frac{|\nabla w|^2}{w}I) \in \Gamma \text{ in } \bar B_1^+(0),\\
w_n = \eta(x,w) \text{ on } B_1(0) \cap \{x_n = 0\},\\
w+w^{-1} + |\nabla w| \leq C_1 \text{ in } B_1^+(0),
\end{array}\right.
\]
for some constant $b$ and
some positive smooth function $\psi$, can one say
\[
[\nabla w]_{0,\alpha; B_{1/2}^+(0)} \leq C(n,f,\Gamma, \psi, \eta, b,C_1)
\]
for some $0$ $<$ $\alpha$ $<$ $1$? (To see the connection between the above equation and \eqref{C2Eqn},
simply put $w^{-2}\,g_{flat}$ $=$ $u^{\frac{4}{n-2}}\,g$, then
$b =-1/2$ and $\eta(x,w)\equiv $constant correspond
to the conformally invariant one.)
\end{question}

An important property of \eqref{Main-Eqn} is its invariance under conformal transformations. More precisely, if $\psi:$ $(\tilde M,\tilde g)$ $\rightarrow$ $(M,g)$ is a bijective conformal transformation and $u$ solves \eqref{Main-Eqn}, then $\tilde u$ $=$ $|Jac_\psi|^{\frac{n-2}{2n}}\,u \circ \psi$ satisfies
\[
\left\{\begin{array}{l}
f(\lambda(A_{\tilde u^\frac{4}{n-2}\tilde g})) = 1 \text{ in }\tilde M,\\
\lambda(A_{\tilde u^\frac{4}{n-2}\tilde g}) \in \Gamma \text{ and } \tilde u > 0 \text{ in } \tilde M,\\
\frac{\partial \tilde u}{\partial\nu} + \frac{n-2}{2}h_{\tilde g}\,\tilde u = c_k\,\tilde u^{\frac{n}{n-2}} \text{ on } \psi^{-1}(N_k), k = 1 \ldots m.
\end{array}\right.
\]
In \cite{L-L-CPAM}, it was shown that any conformally invariant differential operator of second order on a Euclidean domain must be of the form $f(\lambda(A_{u^{\frac{4}{n-2}}g_{flat}}))$ where $g_{flat}$ is the flat metric on $\RR^n$
for some symmetric function $f$. It turns out that a similar classification for conformally invariant boundary differential operators can be carried out.

Since the conformally invariant property involves not only the geometry of the boundary but also that of the ambient space, we model a boundary operator on Euclidean domains by
\begin{align*}
B: \RR^n \times \RR \times \RR^n \times \SS^n \times \textrm{Sym}^{n\times n}
    &\rightarrow \RR\\
(x,u,\nabla u, \nu, H)
    &\mapsto B(x,u,\nabla u, \nu, H).
\end{align*}
Here $\nu$ and $H$ play the roles of the outer unit normal and the Weingarten map of the boundary respectively, and $\textrm{Sym}^{n \times n}$ denotes the set of symmetric $n \times n$ matrices.

\begin{definition}
We say that $B$ is conformally invariant if for any M\"{o}bius transformation $\psi$ of $\RR^n$, $u$ $\in$ $C^1(\RR^n)$, $\nu$ $\in$ $C^0(\RR^n;\SS^n)$, and $H$ $\in$ $C^0(\RR^n;\textrm{Sym}^{n\times n})$, there holds
\begin{equation}
B(\cdot , u_\psi, \nabla u_\psi, \nu, H) = B(\psi, u \circ \psi, (\nabla u) \circ \psi, \nu_\psi, H_\psi),
\label{BdryDiffOp-ConfInv}
\end{equation}
where
\begin{align*}
u_\psi(x)
    &= |Jac_\psi|^{\frac{n-2}{2n}}(x)\,u \circ \psi(x),\\
\nu_\psi(x)
    &= \frac{(d\psi(x))_*(\nu(x))}{|(d\psi(x))_*(\nu(x))|},\\
H_\psi(x)
    &= |Jac_\psi(x)|^{-\frac{1}{n}}(H + \frac{1}{n}\frac{\nabla|Jac_\psi(x)| \cdot \nu(x)}{|Jac_\psi(x)|}I).
\end{align*}
Note that if $H$ is the Weingarten map of a submanifold $N^{n-1}$ $\subset$ $\RR^n$ with unit normal vector $\nu$, then $H_\psi$ is the Weingarten map of $\psi(N)$.
\end{definition}

\begin{theorem}\label{BdryOpClass}
If $B$ is a conformally invariant boundary operator then
\[
B(x,s,p,\nu,H) = B\Big(0,1,0,e,s^{-\frac{2}{n-2}}\big(H + \frac{2}{n-2}\frac{p \cdot \nu}{s}I\big)\Big),
\]
where $e$ is any arbitrary unit vector. In particular, prescribing a conformally invariant ``boundary condition'' on an umbilic boundary is equivalent to prescribing its mean curvature.
\end{theorem}

The rest of the paper is organized as follows. We start with the proof of Theorem \ref{Structure} in Section \ref{Struct}. Section \ref{C0Est} is devoted to the proofs of Theorems \ref{Main-Thm} and
 \ref{By-product}. The proof of Theorem \ref{Existence} and of Proposition \ref{PositiveView-High-k} are carried out in Section \ref{ExistProof}. In Section \ref{Counterex}, we present counterexamples to local and global $C^2$ estimates discussed above. Some simple relevance results to the notion of $\Gamma$-type are done in Section \ref{GammaType}. The classification result in Theorem \ref{BdryOpClass} is done in Appendix \ref{BdryOpClassProof}. Finally, in Appendix \ref{App-Degree} we define a degree theory for second order nonlinear elliptic equations with nonlinear oblique boundary conditions.

\section{The holonomy covering of a locally conformally flat manifold with umbilic boundary}\label{Struct}

In this section, we prove Theorem \ref{Structure}. The idea is to attach another copy of $M$ to $M$ along its boundary and to use the corresponding result of Schoen and Yau \cite{S-Y} on closed manifolds.

\medskip
\noindent{\bf Proof of Theorem \ref{Structure}.} Since $\lambda_1(M,g)$ $>$ $0$, we can assume that $R_g$ $>$ $0$ in $M$ and $h_g$ $=$ $0$ on $\partial M$. This can
be achieved by working with $\hat g=\varphi^{\frac 4{n-2}}
g$ for a positive eigenfunction $\varphi$ of (\ref{7-1}).
 Let $M_2$ be the double of $M$ obtained by attaching a second copy of $M$ to $M$ along its boundary. Extend the metric $g$ to $M_2$ by an even extension. Then, as $\partial M$ is umbilic and minimal, $(M_2,g)$ is a closed $C^{2,1}$-regular locally conformally flat manifold whose scalar curvature is a positive function on $M_2$.

Let $\tilde M_2$ be the universal covering of $M_2$ with covering map $\pi$. Equip $\tilde M_2$ with the metric $\tilde g$ inherited from $g$. By a
deep result of Schoen and Yau (see \cite[Theorems 4.5, 4.7]{S-Y}), there exists an injective conformal map $\Phi:$ $\tilde M_2$ $\rightarrow$ $\SS^n$ such that $\partial\Phi(\tilde M_2)$ is the same as $\SS^n \setminus \Phi(\tilde M_2)$ and has Hausdorff dimension at most $\frac{n-2}{2}$.

Let $\hat M$ be any connected component of $\pi^{-1}(M)$ $\subset$ $\tilde M_2$. Then $\pi:$ $\hat M$ $\rightarrow$ $M$ is a covering map. In particular, $(\hat M,\tilde g)$ is a complete locally conformally flat manifold with umbilic boundary. We will show that $\Phi(\hat M)$ is of the form $\SS^n \setminus (\cup B_\alpha \cup \Lambda)$ where $B_\alpha$ and $\Lambda$ are as described in the statement of the theorem.

We claim that there is an identification of each component $\hat N$ of $\partial \hat M$ with a geodesic $(n-1)$-sphere $\FS_{\hat N}$ of $\SS^n$ such that the following holds:
\begin{enumerate}[(a)]
\item $\Phi(\hat N)$ is contained and dense in $\FS_{\hat N}$, $\Phi(\hat N)$ is open in the relative topology of $\FS_{\hat N}$, and $\FS_{\hat N} \setminus \Phi(\hat N)$ has Hausdorff dimension at most $\frac{n-2}{2}$,
\item For any two different components $\hat N_1$ and $\hat N_2$ of $\partial\hat M$, $\FS_{\hat N_1} \cap \FS_{\hat N_2}$ cannot contain more than one point,
\item $\Phi(\hat M\setminus
\hat N)$ belongs to only one component of $\SS^n\setminus\FS_{\hat N}$.
\end{enumerate}

First, if $\hat N$ is a component of the boundary of $\hat M$, then, by umbilicity, $\Phi(\hat N)$ is contained in some geodesic $(n-1)$-sphere of $\SS^n$, which will be denoted by $\FS_{\hat N}$. Clearly, $\Phi(\hat N)$ is open in $\FS_{\hat N}$. Let $Z_{\hat N}$ denote the boundary of $\Phi(\hat N)$ in $\FS_{\hat N}$ and consider a point $p$ in $Z_{\hat N}$. Then there exists $\hat x_k$ $\in$ $\hat N$ such that $\Phi(\hat x_k)$ $\rightarrow$ $p$ in the topology of $\FS_{\hat N}$, and so in the topology of $\SS^n$. If $p$ $=$ $\Phi(\tilde x)$ for some $\tilde x$ $\in$ $\tilde M_2$, the injectivity and local invertability of $\Phi$ implies that $\hat x_k$ $\rightarrow$ $\tilde x$, and so $\tilde x$ $\in$ $\hat N$ and so $p$ $\in$ $\Phi(\hat N)$, which contradicts the choice of $p$. We
conclude that the boundary
 $Z_{\hat N}$ must be a subset of $\SS^n\setminus \Phi(\tilde M_2)$, and thus must have Hausdorff dimension at most $\frac{n-2}{2}$ $<$ $n-2$.

We claim that
$\Phi(\hat N)$ is dense in
$\FS_{\hat N}$.
 Indeed, if not, then exist
$\epsilon>0$ and $p, q\in
\FS_{\hat N}$ such that $B_{2\epsilon}(p)
\subset \Phi(\hat N)$ and
$B_{2\epsilon}(q)\cap \Phi(\hat N)=\emptyset$, where $B_{2\epsilon}(p)$
and $B_{2\epsilon}(q)$ denote the geodesic balls
in $\FS_{\hat N}$.
We can then easily construct a Lipschitz map
$\psi: Z_{\hat N}\to \partial B_\epsilon(p)$
which is onto.  To see the existence
of such $\psi$, we may assume (modulo
a self-diffeomorphism of
$\FS_{\hat N}$) that $p$ and $q$ are antipodal
points of
$\FS_{\hat N}$.  Then along each geodesic connecting $p$ and $q$ there
must be a point in $Z_{\hat N}$.
Then the map projects $Z_{\hat N}$ along geodesics to $\partial
B_\epsilon(p)$ is clearly Lipschitz and onto.  The existence
of such a $\psi$ implies that the Hausdorff dimension
of $Z_{\hat N}$ is not smaller than that of
$\partial B_\epsilon(p)$ which is $n-2$, violating the fact that the
Hausdorff dimension of $Z_{\hat N}\le \frac {n-2}2$.
We have proved that $\Phi(\hat N)$ is dense in
$\FS_{\hat N}$.  It follows that
$Z_{\hat N}=\FS_{ \hat N}\setminus \Phi(\hat N)$.  We have established (a).

Next, consider two different components $\hat N_1$ and $\hat N_2$ of the boundary of $\hat M$ where $\FS_{\hat N_1} \cap \FS_{\hat N_2}$ is non-empty. Then either $\FS_{\hat N_1} \cap \FS_{\hat N_2}$ consists of a single point or it is an $(n-2)$-sphere. Pick any point $p$ in this intersection. Then we can pick $\hat x_k$ $\in$ $\hat N_1$ and $\hat y_k$ $\in$ $\hat N_2$ such that $\Phi(\hat x_k)$ and $\Phi(\hat y_k)$ converge to $p$ in the topology of $\SS^n$. Therefore, if $p$ $=$ $\Phi(\tilde x)$ for some $\tilde x$ $\in$ $\tilde M_2$, we can argue as in the previous paragraph to get $\hat x_k$, $\hat y_k$ $\rightarrow$ $\tilde x$, which implies $\tilde x$ $\in$ $\hat N_1 \cap \hat N_2$ contradicting our initial assumption that $\hat N_1$ and $\hat N_2$ are different connected components of $\partial\hat M$. We infer that $\FS_{\hat N_1} \cap \FS_{\hat N_2}$ is a subset of $\SS^n \setminus \Phi(\tilde M_2)$ and so has Hausdorff dimension at most $\frac{n-2}{2}$ $<$ $n-2$. This implies that $\FS_{\hat N_1} \cap \FS_{\hat N_2}$ consists of a single point.
We have established (b).

Consider a sphere $\FS_{\hat N}$ where $\hat N$ is some component of $\partial\hat M$. Then $\FS_{\hat N}$ separates $\SS^n$ into two components. We
claim that $\Phi(\hat M\setminus \hat N)$
 is contained in only one of those two components. Arguing indirectly, assume for some $\hat x$ and $\hat y$ in $\hat M\setminus \hat N$ that $\Phi(\hat x)$ and $\Phi(\hat y)$ belong to different components of $\SS^n \setminus \FS_{\hat N}$. Let $\gamma$ be a path in $\hat M\setminus \hat N$ that connects
 $\hat x$ to $\hat y$. Then $\Phi \circ \gamma$ is a path in
 $\Phi(\hat M\setminus \hat N)$ that connects $\Phi(\hat x)$ to $\Phi(\hat y)$.
Since $\gamma(t_0)$ is not in $\hat N$, there exists
an open neighborhood of $\gamma(t_0)$ which
has no intersection with $\hat N$, so
$\Phi \circ \gamma(t_0)$ is not in the closure of
$\Phi(\hat N)$, violating (a).
Part (c) is proved.

By the above claim, for each component $\hat N$ of $\partial \hat M$, there is a component of $\SS^n\setminus\FS_{\hat N}$, which we will denote by $B_{\hat N}$, that does not intersect $\Phi(\hat M)$. We would like to say that $\{B_\alpha\}$ $:=$ $\{B_{\hat N}\}$ and $\Lambda$ $:=$ $[\SS^n\setminus \Phi(\tilde M_2)] \cap [\SS^n\setminus(\cup B_{\hat N})]$ meet the requirements in our theorem.

As shown before, if two spheres $\FS_{\hat N_1}$ and $\FS_{\hat N_2}$ intersect non-trivially, the intersection consists of a single point. Hence, the balls $B_\alpha$ are non-overlapping. Also, by construction, $\Lambda$ is closed in $\SS^n$ and has Hausdorff dimension at most $\frac{n-2}{2}$. In addition, (ii) is evident according to our earlier consideration.

To show (iii), we assume that $\{B_{\hat N_j}\}$ converges to a point $p$ in $\SS^n$ where $\{\hat N_j\}$ is a sequence of  distinct components of $\partial\hat M$. Assume by contradiction that
$p$ does not belong to $\Lambda$,  clearly
 $p$ $=$ $\Phi(\tilde x)$ $\in$ $\Phi(\tilde M_2)$
for some $\tilde x\in \tilde M_2$. Then the injectivity and local invertability of $\Phi$ implies that for any $\epsilon$ $>$ $0$, all $\hat N_j$ must eventually lie in an $\epsilon$-ball centered at $\tilde x$. On the other hand, we can always pick some $\delta$ sufficiently small so that any $\delta$-neighborhood of $\pi(x)$ $\in$ $M_2$ cannot intersect more than one component of $\partial M$. Since $\pi:$ $\tilde M_2$ $\rightarrow$ $M_2$ is a local isometry, we get a contradiction for $\epsilon$ $<$ $\delta$.

Before proving (i), we claim that the set $G$ $:=$ $\SS^n \setminus (\cup B_\alpha \cup \Lambda)$ is connected.
The connectedness of
 $\SS^n\setminus \Lambda$
and $\partial B_\alpha \setminus \Lambda$ is a consequence of the
fact that the Hausdorff dimension of $\Lambda$ is less than
$n-2$, by an argument used earlier.  Next we show that $G$ is connected.
Pick $p$ and $q$ in $G$ and a path $\gamma:$ $[0,1]$ $\rightarrow$ $\SS^n\setminus \Lambda$ which connects $p$ to $q$.
 If $\gamma$ intersect infinitely many $B_\alpha$'s, then as these balls are non-overlapping, we can pick a subsequence of them that shrinks to a point on $\gamma$ which contradicts (iii) and the definition of $\gamma$. Hence $\gamma$ can only intersect at most finitely many $B_\alpha$. Therefore, it suffices to show that for any ball $B_\alpha$ that intersect $\gamma$, we can modify $\gamma$ to $\tilde \gamma$ that does not intersect $B_\alpha$. Define
\begin{align*}
\underbar{t}_\alpha
    &= \sup\{0 \leq t \leq 1: \gamma(s) \notin B_\alpha \text{ for } 0 \leq s \leq t\},\\
\bar{t}_\alpha
    &= \inf\{0 \leq t \leq 1: \gamma(s) \notin B_\alpha \text{ for } t \leq s \leq 1\}.
\end{align*}
Then these numbers are well-defined and $\underbar{t}_\alpha$ $<$ $\bar t_\alpha$. Moreover, we must have $\gamma(\underbar{t}_\alpha)$ and $\gamma(\bar t_\alpha)$ belongs to $\partial B_\alpha \setminus \Lambda$. We then obtain $\tilde \gamma$ from $\gamma$ by replacing $\gamma([\underbar{t}_1,\bar t_1])$ by a path in $\partial B_1 \setminus \Lambda$ with same endpoints. The claim follows.

To prove (i), we need to show that $G$ $=$ $\Phi(\hat M)$. As $B_\alpha$ does not intersect $\Phi(\hat M)$, we must have $\Phi(\hat M)$ $\subset$ $G$. Since $G$ is also connected, it suffices to show that $\Phi(\hat M)$ is both closed and open in $G$. To see that $\Phi(\hat M)$ is closed in $G$, pick $p_k$ $=$ $\Phi(\hat x_k)$ $\in$ $\Phi(\hat M)$ such that $p_k$ $\rightarrow$ $p$ $\in$ $G$. By
the  definition of $p$, $p$ $=$ $\Phi(\tilde x)$ for some $\tilde x$ $\in$ $\tilde M_2$. The injectivity and local invertability of $\Phi$ then implies that $\hat x_k$ $\rightarrow$ $\tilde x$, which shows $\tilde x$ $\in$ $\hat M$ and so $p$ $\in$ $\Phi(\hat M)$. To see that $\Phi(\hat M)$ is open in $G$, pick a point $p$ $=$ $\Phi(\hat x)$ $\in$ $\Phi(\hat M)$. If $\hat x$ $\in$ $\hat M^\circ$, then by the injectivity and local invertability of $\Phi$, $p$ is an interior point of $\Phi(\hat M)$ relative to $G$. Hence, it is enough to consider $\hat x$ $\in$ $\partial M$, i.e. $p$ $\in$ $\partial B_{\alpha_0}$ for some $\alpha_0$. By (iii), using the fact that $p$ is not in $\Lambda$, we can pick a ball $B_\delta(p)$ in $\SS^n$ such that $B_\delta(p) \cap \partial B_\alpha$ $=$ $\emptyset$ for all $\alpha$ $\neq$ $B_{\alpha_0}$, and $B_\delta(p) \setminus \partial B_{\alpha_0}$ has exactly two components, one is $B_\delta(p) \cap (G \setminus \partial B_{\alpha_0})$ and the other is $B_\delta(p) \cap B_{\alpha_0}$. By lowering $\delta$ if necessary, we can further assume that $\Phi|_{\Phi^{-1}(B_\delta(p))}:$ $\Phi^{-1}(B_\delta(p))$ $\rightarrow$ $B_\delta(p)$ is a bijection. Using the definition of $M_2$, $\hat M$ and $\tilde M_2$, we infer that $\Phi|_{\Phi^{-1}(B_\delta(p)) \cap \hat M}:$ $\Phi^{-1}(B_\delta(p)) \cap \hat M$ $\rightarrow$ $B_\delta(p) \cap G$ is also a bijection, which implies that $p$ is an interior point of $\Phi(\hat M)$ relative to $G$. Hence $\Phi(\hat M)$ is open (and closed) in $G$. Assertion (i) follows with $\Psi$ $=$ $\pi \circ \Phi^{-1}$.

Finally, we show (iv). Let $\tilde w^{\frac{4}{n-2}}g_{\SS^n}$ be the pull-back of $\tilde g$ to $\Phi(\tilde M_2)$ by $\Phi^{-1}$. Then $(\Phi(\tilde M_2),\tilde w^{\frac{4}{n-2}}g_{\SS^n})$ is a complete manifold. Hence, by an application of the Harnack inequality (see \cite[Proposition 2.6]{S-Y}), we have
\[
\tilde w(p) \geq c\,\dist(p,\partial\Phi(\tilde M_2))^{-\frac{n-2}{2}},\qquad p \in \Phi(\tilde M_2),
\]
where $c$ is some positive constant. Since $w$ $=$ $\tilde w|_{\Omega}$,
(iv) follows.  Theorem \ref{Structure} is established.
\eproof

\section{$C^0$ estimates}\label{C0Est}

In this section, we first prove Theorem \ref{By-product} and then use it in conjunction with Theorem \ref{Structure} to prove Theorem \ref{Main-Thm}. Throughout the section, unless otherwise stated, we will use $B_r(x)$ to denote the open ball of radius $r$ centered at $x$ in $\RR^n$, $n$ $\geq$ $3$.

For a positive $C^2$ function $u$, define
\[
A^u = -\frac{2}{n-2} u^{-\frac{n+2}{n-2}}\,\nabla^2 u + \frac{2n}{(n-2)^2}\,u^{-\frac{2n}{n-2}}\,\nabla u \otimes \nabla u - \frac{2}{(n-2)^2}\,u^{-\frac{2n}{n-2}}\,|\nabla u|^2\,I.
\]
As noted in the introduction, if $g_{flat}$ is the flat metric of $\RR^n$ and $\tilde g$ $=$ $u^{\frac{4}{n-2}}\,g_{flat}$, then $A_{\tilde g}$ $=$ $u^{\frac{4}{n-2}}A^u_{ij}\,dx^i\,dx^j$.

For the ease of exposition, we define a pair $(F,U)$ by
\begin{align*}
U
    &= \{M \in \textrm{Sym}^{n \times n}: \lambda(M) \in \Gamma\},\\
F(M)
    &= f(\lambda(M)).
\end{align*}
Here $\lambda(M)$ denotes the set of (real) eigenvalues of a symmetric matrix $M$ and $\textrm{Sym}^{n \times n}$ denotes the set of symmetric $n \times n$ matrices. Then the equation
\[
f(\lambda(A_{u^{\frac{4}{n-2}}\,g_{flat}})) = 1 \text{ and } \lambda(A_{u^{\frac{4}{n-2}}\,g_{flat}}) \in \Gamma
\]
is equivalent to
\[
F(A^u) = 1 \text{ and } A^u \in U.
\]
It is standard to check that \eqref{ConvexCone}-\eqref{f-homogeneity}, imply
\begin{align}
&\parbox{.8\textwidth}{$U$ is invariant under orthogonal conjugation, i.e.  $O^t\,U\,O$ $=$ $U$ for any orthogonal matrix $O$,} \label{U-OrthoConjInv}\\
&U \cap \{M + tN: t > 0\} \text{ is convex for any } M, N \in \textrm{Sym}^{n \times n}, N \geq 0, \label{U-Convex}\\
&{\rm tr}(M) := \sum_i M_{ii} \geq 0 \text{ for any } M \in U, \label{U-SuperHar}\\
&\parbox{.8\textwidth}{$F$ is invariant under orthogonal conjugation, i.e.  $F(O^t\,M\,O)$ $=$ $F(M)$ for any orthogonal matrix $O$ and $M$ $\in$ $U$,} \label{F-OrthoConjInv}\\
&\text{$F$ is homogeneous of degree one}, \label{F-Homog}\\
&F \in C^1(U) \text{ and }(F_{ij}(M)) > 0 \text{ for all } M \in U \text{ where } F_{ij}(M) = \frac{\partial F}{\partial M_{ij}}(M), \label{F-Ellipticity}
\end{align}
\begin{align}
&\text{there exists $\delta$ $>$ $0$ such that} \nonumber\\
&\qquad\qquad F(M) < 1 \text{ for all } M \in U, \|M\| = \sqrt{\sum M_{ij}^2} < \delta. \label{Nontriviality}
\end{align}

We will use the method of moving spheres,
a variant of the method of moving planes developed through the works of Alexandrov \cite{Alexandrov}, Serrin \cite{Se} and Gidas, Ni and Nirenberg \cite{G-N-N-1979}, \cite{G-N-N-1981}.
For a continuous function $w$, a point $x$ and a real number $\lambda$ $>$ $0$, let $w_x^\lambda$ denote its Kelvin transformation with respect to the sphere $\partial B_\lambda(x)$, i.e.
\[
w_x^\lambda(y) = \frac{\lambda^{n-2}}{|y-x|^{n-2}}\,w\Big(x + \frac{\lambda^2(y-x)}{|y-x|^2}\Big) \text{ wherever the expression makes sense}.
\]
We will also use the notation
\[
\psi_x^\lambda(y) = x + \frac{\lambda^2(y - x)}{|y - x|^2}.
\]

The proof of Theorem \ref{By-product} is split into two parts:
the ``interior'' estimate and the
``boundary'' estimate. We first treat the interior estimate, as it is easier to present and already contains
most of the ideas of the proof for the other part.

\subsection{Interior $C^0$ estimate}\label{InteriorEst}

\begin{proposition}\label{DistanceEst-Sing}
Let $(F,U)$ satisfy \eqref{U-OrthoConjInv}-\eqref{F-OrthoConjInv}, \eqref{F-Ellipticity}, and \eqref{Nontriviality} and $(\Omega,\Lambda)$ be as in Theorem \ref{By-product},
though allowing
$(\cup B_\alpha)\cup \Lambda=\emptyset$.
For any $\beta$ $>$ $0$, there exists a constant $C$ depending only on $n$, $(F, U)$ and $\beta$ such that if $u$ $\in$ $C^2(\Omega) \cap C^1(\bar\Omega\setminus\Lambda)$ is a positive solution of \eqref{Main-Eqn-Sing-Gen} for some $c(B)$, $c(B_\alpha)$ $>$ $-\beta$, there holds
\[
u(x) \leq C\,\dist(x,\partial\Omega)^{-\frac{n-2}{2}}\text{ for all } x \in \Omega.
\]
\end{proposition}

We first sketch the proof and then fill in the details
 later. Fix a point $x$ $\in$ $\Omega$. We rescale to get a ``unit size'' solution by
\begin{equation}
\left\{\begin{array}{rl}
\hat\Omega
    &= \{y: x + u(x)^{-\frac{2}{n-2}}\,y \in \Omega\},\,\hat\Lambda = \{y: x + u(x)^{-\frac{2}{n-2}}\,y \in \Lambda\},\\
\hat B
    &= \{y: x + u(x)^{-\frac{2}{n-2}}\,y\in B\},\, \hat B_{\alpha} = \{y: x + u(x)^{-\frac{2}{n-2}}\,y\in B_\alpha\},\\
\hat u(y)
    &= \displaystyle\frac{1}{u(x)} u\big(u(x)^{-\frac{2}{n-2}}\,y\big).
\end{array}\right.
\label{Rescale::Interior}
\end{equation}
By the conformally invariant property, $\hat u$ satisfies
\begin{equation}
\left\{\begin{array}{l}
F(A^{\hat u}) = 1 \text{ in }\hat \Omega,\\
A^{\hat u} \in U \text{ and } \hat u > 0 \text{ in }
\overline {\hat \Omega},\\
\hat u(x) \rightarrow \infty \text{ as } x \rightarrow \hat \Lambda,\\
\frac{\partial \hat u}{\partial\nu} + \frac{n-2}{2}\,\hat h\,\hat u \geq -\beta\,\hat u^{\frac{n}{n-2}} \text{ on }\partial \hat\Omega \setminus \hat \Lambda,\\
\hat u(0) = 1.
\end{array}\right.
\label{RescaleEqn::Interior}
\end{equation}
Here $\hat h$ is the mean curvature of $\partial\hat\Omega \setminus \hat\Lambda$, i.e.
\begin{equation}
\hat h = \left\{\begin{array}{l}
\frac{1}{u(x)^{\frac{2}{n-2}}\,r} \text{ on } \partial \hat B,\\
-\frac{1}{u(x)^{\frac{2}{n-2}}\,r_\alpha} \text{ on } \partial \hat B_\alpha.
\end{array}\right.
\label{MeanCurv}
\end{equation}
To finish the proof, we need to show that
\[
\dist(0,\partial\hat\Omega)= u(x)^{\frac 2{n-2}}\dist(x, \partial \Omega) \leq C(n,F,U,\beta).
\]

For this end, we employ the method of moving spheres. Define
\[
\bar\lambda = \sup\{0 < \lambda < \dist(0,\partial \hat\Omega) : \hat u_0^\lambda(y) \leq \hat u(y) \text{ for all } y \in \hat\Omega\setminus B_\lambda(0)\}.
\]
Note that $\bar\lambda$ is well-defined in light of the following lemma.
\begin{lemma}[\cite{L-Z}]\label{InitMovingSphere::Interior}
Let $D$ be an open subset of $\RR^n$ and $w$ $\in$ $C^{0,1}(D)$ satisfy
\[
\inf_{B_R(0) \cap D} w > 0 \text{ for any } R > 0.
\]
If $D$ is unbounded, assume in addition that
\[
\liminf_{|y| \rightarrow \infty} |y|^{n-2}w(y) > 0.
\]
For any $x$ $\in$ $D$, there exists $\lambda_1$ $>$ $0$ such that $B_{\lambda_1}(x)$ $\subset$ $D$ and for any $0$ $<$ $\lambda$ $<$ $\lambda_1$ and $y$ $\in$ $D\setminus B_\lambda(x)$, we have
\[
w_x^\lambda(y) \leq w(y).
\]
\end{lemma}

The remain of the proof is split into two independent parts.
\begin{enumerate}[\underline{Step \arabic{enumi}:}]
\item $\bar\lambda$ $\leq$ $C(n,F,U)$;
\item $\dist(0,\partial\hat\Omega)$ $\leq$ $C(n,\beta)(\bar\lambda^2 + \bar\lambda)$.
\end{enumerate}
Evidently, these imply the result.

In the proof of Step 2, we need the following gradient estimate:
\[
|\nabla\ln \hat u| \leq \frac{C(n)}{\bar\lambda} \text{ in } B_{\bar\lambda/2}(0).
\]
For the ease in exposition, we first present the proof of this estimate.

\subsubsection{Gradient estimates}

Interior gradient estimates were
 established in \cite{Li-CPAM}. Here we provide a simpler and self-contained proof for a somewhat different version of interior gradient estimates which
suffice for our purpose. The proof of Proposition \ref{DistanceEst-Sing} can be made shorter if we use the interior gradient estimates in \cite{Li-CPAM}, but the arguments therein
are considerably more involved.

\begin{lemma}\label{MovingSphere-GradEst}
Let $w$ $\in$ $C^0(B_{\lambda}(0))$, $\lambda$ $>$ $0$, be a positive function. Assume that
\[
w_z^\eta(y) \leq w(y) \text{ for any } B_\eta(z) \subset B_{\lambda}(0) \text{ and } y \in B_{\lambda}(0) \setminus B_\eta(z).
\]
Then $\ln w$ is locally Lipschitz in $B_{\lambda}(0)$ and
\[
|\nabla \ln w(x)| \leq \frac{n-2}{\lambda - |x|} \text{ for a.e. } x \in B_{\lambda}(0).
\]
\end{lemma}

\begin{remark}
Under a stronger assumption that $w$ is differentiable, this lemma was proved in \cite[Lemma A.2]{L-L-Acta}.
\end{remark}

\bproof Write $v$ $=$ $\ln w$. For $x$ in $B_{\lambda}(0)$, $r$ $=$ $\lambda - |x|$, and $0$ $<$ $\epsilon$ $<$ $\frac{r}{8}$, we will show that
\[
|v(y) - v(x)| \leq \frac{n-2}{2(\frac{r}{2} - 2\epsilon)}|y - x| \text{ for all } y \in B_{\epsilon}(x).
\]
Indeed, for $y$ $\in$ $B_{\epsilon}(x)$, $y$ $\neq$ $x$, let
\[
z_\pm = x \pm (\frac{r}{2} - \epsilon)\frac{y - x}{|y - x|} \text{ and } s_\pm = |y - z_\pm|.
\]
Note that
\[
y \in B_{\epsilon}(x) \subset B_{\frac{r}{2}}(z_\pm) \subset B_{\lambda}(x).
\]
Hence, for any $0$ $<$ $\eta$ $<$ $\frac{r}{2}$, we have
\[
w_{z_\pm}^\eta(p) \leq w(p) \text{ for any } p \in B_{\lambda}(0) \setminus B_\eta(z_\pm),
\]
which implies
\[
w\Big(z_\pm + s_1\frac{x - z_\pm}{|x - z_\pm|}\Big) \leq \Big(\frac{s_2}{s_1}\Big)^{\frac{n-2}{2}}\,w\Big(z_\pm + s_2\frac{x - z_\pm}{|x - z_\pm|}\Big) \text{ for any } 0 < s_1 \leq s_2 \leq \frac{r}{2}.
\]
Therefore, if we define
\[
g_\pm(s) = \frac{n-2}{2}\ln s + v\Big(z_\pm + s\frac{x - z_\pm}{|x - z_\pm|}\Big), 0 \leq s \leq \frac{r}{2},
\]
then $g_\pm(s)$ is increasing for $s$ $\in$ $[0,\frac{r}{2}]$. In particular, as $s_+$ $\leq$ $\frac{r}{2} - \epsilon$ $\leq$ $s_-$, we have
\[
g_+(s_+) \leq g_+\big(\frac{r}{2} - \epsilon\big) \text{ and } g_-\big(\frac{r}{2} - \epsilon\big) \leq g_-(s_-),
\]
which is equivalent to
\begin{multline*}
\frac{n-2}{2}\ln\big(\frac{r}{2} - \epsilon - |y - x|\big) + v(y) \leq \frac{n-2}{2}\ln\big(\frac{r}{2} - \epsilon\big) + v(x)\\
    \leq  \frac{n-2}{2}\ln\big(\frac{r}{4} - \epsilon + |y - x|\big) + v(y).
\end{multline*}
This implies that
\[
|v(y) - v(x)| \leq \frac{n-2}{2(\frac{r}{2} - 2\epsilon)}|y - x|.
\]
We conclude that $v$ is locally Lipschitz and thus differentiable almost every in $B_\lambda(0)$. Moreover, at point where $v$ is differentiable, we have
\[
|\nabla v|(x) \leq \frac{n-2}{r} = \frac{n-2}{\lambda - |x|}.
\]
This completes the proof.
\eproof

The following lemma is crucial in our approach to gradient estimates.

\begin{lemma}\label{MovingSphere-InteriorSphere}
Assume that $(F,U)$ satisfies \eqref{U-OrthoConjInv}, \eqref{U-Convex}, \eqref{F-OrthoConjInv} and \eqref{F-Ellipticity}. Let $D$ be an open subset of $\RR^n$, $n$ $\geq$ $3$, and let $u$ $\in$ $C^2(D)$ be a positive solution of
\begin{equation}
F(A^u) = 1 \text{ and } A^u \in U \text{ in } D.
\label{FNLEqn}
\end{equation}
Assume for some $\bar B_\lambda(x)$ $\subset$ $D$ that
\[
u_x^\lambda(y) \leq u(y) \text{ for all } y \in D\setminus B_\lambda(x).
\]
Then for any $B_\eta(z)$ $\subset$ $B_\lambda(x)$, there holds
\[
u_z^\eta(y) \leq u(y) \text{ for all } y \in D\setminus B_\eta(z).
\]
\end{lemma}

\bproof Define the function $\tilde u$ on $\RR^n$ by
\[
\tilde u(y) = \left\{\begin{array}{ll}
u(y) &\text{if } y \in B_\lambda(x),\\
u_x^\lambda(y) & \text{elsewhere}.
\end{array}\right.
\]
By conformal invariance, $\tilde u$ satisfies \eqref{FNLEqn} in $B_\lambda(x)$ and $\RR^n \setminus \partial B_\lambda(x)$. Moreover, $\tilde u$ $\leq$ $u$ in $D\setminus B_\lambda(x)$. As $\tilde u$ $=$ $u$ on $\partial B_\lambda(x)$, it follows that $\frac{\partial \tilde u}{\partial \nu}$ $\leq$ $\frac{\partial u}{\partial \nu}$ on $\partial B_\lambda(x)$, where $\nu$ is the outer unit normal to $\partial B_\lambda(x)$.

Fix $z$ $\in$ $B_\lambda(x)$. By Lemma \ref{InitMovingSphere::Interior}, there exists $\epsilon$ $>$ $0$ such that $\tilde u_z^\eta$ $\leq$ $\tilde u$ on $\RR^n\setminus B_\eta(z)$ for any $\eta$ $<$ $\epsilon$. Let
\[
\bar\eta = \sup\{\mu > 0: \tilde u_z^\eta \leq \tilde u \text{ on }\RR^n\setminus B_\eta(z) \text{ for all } 0 < \eta < \mu\}.
\]
To finish the proof, we only need to show that
\[
\bar\eta \geq \lambda - |z - x|.
\]
Arguing by contradiction, assume the converse so that $\bar\eta$ $<$ $\lambda - |z - x|$. We have
\[
\tilde u_z^{\bar\eta} \leq \tilde u \text{ on }\RR^n\setminus B_{\bar\eta}(z),
\]
which implies
\[
\tilde u_z^{\bar\eta} \geq \tilde u \text{ on } B_{\bar\eta}(z).
\]

We claim that we can find $q$ $\in$ $\partial B_\lambda(x)$ such that $\tilde u_z^{\bar\eta}(q)$ $=$ $\tilde u(q)$. Indeed, by the maximality of $\bar\eta$, one of the following two cases must occur:
\begin{enumerate}[(i)]
\item there exists $y$ $\in$ $B_{\bar \eta}(z)$ such that $\tilde u_z^{\bar\eta}(y)$ $=$ $\tilde u(y)$;
\item there exists $y$ $\in$ $\partial B_{\bar \eta}(z)$ such that $\frac{\partial\tilde u_z^{\bar\eta}}{\partial \nu_z}(y)$ $=$ $\frac{\partial\tilde u}{\partial \nu_z}(y)$ where $\nu_z$ is the outer unit normal to $\partial B_{\bar\eta}(z)$.
\end{enumerate}
Let
\[
B = \Big\{\zeta: \Big|z
 + \frac{\bar\eta^2(\zeta - z)}{|\zeta - z|^2} - x\Big| > \lambda\Big\} =
 (\psi_z^{\bar\eta})^{-1}(\RR^n\setminus \overline{ B_\lambda(x) }).
\]
Note that $B$ is a ball. By conformal invariance, we have
\[
\left\{\begin{array}{l}
F(A^{\tilde u}) = F(A^{\tilde u_z^{\bar \eta}}) = 1 \text{ in } B_{\bar\eta}(z) \setminus \partial B,\\
A^{\tilde u}, A^{\tilde u_z^{\bar\eta}} \in U \text{ in }
\overline{ B_{\bar\eta}(z)} \setminus \partial B,\\
\tilde u \leq \tilde u_z^{\bar \eta} \text{ in } B_{\bar \eta}(z),\\
\tilde u = \tilde u_z^{\bar \eta} \text{ on } \partial B_{\bar \eta}(z).
\end{array}\right.
\]
If Case (i) holds and $y$ $\notin$ $\partial B$, a standard argument using the strong maximum principle shows that $\tilde u$ and $\tilde u_z^{\bar \eta}$ are identically the same in $B_{\bar\eta}(z) \setminus B$ or in $B$. Indeed, we have $\tilde u(y)$ $=$ $\tilde u_z^{\bar\eta}(y)$, $D\tilde u(y)$ $=$ $D\tilde u_z^{\bar\eta}(y)$, $D^2\tilde u(y)$ $\leq$ $D^2\tilde u_z^{\bar\eta}(y)$, and so $A^{\tilde u}(y)$ $\geq$ $A^{\tilde u_z^{\bar\eta}}(y)$. Thus, near $y$, the function $w$ $=$ $\tilde u_z^{\bar\eta} - \tilde u$ $\geq$ $0$ satisfies
\[
0 = F(A^{\tilde u}) - F(A^{\tilde u_z^{\bar\eta}}) = a_{ij}\nabla_{ij} w + b_i\,\nabla_i w + cw := Lw,
\]
where $(a_{ij})$ $>$ $0$, $b_i$ and $c$ are continuous near $y$. Here we have used \eqref{U-Convex} and \eqref{F-Ellipticity}. If Case (ii) holds, another standard argument using the Hopf lemma implies that $\tilde u$ and $\tilde u_z^{\bar \eta}$ are identically the same in $B_{\bar\eta}(z) \setminus B$. To see this, it suffices to notice that $\tilde u$ $=$ $\tilde u_z^{\bar\eta}$ on $\partial B_{\bar\eta}(z)$ which implies $D\tilde u(y)$ $=$ $D\tilde u_z^{\bar\eta}(y)$, $D^2\tilde u(y)$ $\leq$ $D^2\tilde u_z^{\bar\eta}(y)$, and $A^{\tilde u}(y)$ $\geq$ $A^{\tilde u_z^{\bar\eta}}(y)$. In either case, we must have that $\tilde u(y)$ $=$ $\tilde u_z^{\bar \eta}(y)$ for some $y$ $\in$ $\partial B$, which implies the claim.

We have showed that $\tilde u_z^{\bar \eta}$ must touch $\tilde u$ from below at some point $q$ $\in$ $\partial B_\lambda(x)$, i.e. $\tilde u_z^{\bar \eta}(q)$ $=$ $\tilde u(q)$ and $\tilde u_z^{\bar \eta}$ $\leq$ $\tilde u$ near $q$. Recalling the definition of $\tilde u$ and noting that $\frac{\partial \tilde u}{\partial \nu}$ $\leq$ $\frac{\partial u}{\partial \nu}$ on $\partial B_\lambda(x)$
and $u^{\bar\eta}_z$ is $C^1$ near $q$, we infer that $\frac{\partial \tilde u}{\partial \nu}(q)$ $=$ $\frac{\partial u}{\partial \nu}(q)$. As
\[
\left\{\begin{array}{l}
F(A^{\tilde u}) = F(A^u) = 1 \text{ in } D\setminus \bar B_{\lambda}(x),\\
A^{\tilde u}, A^u \in U \text{ in } D\setminus B_{\lambda}(x),\\
\tilde u \leq u \text{ in } D\setminus B_{\lambda}(x),\\
\tilde u = u \text{ on } \partial B_\lambda(x),
\end{array}\right.
\]
the Hopf lemma again implies that $u$ and $\tilde u$ are identically the same near $\partial B_\lambda(x)$ (in $D\setminus B_\lambda(x)$). Therefore, $\tilde u$ $\equiv$ $\tilde u_x^\lambda$ is $C^2$ and is an entire solution of \eqref{FNLEqn}.

We thus have
\[
\left\{\begin{array}{l}
F(A^{\tilde u}) = F(A^{\tilde u_z^{\bar \eta}}) = 1 \text{ in } \RR^n \setminus B_{\bar\eta}(z),\\
A^{\tilde u}, A^{\tilde u_z^{\bar\eta}} \in U \text{ in } \RR^n \setminus B_{\bar\eta}(z),\\
\tilde u \geq \tilde u_z^{\bar \eta} \text{ in } \RR^n \setminus B_{\bar\eta}(z),\\
\tilde u(q) = \tilde u_z^{\bar \eta}(q) \text{ for some } q \in \partial B_\lambda(x).
\end{array}\right.
\]
By the strong maximum principle, we infer that $\tilde u$ $\equiv$ $\tilde u_z^{\bar\eta}$ in $\RR^n \setminus B_{\bar\eta}(z)$, and so in $\RR^n$.

In the remaining, we show that $\tilde u$ $\equiv$ $\tilde u_x^\lambda$ $\equiv$ $\tilde u_z^{\bar\eta}$ leads to a contradiction. We compute for $\xi$ $\in$ $B_\lambda(x)$ and $p$ $:=$ $x + \frac{\lambda^2(\xi - x)}{|\xi - x|^2}$ $=$ $\psi_x^\lambda(\xi)$,
\begin{align*}
\tilde u(\xi)
    &= (\tilde u_z^{\bar\eta})_x^\lambda(\xi) = \frac{\lambda^{n-2}}{|\xi - x|^{n-2}}\frac{\bar\eta^{n-2}}{|p - z|^{n-2}}\,\tilde u\Big(\psi_z^{\bar\eta}(p)\Big)\\
    &= \Big[\frac{|p - x|}{\lambda}\frac{\bar\eta}{|p - z|}\Big]^{n-2}\,\tilde u\Big(\psi_z^{\bar\eta}(p)\Big)\\
    &\leq \Big[\frac{|p - x|}{\lambda}\frac{\bar\eta + |z - x|}{|p - z| + |z - x|}\Big]^{n-2}\,\tilde u\Big(\psi_z^{\bar\eta}(p)\Big)\\
    &\leq \Big[\frac{\bar\eta + |z - x|}{\lambda}\Big]^{n-2}\,\tilde u\Big(\psi_z^{\bar\eta}(p)\Big).
\end{align*}
Observe that the map $\xi$ $\mapsto$ $\psi_z^{\bar\eta}(p)$ $=$ $\psi_z^{\bar\eta} \circ \psi_x^\lambda (\xi)$ maps $\bar B_{\bar\eta}(z)$ into itself and so has a fixed point $\xi_*$ $\in$ $\bar B_{\bar\eta}(z)$. It follows that
\[
\tilde u(\xi_*)
    \leq \Big[\frac{\bar\eta + |z - x|}{\lambda}\Big]^{n-2}\,\tilde u(\xi_*),
\]
which contradicts our earlier assumption that $\bar\eta + |z - x|$ $<$ $\lambda$.
\eproof

As a consequence of Lemma \ref{MovingSphere-GradEst} and Lemma \ref{MovingSphere-InteriorSphere}, we have:
\begin{corollary}\label{MovingSphere-FinerGradEst}
Assume that $(F,U)$ satisfies \eqref{U-OrthoConjInv}, \eqref{U-Convex}, \eqref{F-OrthoConjInv} and \eqref{F-Ellipticity}. Let $D$ be an open subset of $\RR^n$, $n$ $\geq$ $3$, and let $u$ $\in$ $C^2(D)$ be a positive solution of
\[
F(A^u) = 1 \text{ and } A^u \in U \text{ in } D.
\]
If for some $\bar B_\lambda(x)$ $\subset$ $D$ there holds
\[
u_x^\lambda(y) \leq u(y) \text{ for all } y \in D\setminus B_\lambda(x),
\]
then
\[
|\nabla u(y)| \leq \frac{n-2}{\lambda - |y - x|}\,u(y) \text{ for all } y \in B_{\lambda}(x).
\]
\end{corollary}

\subsubsection{Step 1 of the proof of Proposition \ref{DistanceEst-Sing}}

\begin{lemma}\label{NontrivialSolution}
Assume that $(F,U)$ satisfies \eqref{U-SuperHar} and \eqref{Nontriviality}. Let $u$ $\in$ $C^2(B_s(0))$, $s$ $>$ $0$, be a positive solution of
\[
F(A^u) \geq 1 \text{ and } A^u \in U \text{ in } B_s(0).
\]
Assume for some positive constant $C_1$ that
\[
u(y) \geq C_1 \text{ in } B_s(0).
\]
Then
\[
s \leq C(n)\delta^{-1/4}\,C_1^{-\frac{2}{n-2}},
\]
where $\delta$ is the constant in \eqref{Nontriviality}.
\end{lemma}

\begin{remark}
Under an additional hypothesis that $u(y)$ $\leq$ $C_2$ in $B_s(0)$ for some $C_2$, it was shown in \cite{L-L-Acta} that $s$ $\leq$ $C(n,C_1,C_2,\delta)$.
\end{remark}

Lemma \ref{NontrivialSolution} is equivalent to:
\setcounter{lemm}{\thelemma}
\addtocounter{lemm}{-1}
\begin{lemm}
Assume that $(F,U)$ satisfies \eqref{U-SuperHar} and \eqref{Nontriviality}. Let $u$ $\in$ $C^2(B_1(0))$ be a positive solution of
\[
F(A^u) \geq 1 \text{ and } A^u \in U \text{ in } B_1(0).
\]
Then
\[
\inf_{B_1(0)} u \leq C(n,\delta).
\]
\end{lemm}

Indeed, the equivalence follows by considering $\tilde u(x)$ $=$ $s^{\frac{n-2}{2}}\,u(sx)$ and using the fact that $A^{\tilde u}(x)$ $=$ $A^{u}(sx)$.

Before giving the proof of Lemma \ref{NontrivialSolution} we point out that under an additional hypothesis that there exists some $\alpha$ $>$ $0$ such that
\[
{\rm tr}(M) \geq \alpha \text{ for all } M \in U \text{ satisfying } F(M) \geq 1,
\]
Lemma \ref{NontrivialSolution} can be treated simplier by ODE method. For in this case we must have
\[
-\Delta u \geq \frac{n-2}{2}\alpha\,u^{\frac{n+2}{n-2}} \text{ and } u \geq C_1 \text{ in } B_s(0),
\]
which implies the radial average $\bar u(r)$ $:=$ $\frac{1}{|\partial B_r|}\int_{\partial B_r} u\,dx$ satisfies
\[
-\Delta \bar u \geq \frac{n-2}{2}\alpha\,\bar u^{\frac{n+2}{n-2}} \text{ and } \bar u \geq C_1 \text{ in } B_s(0).
\]

\medskip

\bproof Let $\tau$ be the largest positive number such that
\[
\xi(y) := \frac{\tau}{s^4}(s^2 - |y|^2)^2 \leq u \text{ in } B_s(0).
\]
Then for some $|\bar y|$ $<$ $s$,
\[
u(\bar y) = \xi(\bar y),
\]
and so
\begin{align*}
\nabla u(\bar y)
    &= \nabla \xi(\bar y) = \frac{-4\tau}{s^4}(s^2 - |\bar y|^2)\bar y,\\
\nabla^2 u(\bar y)
    &\geq \nabla^2\xi(\bar y) = \frac{8\tau}{s^4}\bar y \otimes \bar y - \frac{4\tau}{s^4}(s^2 - |\bar y|^2)I.
\end{align*}
Therefore
\begin{align*}
A^\xi(\bar y)
    &= \frac{2}{n-2}\,\xi(\bar y)^{-\frac{n+2}{n-2}}\Big[-\nabla^2\xi(\bar y) + \frac{n}{n-2}\frac{\nabla\xi(\bar y) \otimes \nabla\xi(\bar y)}{\xi(\bar y)} - \frac{1}{n-2}\frac{|\nabla\xi(\bar y)|^2}{\xi(\bar y)}I\Big]\\
    &= \frac{16\tau}{(n-2)s^4}\,\xi(\bar y)^{-\frac{n+2}{n-2}}\Big[\frac{n+2}{n-2}\,\bar y \otimes \bar y + \Big(\frac{s^2 - |\bar y|^2}{2} - \frac{2}{n-2}|\bar y|^2\Big)I\Big].
\end{align*}
Using \eqref{U-SuperHar}, we get
\[
0 \leq {\rm tr}(A^u(\bar y)) \leq {\rm tr}(A^\xi(\bar y)) = \frac{8\tau}{(n-2)s^4}\,\xi(\bar y)^{-\frac{n+2}{n-2}}\Big[ns^2 - (n+2)|\bar y|^2\Big],
\]
which implies
\[
|\bar y|^2 \leq \frac{n}{n+2}s^2 \text{ and therefore } \xi(\bar y) \geq
\frac{4}{(n+2)^2}\tau.
\]
As we also have $\xi(\bar y)$ $=$ $u(\bar y)$ $\geq$ $C_1$, it follows that
\[
\xi(\bar y) \geq
(\frac{4}{(n+2)^2}\tau)^{ \frac{n-2}{n+2}}
(C_1)^{ 1-\frac{n-2}{n+2}}=
 C(n)\,C_1^{\frac{4}{n+2}}\,\tau^{\frac{n-2}{n+2}}.
\]
Consequently, by our previous calculation,
\begin{align*}
A^u(\bar y)
    &\leq A^\xi(\bar y) = \frac{16\tau}{(n-2)s^4}\,\xi(\bar y)^{-\frac{n+2}{n-2}}\Big[\frac{n+2}{n-2}\,\bar y \otimes \bar y + \Big(\frac{s^2 - |\bar y|^2}{2} - \frac{2}{n-2}|\bar y|^2\Big)I\Big]\\
    &\leq C(n)\frac{C_1^{-\frac{4}{n-2}}}{s^2}\,I.
\end{align*}
On the other hand, from $F(A^u(\bar y))$ $\geq$ $1$, \eqref{U-SuperHar} and \eqref{Nontriviality}, we have
\[
\|A^u\| \geq \delta \text{ and } {\rm tr}(A^u) \geq 0.
\]
Combining the above estimates we obtain
\[
\frac{C_1^{-\frac{4}{n-2}}}{s^2} \geq C(n)\sqrt{\delta} > 0.
\]
The assertion follows.
\eproof

\begin{remark}
The same proof yields the following conclusion: If $u$ $\in$ $C^2(B_s(0))$ satisfies for some $p$ $\geq$ $1 + \epsilon$, $\epsilon$ $>$ $0$, and $C_1$ $>$ $0$
\[
F(u^{\frac{n+2}{n-2} - p}A^u) \geq 1, u^{\frac{n+2}{n-2} - p}A^u \in U  \text{ and } u \geq C_1 \text{ in } B_s(0),
\]
then
\[
s \leq C(n,\delta,\epsilon)\,C_1^{-\frac{p-1}{2}}.
\]
\end{remark}

\subsubsection{Step 2 of the proof of Proposition \ref{DistanceEst-Sing}}

The following lemma describes the contact set when we carry out the method of moving spheres at an interior point.

\begin{lemma}\label{MovingSphere-ContactSet::Interior}
Assume that $(F,U)$ satisfies \eqref{U-OrthoConjInv}, \eqref{U-Convex}, \eqref{F-OrthoConjInv} and \eqref{F-Ellipticity}. Let $D$ be a connected open subset of $\RR^n$, $\Lambda$ a (possibly empty) closed subset of $\partial D$, and let $u$ $\in$ $C^2(D) \cap C^0(\bar D \setminus \Lambda)$ be a positive solution of
\[
\left\{\begin{array}{l}
F(A^u) = 1 \text{ and } A^u \in U \text{ in } D,\\
\displaystyle\lim_{y \rightarrow \Lambda} u = +\infty \text{ if } \Lambda \neq \emptyset.
\end{array}\right.
\]
Assume for some $\bar B_\eta(z)$ $\subset$ $D$ that
\begin{equation}
u_z^\eta(y) \leq u(y) \text{ for all } y \in D\setminus B_\eta(z).
\label{Moving::Interior}
\end{equation}
Let $C$ $=$ $C_\eta(z)$ denote the contact set
\[
C = \{y \in \bar D\setminus \bar B_\eta(z): u_z^\eta(y) = u(y)\}.
\]
Then one of the following (mutually exclusive) cases must occur:
\begin{enumerate}[(a)]
\item $C$ is empty and
\[
\frac{\partial u}{\partial\nu} > \frac{\partial u_z^\eta}{\partial\nu} \text{ on } \partial B_\eta(z),
\]
where $\nu$ is the outward unit normal to $\partial B_\eta(z)$;
\item $C$ $=$ $\bar D \setminus \bar B_\eta(z)$ and $\Lambda$ $=$ $\emptyset$;
\item $\emptyset$ $\neq$ $C$ $\subset$ $\partial D\setminus \Lambda$.
\end{enumerate}
Moreover, if (a) holds and $\eta$ is the largest number such that \eqref{Moving::Interior} is satisfied, then $D$ is unbounded and ``touching occurs at infinity'', i.e.
\[
\liminf_{|y| \rightarrow \infty} |y|^{n-2}[u(y) - u_x^\eta(y)] = 0.
\]
\end{lemma}

\bproof We follow the argument given in \cite{L-L-CPAM} for solutions of conformally invariant equations. The idea is to use the
strong maximum principle and the Hopf lemma as in the proof of Lemma \ref{MovingSphere-InteriorSphere}.

If $C$ is non-empty, then $u_z^\eta$ touches $u$ from below somewhere. If interior touching happens, the strong maximum principle implies $u$ $\equiv$ $u_z^\eta$ outside $B_\eta(z)$, which gives (b). If interior touching does not happen, (c) holds.

If $C$ is empty, (a) follows from the Hopf lemma. Moreover if $D$ is bounded or if $D$ is unbounded and $\liminf_{|y| \rightarrow \infty} |y|^{n-2}[u(y) - u_x^\eta(y)]$ $>$ $0$, we can find by a continuity and compactness argument some $\epsilon$ $>$ $0$ such that
\[
u_x^\lambda(y) \leq u(y) \text{ for all } \eta < \lambda < \eta + \epsilon \text{ and } y \in D \setminus B_\lambda(x).
\]
The last assertion follows.
\eproof

Before going on to the proof of Proposition \ref{DistanceEst-Sing}, we state and prove an elementary result concerning the geometry of a sphere.

\begin{lemma}\label{SimpleGeom}
Let $B$ $=$ $B_r(x)$ be a ball in $\RR^n$, $n$ $\geq$ $1$. Let $\nu^-$ denote the inward unit normal along the boundary of $B$.

\medskip\noindent
(a) If $z$ $\notin$ $B$, then for any $y$ $\in$ $\partial B$,
\[
\frac{|y-z|^2}{2r} + (y-z) \cdot \nu^-(y) \geq \dist(z,\partial B).
\]
The equality holds if and only if $z$ $\in$ $\partial B$.

\medskip\noindent
(b) If $z$ $\in$ $\bar B$, then for any $y$ $\in$ $\partial B$,
\[
-\frac{|y-z|^2}{2r} - (y-z) \cdot \nu^-(y) \geq \frac{\dist(z,\partial B)}{2}.
\]
The equality holds if and only if $z$ $\in$ $\partial B \cup \{x\}$.
\end{lemma}

\bproof We can assume without loss of generality that $z$ is the origin.

\medskip\noindent
(a) We calculate using the cosine law,
\begin{align*}
\frac{|y|^2}{2r} + y \cdot \nu^-(y)
    &= \frac{|y|^2}{2r} + y \cdot \frac{x-y}{r} = \frac{-|y|^2 + 2 x \cdot y}{2r}\\
    &= \frac{-|y|^2 + (|x|^2 + |y|^2 - r^2)}{2r} = \frac{|x|^2 - r^2}{2r} = |x| - r + \frac{(|x| - r)^2}{2r}\\
    &\geq |x| - r = \dist(0,\partial B).
\end{align*}

\medskip\noindent
(b) Similarly,
\begin{align*}
-\frac{|y|^2}{2r} - y \cdot \nu^-(y)
    &= -\frac{|y|^2}{2r} - y \cdot \frac{x-y}{r} = \frac{|y|^2 - 2 x \cdot y}{2r}\\
    &= \frac{|y|^2 - (|x|^2 + |y|^2 - r^2)}{2r} = \frac{-|x|^2 + r^2}{2r} = \frac{r - |x|}{2} + \frac{|x|(r - |x|)}{2r}\\
    &\geq \frac{r - |x|}{2} = \frac{\dist(0,\partial B)}{2}.
\end{align*}
\eproof

\noindent{\bf Proof of Proposition \ref{DistanceEst-Sing}.} Fix a point $x$ $\in$ $\Omega$. Define $\hat\Omega$, $\hat\Lambda$, $\hat B$, $\hat B_\alpha$, and $\hat u$ by \eqref{Rescale::Interior}.
Then, by conformal invariance, $\hat u$ satisfies \eqref{RescaleEqn::Interior}. We need to show that
\[
u(x)^{\frac{2}{n-2}}\dist(x,\partial\Omega) = \dist(0,\partial\hat\Omega) \leq C(n,\delta,\beta).
\]

Using Lemma \ref{InitMovingSphere::Interior}, we define
\[
\bar\lambda = \sup\{0 < \lambda < \dist(0,\partial \hat\Omega) : \hat u_0^\lambda(y) \leq \hat u(y) \text{ for all } y \in \hat\Omega\setminus B_\lambda(0)\}.
\]
By Corollary \ref{MovingSphere-FinerGradEst}, we have
\begin{equation}
|\nabla\ln \hat u| \leq \frac{C(n)}{\bar\lambda} \text{ in } B_{\bar\lambda/2}(0).
\label{DistEst-Sing::GradEst}
\end{equation}
As $\hat u(0)$ $=$ $1$, this implies that
\begin{equation}
C(n) \geq \hat  u \geq C(n)^{-1} > 0 \text{ in } B_{\bar\lambda/2}(0).
\label{DistEst-Sing::C0Bound}
\end{equation}

\noindent\underline{Step 1:} By Lemma \ref{NontrivialSolution}, we must have
\[
\bar\lambda \leq C(n,\delta).
\]

\noindent\underline{Step 2:} If $\dist(0,\partial\hat\Omega)$ $\leq$ $2C(n,\delta)$, with the same
$C(n, \delta)$ above, we are done. Otherwise,
$\dist(0, \partial \hat \Omega)\ge 2\bar\lambda$, and
 by Lemma \ref{MovingSphere-ContactSet::Interior}, we can find $y$ $\in$ $\partial\hat\Omega$, $|y|$ $>$ $\bar\lambda$ such that $\hat u(y)$ $=$ $\hat u_{0}^{\bar\lambda}(y)$. Let $B_*$ be the ball in the family $\{\hat B\} \cup \{\hat B_\alpha\}$ such that  $y\in \partial B_*$,  and let $h_*$ be the mean curvature
of $\partial B_*$ at $y$ (with respect to the inner normal). Then
\[
\bar\lambda \leq \frac{1}{2}\dist(0,\partial\hat\Omega) \leq \frac{1}{2}\dist(0,\partial\hat B_*),
\]
and so
\[
\Big|\bar\lambda^2\,\frac{y}{|y|^2}\Big| \leq \frac{\bar\lambda}{2}.
\]
Hence, by \eqref{DistEst-Sing::GradEst} and \eqref{DistEst-Sing::C0Bound},
\[
\hat u\big(\bar\lambda^2\,\frac{y}{|y|^2}\big) \leq C(n) \text{ and } \Big|D\ln\hat u\big(\bar\lambda^2\,\frac{y}{|y|^2}\big)\Big| \leq \frac{C(n)}{\bar\lambda}.
\]

For simplicity, we write $\hat u^{\bar\lambda}$ $\equiv$ $\hat u_0^{\bar\lambda}$. By a straightforward calculation using the expression for $\hat u^{\bar\lambda}$ and the above gradient estimate, we get
\begin{align*}
\frac{\partial \hat u^{\bar\lambda}}{\partial \nu}(y) + \frac{n-2}{2}h_*\,\hat u^{\bar\lambda}(y)
    &\leq \hat u^{\bar\lambda}(y)\Big[-(n-2)\frac{y \cdot \nu}{|y|^2} + C(n)\,\frac{1}{\bar\lambda}\,\frac{\bar\lambda^2}{|y|^2}\Big] + \frac{n-2}{2}h_*\,\hat u^{\bar\lambda}(y)\\
    &\leq -\frac{(n-2)}{|y|^2}\,\hat u^{\bar\lambda}(y)\Big[y \cdot \nu - \frac{|y|^2}{2}h_* - C(n)\,\bar\lambda\Big].
\end{align*}
Applying Lemma \ref{SimpleGeom}, we arrive at
\[
\frac{\partial \hat u^{\bar\lambda}}{\partial \nu}(y) + \frac{n-2}{2}h_*\,\hat u^{\bar\lambda}(y)
    \leq -\frac{(n-2)}{2|y|^2}\,\hat u^{\bar\lambda}(y)\Big[\dist(0,\partial B^*) - C(n)\,\bar\lambda\Big].
\]
Recalling the expression for $u^{\bar\lambda}$ and the upper bound for $u$ in $B_{\bar\lambda/2}(x)$, we infer that
\[
\frac{\partial \hat u^{\bar\lambda}}{\partial \nu}(y) + \frac{n-2}{2}h_*\,\hat u^{\bar\lambda}(y)
    \leq -\frac{C_1(n)}{\bar\lambda^2}\,[\hat u^{\bar\lambda}(y)]^{\frac{n}{n-2}}\Big[\dist(0,\partial B^*) - C(n)\,\bar\lambda\Big].
\]

On the other hand, as $\hat u^{\bar\lambda}$ $\leq$ $\hat u$ near $y$ and $\hat u^{\bar\lambda}(y)$ $=$ $\hat u(y)$,
\[
\frac{\partial \hat u^{\bar\lambda}}{\partial \nu}(y) + \frac{n-2}{2}h_*\,\hat u^{\bar\lambda}(y)
    \geq \frac{\partial \hat u}{\partial \nu}(y) + \frac{n-2}{2}h_*\,\hat u(y) \geq -\beta\,\hat u(y)^{\frac{n}{n-2}} = -\beta\,[\hat u^{\bar\lambda}(y)]^{\frac{n}{n-2}}.
\]
Combining with the preceeding estimate, we get
\[
-\beta \leq -\frac{C_1(n)}{\bar\lambda^2}\,\Big[\dist(0,\partial B^*) - C(n)\,\bar\lambda\Big],
\]
which implies
\[
\dist(0,\partial B^*) \leq \frac{\beta}{C_1(n)}\,\bar\lambda^2 + C_1(n)C(n)\,\bar\lambda \leq C(n,F,U,\beta).
\]
It follows that
\[
\dist(0,\partial\hat \Omega) \leq C(n,F,U,\beta),
\]
and the conclusion follows readily.
\eproof

\subsection{The proof of Theorem \ref{By-product}}

Loosely speaking, by Proposition \ref{DistanceEst-Sing}, in order to establish Theorem \ref{By-product}, we only need to focus on establishing $C^0$ bound near the boundary. The idea is to adapt the proof of Proposition \ref{DistanceEst-Sing}. The main difference is that here we do not have a simple gradient estimate as in the interior case. We have to resort to the gradient estimate established in \cite{Li-CPAM}.
We will present the proof by contradiction arguments, though
a direct argument can be given similarly.

\medskip
\noindent{\bf Proof of Theorem \ref{By-product}.} Suppose the contrary, then for some $\beta>0$ and $\epsilon>0$, there exists a sequence of solutions $\{u_i\}
\subset C^3(\Omega)\cap C^2(\overline \Omega\setminus \Lambda)$ satisfying (\ref{Main-Eqn-Sing-Gen}), with $c_i(B),
c_i(B_\alpha)>-\beta$, and for some $\bar x_i\in \overline \Omega$,
\begin{equation}
\dist(\bar x_i, \Lambda)\ge \epsilon,
\label{E1-1}
\end{equation}
\begin{equation}
u_i(\bar x_i)=\max \{ u_i(x)\ :\
x\in \overline \Omega, \dist (x, \Lambda)\ge \epsilon\}
\ \to \infty.
\label{E1-2}
\end{equation}

By the interior estimate (Proposition \ref{DistanceEst-Sing}),
\begin{equation}
u_i(\bar x_i) \dist(\bar x_i, \partial \Omega)^{ \frac{n-2}2}
\le C.
\label{E2-0}
\end{equation}
Here and below, we use $C>1$ to denote some
positive constant independent of $i$.
It follows that
\begin{equation}
\dist(\bar x_i, \partial \Omega)\to 0.
\label{E2-1}
\end{equation}

Because of (\ref{E1-1}), there exists
$\{B_i:=B_{\alpha_i}\}\subset \{B_\alpha\}$,
$\bar x_i'\in \partial B_i$, such that
\begin{equation}
|\bar x_i -\bar x_i'|=
\dist(\bar x_i, \partial \Omega)\to 0.
\label{E2-2}
\end{equation}

By (\ref{E1-1}) and (\ref{E2-2}),
$\{\bar x_i'\}$ is $\epsilon/2-$distance away
from $\Lambda$ for large $i$.
It follows that
\begin{equation}
r_i:= \text{radius of}\ B_i\ge \frac 1C,
\label{E3-1}
\end{equation}
\begin{equation}
\dist(\bar x_i', \partial\Omega\setminus \partial B_i)\ge \frac 1C.
\label{E3-2}
\end{equation}
The reason is that otherwise a subsequence of $\{\bar x_i'\}$
together with
a subsequence of $\{B_\alpha\}$ would converge to a point, which has to be in
$\Lambda$.  This would violate the fact that
$\{\bar x_i'\}$ is of a fixed distance away from $\Lambda$.

Let $\psi_i$ be a M\"{o}bius transformation which maps
$\RR^n\setminus \overline {B_i^+}$ to
$$
\RR^n_+:= \{x\ :\
x=(x', x_n)=(x_1, \cdots,
x_{n-1}, x_n)\in \RR^n, x_n>0\}
$$
such that
$$
\psi_i(\bar x_i')=0,\
\psi_i(\bar x_i)=\tilde x_i:= (0', T_i), \ T_i:=\frac{r_i}{r_i + |\bar x_i - \bar x_i'|}\,|\bar x_i-\bar x_i'|,
\text{and }\psi_i(\bar x_i'')=\infty,
$$
where $\bar x_i''$ is the antipodal point of $\bar x_i'$ on $\partial B_i$.

Set
\begin{align*}
\tilde \Lambda_i
    &=\psi_i(\Lambda\setminus \{\bar x_i''\}),\\
\tilde \Omega_i
    &= \psi_i(\Omega\setminus\Lambda),\\
\tilde u_i
    &= \tilde u_{\psi_i^{-1}}= |Jac_{ \psi_i^{-1}}|^{ \frac {n-2}{2n}}u\circ \psi_i^{-1}.
\end{align*}
Then, depending on whether there is another ball in $\{B_\alpha\}$ that touches $B_i$ at $\bar x_i''$ or not, we can write $\tilde\Omega_i$ as either
\[
\RR^n_+ \setminus (\cup \tilde B_\alpha^i \cup \tilde\Lambda_i)
\]
or
\[
\RR^n_+ \setminus (\{(x',x_n) \in \RR^n_+, x_n \geq H_i\} \cup \tilde B_\alpha^i \cup \tilde\Lambda_i),
\]
where in the first case $\tilde B_\alpha^i$ are non-overlapping balls contained in $\RR^n_+$ and in the second case they are non-overlapping balls contained in $\{(x',x_n) \in \RR^n_+, x_n < H_i\}$. By abuse of notation, we will ambiguously write
\[
\tilde\Omega_i = \RR^n_+ \setminus (\cup \tilde B_\alpha^i \cup \tilde\Lambda_i),
\]
by which we view half-spaces as balls of infinite radius and centered at infinity.

By (\ref{E3-1}) and (\ref{E3-2}),
\begin{equation}
\dist(0, (\partial \tilde \Omega_i\setminus
\partial  \RR^n_+)\cup \tilde \Lambda_i)
\ge 1/C.
\label{E6-1}
\end{equation}

By the conformal invariance of the equation of $u$,
$\tilde u_i$
 $\in$ $C^3(\tilde\Omega_i) \cap C^2(\overline{\tilde\Omega_i}
\setminus\tilde\Lambda_i)$ satisfies
\[
\left\{\begin{array}{l}
F(A^{\tilde u_i}) = 1 \text{ in } \tilde\Omega_i,\\
A^{\tilde u_i} \in U \text{ and } \tilde u_i > 0\text{ in }
\overline{\tilde \Omega_i} \setminus \tilde \Lambda_i,\\
\lim_{ x\to x_0}
\tilde u_i(x)=  \infty \ \forall\ x_0\in
 \tilde\Lambda_i,\\
\frac{\partial \tilde u_i}{\partial x_n}
 =  c_i(B_i)\,\tilde u_i^{\frac{n}{n-2}}
\text{ on } \partial  \RR^n_+ \setminus \tilde\Lambda_i,\\
\frac{\partial \tilde u_i}{\partial\nu} + \frac{n-2}{2\tilde r_\alpha^i}\,
\tilde u_i \geq - \beta\,\tilde u_i^{\frac{n}{n-2}}
\text{ on }\partial\tilde B_\alpha^i \setminus \tilde\Lambda_i,\\
\displaystyle \liminf_{|y| \rightarrow \infty} |y|^{n-2}\,\tilde u_i(y)
 = +\infty, \text{ or }\\
\text{$(\tilde u_i)_0^{R_i}$
 extends to a $C^2$ positive function in
 $\bar B_{R_i}^+(0)$ for some $R_i>0$}.
\end{array}\right.
\]
In the above, $\tilde r _\alpha^i=$ radius of $\tilde B_\alpha^i$.

It is clear that
\begin{equation}
\tilde u_i (x)\le C \tilde u_i(0', T_i),
\quad \forall \ |x|\le \epsilon/C,\
x\in \RR^n_+,
\label{E8-2}
\end{equation}
\begin{equation}
C u_i(\bar x_i)\ge \tilde u_i(0', T_i)\ge \frac 1C
u_i(\bar x_i)\to \infty.
\label{E8-3}
\end{equation}
Moreover, by \eqref{E2-0},
\begin{equation}
\tilde u_i(0', T_i) T_i^{ \frac {n-2}2 }\le C,
\label{E8-1}
\end{equation}

Let
$$
\tilde M_i= \tilde u_i(0', T_i).
$$
We define
\begin{align*}
\hat\Omega_i
    &= \{y: \tilde M_i^{-\frac{2}{n-2}}\,y \in \tilde\Omega_i\},\\
\hat\Lambda_i
    &= \{y: \tilde M_i^{-\frac{2}{n-2}}\,y \in \tilde\Lambda_i\},\\
\hat B_{\alpha}^i
    &= \{y: \tilde M_i^{-\frac{2}{n-2}}\,y\in \tilde B_\alpha^i\},\\
\hat u_i(y)
	&= \frac{1}{\tilde M_i}\,\tilde u_i
\big(\tilde M_i^{-\frac{2}{n-2}}\,y\big), \qquad y \in \hat \Omega_i.
\end{align*}
By the conformal invariance, $\hat u_i$ satisfies
\begin{equation}
\left\{\begin{array}{l}
F(A^{\hat u_i}) = 1 \text{ in }\hat\Omega_i,\\
A^{\hat u_i} \in U \text{ and } \hat u_i > 0\text{ in }
 \overline{\hat \Omega_i} \setminus \hat \Lambda_i,\\
\lim_{x\to x_0}
\hat u_i(x)=  \infty \ \forall\
x_0\in  \hat\Lambda_i,\\
\frac{\partial \hat u_i}{\partial x_n}  =  c_i(B_i)
\,\hat u_i^{\frac{n}{n-2}} \text{ on } \partial \RR^n_+
 \setminus \hat\Lambda_i,\\
\frac{\partial \hat u_i}{\partial\nu} + \frac{n-2}{2 \hat r_\alpha^i}\,
\,\hat u_i \geq - \beta\,\hat u_i^{\frac{n}{n-2}} \text{ on }\partial
\hat B_\alpha^i \setminus \hat\Lambda_i,\\
\displaystyle \liminf_{|y| \rightarrow \infty} |y|^{n-2}
\,\hat u_i(y) = +\infty \text{ or }\\
\text{$(\hat u_i)_0^{R_i}$ extends to a $C^2$ positive
function in $\bar B_{R_i}^+(0)$ for some
$R_i>0$}.
\end{array}\right.
\label{E10-1}
\end{equation}
In the above, $\hat r_\alpha^i=$ radius of
$\hat B_\alpha^i$.

We know from (\ref{E8-1}), (\ref{E8-3}), (\ref{E8-2}) and (\ref{E6-1})
that
\begin{align}
&\hat T_i := T_i \tilde M_i^{  \frac 2{n-2} }= T_i \tilde u_i(0', T_i)^{ \frac 2{n-2} }
        \le C T_i u_i(\bar x_i)^{ \frac 2{n-2}} \le C,\nonumber\\
&\hat u_i(0', \hat T_i) =1,\nonumber\\
&\hat u_i(x)\le C \hat u_i(0', \hat T_i) = C, \quad \forall\ |x|\le \frac \epsilon C \tilde M_i ^{ \frac 2{n-2}}, \ x\in \RR^n_+,\nonumber\\
&\dist(0, \cup \hat B_\alpha ^i \cup \hat \Lambda_i) \ge \frac{1}{C}\tilde M_i ^{ \frac 2{n-2} }\to \infty.
\label{E12-1}
\end{align}

By the local gradient estimates in
\cite[Theorem 1.19]{Li-CPAM},
\begin{equation}
|\nabla \hat u_i(x)|\le C
\hat u_i(x)\le C,\quad
\forall\ |x|\le \frac \epsilon{2C} \hat M_i^{ \frac 2{n-2}},\
x\in \RR^n_+.
\label{E13-1}
\end{equation}

Define
$$
\bar \lambda_i=\sup
\{0<\lambda<\dist(0, \cup \hat B_\alpha^i\cup \hat \Lambda_i)\ :\
(\hat u_i)_0^\lambda(y)\le \hat u_i(y)\ \text{for all}\
y\in \hat \Omega_i\setminus B_\lambda(0)\}.
$$
This definition is justified by the following lemma.

\begin{lemma}[\cite{L-Z}]\label{InitMovingSphere::Bdry}
Let $D$ be an open subset of $\RR^n_+$ containing $\{(x',x_n): |x'| < s, 0 < x_n < \epsilon\}$ for some $\epsilon$ $>$ $)$. Let $w$ be a locally Lipschitz function in $D \cup B_s'(0)$, $B_s'(0)$ $:=$ $\{(x',0): |x'| < s\}$, satisfying
\[
\inf_{B_R(0) \cap D} w > 0 \text{ for any } R > 0.
\]
If $D$ is unbounded, assume in addition that
\[
\liminf_{|y| \rightarrow \infty} |y|^{n-2}w(y) > 0.
\]
Then for any $x$ $\in$ $B_s'(0)$, there exists some $\lambda_1$ $>$ $0$ such that $B_{\lambda_1}^+(x)$ $\subset$ $D$ and for any $0$ $<$ $\lambda$ $<$ $\lambda_1$ and $y$ $\in$ $D\setminus B_\lambda(x)$, we have
\[
w_x^\lambda(y) \leq w(y).
\]
\end{lemma}

\noindent{\it Claim.}\ $\displaystyle{
\lim_{i\to \infty}\bar \lambda_i=\infty}.$

\medskip

\noindent{\it Proof of the claim.}\ Suppose not, then
$$
\bar\lambda_i\le C,\qquad \forall\ i.
$$
We know that
\begin{equation}
(\hat u_i)_0^{  \bar \lambda_i }(y)\le \hat u_i(y)\ \
\text{for all}\ y\in \hat \Omega_i\setminus B_{ \bar \lambda_i}(0).
\label{E15-1}
\end{equation}
Arguing as usual, using the strong maximum principle
and the Hopf lemma, we conclude that there exists some point
$y_i\in \partial \hat B_{\alpha_i}^i$ such that
\begin{equation}
(\hat u_i)_0^{ \bar \lambda_i}(y_i)=\hat u_i(y_i).
\label{E16-1}
\end{equation}
For readers' convenience, we outline the arguments.  If
$$
\frac {\partial }{\partial \nu}[\hat u_i-
(\hat u_i)_0^{ \bar \lambda_i}]=0
\ \ \text{at some point of}\ \partial B_{ \bar \lambda_i}(0),
$$
then by the Hopf lemma and the strong maximum principle,
(\ref{E16-1}) must occur (and
$\hat \Lambda_i=\emptyset$).
So we may assume
\begin{equation}
\frac {\partial }{\partial \nu}[\hat u_i-
(\hat u_i)_0^{ \bar \lambda_i}]>0
\ \text{on}\
 \partial B_{ \bar \lambda_i}(0),
\label{E16-A-1}
\end{equation}

If
$$
\liminf_{|y|\to \infty} |y|^{n-2}
[\hat u_i(y)-
(\hat u_i)_0^{ \bar \lambda_i}(y)]=0,
$$
then line 7 of (\ref{E10-1}) occurs.
After a Kelvin transformation, we can apply the strong maximum principle,
and the Hopf lemma if necessary,
to
show that (\ref{E16-1}) must occur (and $\hat \Lambda_i=\emptyset$).
So we may assume
\begin{equation}
\liminf_{|y|\to \infty} |y|^{n-2}
[\hat u_i(y)-
(\hat u_i)_0^{ \bar \lambda_i}(y)]>0.
\label{E16-B-1}
\end{equation}

Using the third line of (\ref{E10-1}) and the strong maximum
principle, we have, in case (\ref{E16-1}) does not occur,
\begin{equation}
\inf_{ y\in \hat \Omega_i,
\bar\lambda_i+\epsilon<|y|<R}
[\hat u_i(y)-
(\hat u_i)_0^{ \bar \lambda_i}(y)]>0,
\quad \forall\ 0<\epsilon<R<\infty.
\label{E16-B-2}
\end{equation}
Using (\ref{E16-A-1}),
(\ref{E16-B-1}) and (\ref{E16-B-2}), we can show that
for some $\epsilon>0$,
$$
(\hat u_i)_0^ \lambda (y)\le \hat u_i(y)\ \
\text{for all}\ y\in \hat \Omega_i\setminus B_\lambda(0)\
\text{and all}\ \lambda\le \bar \lambda_i+\epsilon.
$$
This violates the definition of $\bar \lambda_i$.
So we have proved (\ref{E16-1}).
 Details can be found in \cite[Lemma 3]{L-L-JEMS}. The main difference between this
 lemma and Lemma \ref{MovingSphere-ContactSet::Interior}
is to prevent touching on
 $\partial \RR^n_+$ and to prove a version of the Hopf lemma
on $\partial B_\eta(z) \cap \partial \RR^n_+$.

On the other hand, in view of (\ref{E13-1}) and (\ref{E12-1}),
(\ref{E16-1}) cannot occur by the arguments in the proof of
Proposition \ref{DistanceEst-Sing}.
The claim is proved.

To complete the proof of Theorem \ref{By-product}, we will show that contrary to the above claim, $\lambda_i$ is uniformly bounded. For this purpose we start with a few lemmas. First is an analogue of Lemma \ref{MovingSphere-InteriorSphere}.

\begin{lemma}\label{MovingSphere-InteriorSphere::Bdry}
Assume that $(F,U)$ satisfies \eqref{U-OrthoConjInv}, \eqref{U-Convex}, \eqref{F-OrthoConjInv} and \eqref{F-Ellipticity}. Let $D$ be a connected open subset of $\RR^n_+$ such that $\RR^{n-1}$ is a component of $\partial D$. Let $u$ $\in$ $C^2(\bar D)$ be a positive solution of
\[
\left\{\begin{array}{l}
F(A^u) = 1 \text{ and } A^u \in U \text{ in } \bar D,\\
\frac{\partial u}{\partial x_n} = c\,u^{\frac{n}{n-2}} \text{ on }
 \partial \RR^n_+ \text{ for some constant $c$.}
\end{array}\right.
\]

Assume for some $x$ $\in$ $
\partial \RR^n_+$ and $B_\lambda^+(x)$ $\subset$ $D$ that
\[
u_x^\lambda(y) \leq u(y) \text{ for all } y \in D\setminus B_\lambda^+(x).
\]
Then for any $z$ $\in$ $\partial \RR^n_+
$ and $B_\eta^+(z)$ $\subset$ $B_\lambda^+(x)$, there holds
\[
u_z^\eta(y) \leq u(y) \text{ for all } y \in D\setminus B_\eta^+(z).
\]
\end{lemma}

\bproof The proof is
similar to that of Lemma \ref{MovingSphere-InteriorSphere}. Here one has to prevent touching on $\partial
\RR^n_+$, which can be done by using the Hopf
lemma and the Neumann boundary condition as in
the proof of the claim below Lemma \ref{InitMovingSphere::Bdry}.
\eproof

We retailor Lemma \ref{NontrivialSolution} into the following two lemmas. In effect, the first one is also an improvement of Lemma \ref{NontrivialSolution}.

\begin{lemma}\label{NontrivialSolution::Bdry-EZ}
Assume that $(F,U)$ satisfies \eqref{U-SuperHar} and \eqref{Nontriviality}. Let $u$ $\in$ $C^2(B_s^+(0)) \cap C^0(\bar B_s^+(0))$ be a positive solution of
\[
F(A^u) \geq 1 \text{ and } A^u \in U \text{ in } B_s^+(0).
\]
Assume for some positive constant $C_1$ that
\[
u \geq C_1 \text{ on } B_s(0) \cap \{x_n = 0\}
\]
Then
\[
s \leq C(n,\delta)C_1^{-\frac{2}{n-2}},
\]
where $\delta$ is the constant in \eqref{Nontriviality}.
\end{lemma}

\begin{remark}
Under an additional hypothesis that $u$ $\leq$ $C_2$ in $B_s^+(0)$, this result was proved in \cite{L-L-PrivComm}. Here we use Lemma \ref{NontrivialSolution} to relax this extra hypothesis.
\end{remark}

\bproof Let $\rho$ $\in$ $C^\infty_c(-1,1)$ be a nice cut-off function such that $0$ $\leq$ $\rho$ $\leq$ $1$ and $\rho$ $=$ $1$ in $(-\frac{1}{2},\frac{1}{2})$. Let $\psi$ be the solution of
\begin{align*}
\Delta\psi
    &= 0 \text{ in } B_s^+(0)\\
\psi
    &= 0 \text{ on } \partial B_s(0) \cap \RR^n_+\\
\psi
    &= \rho\big(\frac{|x'|}{s}\big) \text{ on } B_{\frac{s}{2}}(0) \cap \{x_n = 0\}.
\end{align*}
Then $\psi$ $\geq$ $c(n)$ in $B_{\frac{s}{2}}^+(0)$. Hence, by our assumption on the lower bound of $u$ on $B_s(0) \cap \{x_n = 0\}$ and the maximum principle,
\[
u \geq C_1\psi \geq c(n)C_1 \text{ in } B_{\frac{s}{2}}^+(0).
\]
Applying Lemma \ref{NontrivialSolution} to any ball of radius $\frac{s}{8}$ which is contained in $B_{\frac{s}{2}}^+(0)$, we get the assertion.
\eproof

\begin{lemma}\label{NontrivialSolution::Bdry}
Assume that $(F,U)$ satisfies \eqref{U-SuperHar} and \eqref{Nontriviality}. Let $u$ $\in$ $C^2(B_s^+(0)) \cap C^0(\bar B_s^+(0))$ be a positive solution of
\[
F(A^u) = 1 \text{ and } A^u \in U \text{ in } B_s^+(0),
\]
which satisfies, for some constants $C_1$ $>$ $0$ and $0$ $\leq$ $T$ $\leq$ $\frac{s}{8}$,
\[
u(0',T) \geq C_1.
\]
Moreover, assume that
\[
u_{0}^{s/2}(y) \leq u(y) \text{ for any } y \in B_s^+(0)\setminus B_{s/2}^+(0).
\]
Then
\[
s \leq C(n,\delta)\,C_1^{-\frac{2}{n-2}},
\]
where $\delta$ is the constant in \eqref{Nontriviality}.
\end{lemma}

\bproof By Lemma \ref{NontrivialSolution::Bdry-EZ}, it suffices to shows that
\[
u(x',2T) \geq C(n)C_1 \text{ for any } |x'| \leq \frac{s}{8}.
\]
To see this, fix $|x'|$ $\leq$ $\frac{s}{8}$. By Lemma \ref{MovingSphere-InteriorSphere::Bdry}, we have for any $\lambda$ $\leq$ $s-|x'|$,
\[
u_{(-x',0)}^\lambda(y) \leq u(y) \text{ for any } y \notin B_\lambda^+(-x',0).
\]
Therefore, for $e$ $=$ $\frac{(0,T) - (-x',0)}{|(0,T) - (-x',0)|}$, $r^{\frac{n-2}{2}}u(re)$ is monotonically increasing for $r$ $<$ $\frac{s}{4}$. It follows that
\[
u(x',2T) \geq \frac{1}{2^{\frac{n-2}{2}}} u(0,T) \geq c(n)C_1.
\]
The proof is complete.
\eproof

\bigskip
\noindent{\bf Proof of Theorem \ref{By-product} (continued).}\label{ProofCont} Taking $s=\bar \lambda_i$ in
Lemma \ref{NontrivialSolution::Bdry}, we arrive at
$$
\bar \lambda_i\le C.
$$
This contradicts the claim.  Theorem \ref{By-product} is established.
\eproof

As a corollary of the proofs of Lemmas \ref{NontrivialSolution::Bdry-EZ} and \ref{NontrivialSolution::Bdry}, we obtain the following result which is of independent interest, though not needed in this paper.

\begin{corollary}
Assume that $(F,U)$ satisfies \eqref{U-SuperHar} and \eqref{Nontriviality}. Let $u$ $\in$ $C^2(B_s^+(0)) \cap C^0(\bar B_s^+(0))$ be a positive solution of
\[
F(A^u) \geq 1 \text{ and } A^u \in U \text{ in } B_s^+(0),
\]
Moreover, assume for $B_s'(0)$ $:=$ $B_s(0) \cap \{x_n = 0\}$that
\[
u_{x}^{\eta}(y) \leq u(y) \text{ for any } x \in B_s'(0), 0 < \eta < s - |x| \text{ and } y \in B_s^+(0)\setminus B_{\eta}^+(x).
\]
Then
\[
\sup_{B_{s/4}^+(0)} u \leq C(n,\delta)\,s^{\frac{2-n}{2}},
\]
where $\delta$ is the constant in \eqref{Nontriviality}.
\end{corollary}

\bproof It suffices to show that
\[
u(x) < C(n,\delta)\,s^{-\frac{n-2}{2}} \text{ for any } x = (0',x_n) \text{ with } x_n < \frac{s}{4}.
\]

Define, for $R$ $>$ $0$,
\begin{align*}
W_R(x)
    &:= \Big\{y = x + t\,\frac{(z',x_n)}{\sqrt{|z'|^2 + x_n^2}} \text{ for some } 0 \leq t \leq R, |z'| \leq \frac{s}{4}\Big\},\\
\partial_L W_R(x)
    &:= \Big\{y = x + t\,\frac{(z',x_n)}{\sqrt{|z'|^2 + x_n^2}} \text{ for some } 0 \leq t \leq R, |z'| = \frac{s}{4}\Big\}.
\end{align*}
Arguing as in the proof of Lemma \ref{NontrivialSolution::Bdry}, we get
\[
u(y) \geq C(n)\,u(x) \text{ for any } y \in \partial_L W_{\frac{s}{4}}(x).
\]
Using \eqref{U-SuperHar} and comparing $u$ to a harmonic function as in the proof of Lemma \ref{NontrivialSolution::Bdry-EZ}, we get
\[
u(y) \geq C(n)\,u(x) \text{ for any } y \in W_{\frac{s}{8}}(x).
\]
Applying Lemma \ref{NontrivialSolution} to any ball of radius $\frac{s}{100}$ contained in $W_{\frac{s}{8}}(x)$, we get
\[
s \leq C(n,\delta)\,u(x)^{-\frac{2}{n-2}}.
\]
The assertion follows.
\eproof

\subsection{The proof of Theorem \ref{Main-Thm}}

We next move on to the proof of Theorem \ref{Main-Thm}. Recall that $(M,g)$ is a locally compact complete conformally flat manifold with umbilic boundary $\partial M$ whose components are $N_1$, \ldots, $N_m$.
By (a) in Definition \ref{Type},
$\lambda_1(M, \tilde g)>0$, since, in addition to
$h_{\tilde g}\ge 0$ on $\partial M$,
$\lambda(A_{\tilde g})\in \Gamma
\subset \Gamma_1$ implies
$R_{\tilde g}>0$.
Hence $\lambda_1(M,g)$ $>$ $0$. 
By Theorem \ref{Structure}, we can find a ball $B_{-1}$ $\subset$ $\RR^n$ and (possibly empty, finite or infinite) collection of non-overlapping balls $\{B_\alpha\}_{\alpha = 1, 2, \ldots}$, each being contained in $B_{-1}$, and a closed subset $\Lambda$ of $\bar B_{-1}$, whose Hausdorff dimension does not exceed $\frac{n-2}{2}$, such that
\begin{enumerate}[(i)]
\item If $\Omega$ $=$ $B_{-1} \setminus (\cup  B_\alpha \cup \Lambda)$,
then there exists a conformal map $\Psi:$ $\Omega$ $\rightarrow$ $(M,g)$ which is a covering map.
\item If $\bar B_\alpha \cap \bar B_\beta$ $=$ $\{p\}$,
or $\bar B_\alpha\cap \bar B_{-1}=\{p\}$,
then $p$ $\in$ $\Lambda$.
\item If $\{B_{\alpha_j}\}$ $\subset$ $\{B_\alpha\}$ is a sequence of distinct balls ``converging'' to a point $p$, then $p$ $\in$ $\Lambda$.
\item If we write the pull-back metric of $g$ to $\Omega$ by $\Psi$ as $w^{\frac{4}{n-2}}\,\delta_{ij}$, then
\[
w(x) \rightarrow \infty \text{ as } \dist_{\RR^n}(x,\Lambda) \rightarrow 0.
\]
\end{enumerate}

Then for any positive solution $u$ of \eqref{Main-Eqn}, the function
\[
v(y) := u \circ \Psi(y)\,w(y), \qquad y \in \Omega,
\]
satisfies
\begin{equation}
\left\{\begin{array}{l}
F(A^v) = 1\text{ in } \Omega,\\
A^v \in U \text{ and }v > 0 \text{ in } \bar\Omega,\\
v \rightarrow \infty \text{ near } \Lambda,\\
\frac{\partial v}{\partial \nu} + \frac{n-2}{2}\,h\,v = c_k\,v^{\frac{n}{n-2}} \text{ on } \partial B_\alpha \setminus \Lambda \text{ if } \Psi(\partial B_\alpha\setminus\Lambda) = N_k,
\end{array}\right.
\label{Main-Eqn-Sing}
\end{equation}

\begin{proposition}\label{UpperBound}
Let $(M^n,g)$, $n$ $\geq$ $3$, be a smooth compact locally conformally flat Riemannian manifold with umbilic boundary $\partial M$ and $N_1$, \ldots, $N_m$ be the components of $\partial M$. Assume that (M,g) is not conformally equivalent to the standard half-sphere $\SS^n_+$. Let $(f,\Gamma)$ satisfy \eqref{ConvexCone}-\eqref{f-homogeneity} and assume that $\lambda_1(M,g)$ $>$ $0$. For any $\beta$ $>$ $0$, there exists a constant $C$ depending only on $n$, $(M,g)$, $(f,\Gamma)$ and $\beta$ such that if $u$ $\in$ $C^2(M)$ is a positive solution of \eqref{Main-Eqn} for some $c_1$, \ldots, $c_m$ $>$ $-\beta$, there holds
\begin{equation}
\sup_M u \leq C.
\label{35-1}
\end{equation}
\end{proposition}

\bproof As shown above, we can use Theorem \ref{Structure} to find a ball $B_{-1}$ $\subset$ $\RR^n$ and a (possibly empty, finite or infinite) collection of non-overlapping balls $\{B_\alpha\}_{\alpha = 1, 2, \ldots}$, each being contained in $B_{-1}$, and a closed subset $\Lambda$ of $\bar B_{-1}$, whose Hausdorff dimension does not exceed $\frac{n-2}{2}$, such that the properties (i)-(iv) listed above are satisfied.

Next, fix some point $p_0$ $\in$ $\Omega$ and pick $R$ large enough such that
\[
\Psi(E_R) = M,
\]
where $E_R$ $=$ $\{p: \dist_{(\Omega,\Psi^*(g))}(p,p_0) \leq R\}$. Such $R$ is guaranteed to exist by the compactness of $M$. Evidently,
\[
\dist_{\RR^n}(E_R,\Lambda) > \epsilon > 0.
\]

We distinguish two cases: $(\Omega,\Lambda)$ $\not\cong$ $(\SS^n_+,\emptyset)$ or $(\Omega,\Lambda)$ $\cong$ $(\SS^n_+,\emptyset)$.

\medskip
\noindent\underline{Case 1:} $(\Omega,\Lambda)$ $\not\cong$ $(\SS^n_+,\emptyset)$.

Define
\[
v(y) = u \circ \Psi(y)\,w(y), \qquad y \in \Omega.
\]
Then $v$ satisfies \eqref{Main-Eqn-Sing}. By Theorem \ref{By-product}, we have
\[
\sup_{E_R} v \leq C(n,(M,g),(f,\Gamma),\beta),
\]
which implies
\[
\sup_M u \leq C(n,(M,g),(f,\Gamma),\beta)\big[\inf_{E_R} w\big]^{-1} = C(n,(M,g),(f,\Gamma),\beta).
\]

\medskip
\noindent\underline{Case 2:} $(\Omega,\Lambda)$ $\cong$ $(\SS^n_+,\emptyset)$. In particular, we can assume that $\Omega$ $=$ $\RR^n_+$.

Assume that the conclusion fails so that we can construct a family of solution $u_j$ to \eqref{Main-Eqn} on $M$ such that
\[
u_j(x_j) = \sup_{M} u_j \rightarrow \infty.
\]
Furthermore, assume that $x_j$ $\rightarrow$ $x_*$ $\in$ $M$.

Note that since $(M,g)$ is not conformally equivalent to $\SS^n$, for any $x$ $\in$ $M$, $\Psi^{-1}(x)$ $\subset$ $\bar \RR^n_+$ contains at least two points. Hence, by compactness, we can select $y_j$ $\rightarrow$ $y_*$ and $z_j$ $\rightarrow$ $z_*$ in $\bar\RR^n_+$ with $y_*$ $\neq$ $z_*$ such that $\Psi(y_*)$ $=$ $\Psi(z_*)$ $=$ $x_*$.

Now let
\[
v_j(y) = u_j \circ \Psi(y)\,w(y), \qquad y \in \RR^n_+.
\]
Then $v_j$ blows up at $y_*$ and $z_*$.

On the other hand, $v_j$ satisfies
\[
\left\{\begin{array}{l}
F(A^{v_j}) = 1\text{ in } \RR^n_+,\\
A^{v_j} \in U \text{ and }v_j > 0 \text{ in } \bar\RR^n_+,\\
\frac{\partial v_j}{\partial \nu} = c\,v_j^{\frac{n}{n-2}} \text{ on } \RR^{n-1}.
\end{array}\right.
\]
By the Liouville theorem established in \cite{L-L-JEMS}, $v_j$ must be a ``standard bubble'', i.e.
\[
v_j(p) = c(n)\Big(\frac{a_j}{1 + a_j^2|p - p_j|^2}\Big)^{\frac{n-2}{2}} \text{ for some } a_j > 0, p_j \in \RR^n
\]
This implies in particular that $v_j$ can only blow up at a single point. This contradicts our earlier conclusion that $v_j$ blows up at $y_*$ and $z_*$, which are separated.
\eproof

\noindent{\bf Proof of Theorem \ref{Main-Thm}.} In view of Proposition \ref{UpperBound}, it suffices to show that
\[
\inf_M u \geq C > 0.
\]

Consider the metric $g_\epsilon$ $=$ $(1 - \epsilon\,\varphi)^{\frac{4}{n-2}}\,g$ where $\epsilon$ is a small positive number to be determined and $\varphi$ is a smooth function such that $\varphi(x)$ $=$ $\dist_g(x,\partial M)$ near $\partial M$. Then for all $\epsilon$ sufficiently small, $\lambda(A_{g_\epsilon})$ $\in$ $\Gamma$ in $M$. But as $\varphi$ $=$ $0$ and $\partial_\nu \varphi$ $=$ $1$ on $\partial M$, we have
\[
h_{g_\epsilon} = -\frac{2}{n-2}\,\partial_\nu(1-\epsilon\varphi) + h_g = \frac{2\epsilon}{n-2} + h_g > 0 \text{ on } \partial M.
\]
Therefore, by replacing $g$ by $g_\epsilon$ if necessary, we can assume from the beginning that
\[
\inf_{\partial M} h_g \geq C > 0.
\]

We next show that
\[
\sup_M u \geq C > 0.
\]
Let $u(\bar x)$ $=$ $\max_{M} u$. Then either $\bar x$ $\in$ $M$ or $\bar x$ $\in$ $\partial M$.

If $\bar x$ $\in$ $M$, $\nabla u(\bar x)$ $=$ $0$ and $\nabla^2 u(x)$ $\leq$ $0$. Thus,
\[
1 = f\big(\lambda(A_{u^{\frac{4}{n-2}}g}(\bar x))\big) \geq f\big(u^{-\frac{4}{n-2}}(\bar x)\lambda(A_g(\bar x))\big) = u^{-\frac{4}{n-2}}(\bar x)f(\lambda(A_g(\bar x))) .
\]
Since $\lambda(A_g)$ $\in$ $\Gamma$ in $M$, this implies that
\[
u(\bar x) \geq \min_{M} \big[f(\lambda(A_g))\big]^{\frac{n-2}{4}} = C > 0.
\]

If $\bar x$ $\in$ $N_k$ $\subset$ $\partial M$, we have
\[
\frac{\partial u}{\partial \nu}(\bar x) + \frac{n-2}{2}h_g u(\bar x) = c_k u^{\frac{n}{n-2}}(\bar x).
\]
Since $\bar x$ is a maximum point of $u$ and $h_g$ is positive, this implies that $c_k$ $>$ $0$ and
\[
\frac{n-2}{2}\,h_g(\bar x)\,u(\bar x) \leq c_k\,u^{\frac{n}{n-2}}(\bar x),
\]
which in turn gives
\[
u(\bar x) \geq \min_{\partial M}\Big[\frac{n-2}{2c_k}h_g\Big]^{\frac{n-2}{2}} = C > 0.
\]

In either case, we have shown that
\[
\sup_M u = u(\bar x) \geq C > 0.
\]

To finish the proof,
we invoke the gradient estimates in \cite[Theorem 1.19]{Li-CPAM} to
obtain, in view of (\ref{35-1}),
\[
\sup_M |\nabla \ln u| \leq C.
\]
Evidently, this estimate implies
\[
\inf_M u \geq C\,\sup_M u \geq C > 0.
\]
The proof is complete.
\eproof

\section{General existence results}\label{ExistProof}

In this section we give the proof of Theorems \ref{Existence} and of Proposition \ref{PositiveView-High-k}.
 
\medskip
\noindent{\bf Proof of Theorem \ref{Existence}.} First notice that, by the results in \cite[Appendix B]{L-L-Acta}, we can assume without loss of generality that \eqref{f-homogeneity} holds. Then, by \eqref{f-concavity}, \eqref{delta} also holds. 

We first prove the estimate
\eqref{C4-Est}. By Theorem \ref{Main-Thm},
\[
\|u\|_{C^0(M)} + \|u^{-1}\|_{C^0(M)} \leq C.
\]
By known  $C^1$ and $C^2$ estimates as described in the introduction,
see e.g. \cite[Theorem 1.19]{Li-CPAM}, \cite{L-L-CPAM} (the proof of (1.39) there), \cite[Theorems 1.3 and 1.5]{J-L-L}
and \cite[Theorem 3]{Chen}, we have
\[
\|u\|_{C^2(M)} + \|u^{-1}\|_{C^2(M)} \leq C.
\]
Here, for the $C^2$ estimates, we have used (\ref{f-concavity})
and $c_k\ge 0$; the latter is for the boundary
$C^2$ estimates.
 \eqref{C4-Est} follows from
Evans-Krylov's estimates (\cite{Ev}, \cite{Kr}) and the Schauder theory.

We now turn to the existence part. If $M$ is conformally equivalent to the standard half-sphere, there is nothing to prove. We thus assume that $M$ is not conformally equivalent to $\SS^n_+$. Moreover, since $\lambda_1(M,g)$ $>$ $0$, we can assume, after making a conformal change of the metric using a first eigenfunction of \eqref{7-1}, that $R_g$ $>$ $0$ in $M$ and $h_g$ $\geq$ $0$ on $\partial M$. To finish the proof, it suffices to show that
\begin{equation}
\deg((G,B),D,0) = -1,
\label{DegCount}
\end{equation}
where $\deg$ denotes the degree defined in Appendix \ref{App-Degree}, $D$ is an appropriately chosen open bounded subset of $C^{4,\alpha}(M)$ and
\begin{align*}
G[u]
	&= f(\lambda(A_{u^{\frac{4}{n-2}}g})) - 1,\\
B[u](x)
	&= \partial_\nu u + \frac{n-2}{2}h_g(x) u - c_{k}u^{\frac{n}{n-2}} \text{ for } x \in N_k.
\end{align*}

For $0$ $\leq$ $t$ $\leq$ $1$, define, as in \cite{L-L-CPAM},
\begin{align*}
\Gamma_t
    &= \{\lambda \in \RR^n| t\lambda + (1-t)\sigma_1(\lambda)e \in \Gamma\},\\
f_t(\lambda)
    &= f(t\lambda + (1-t)\sigma_1(\lambda)e) \text{ for } \lambda \in \Gamma_t,\\
c_{k,t}
    &= \left\{\begin{array}{l}tc_k + (1-t) \text{ if } c_k > 0,\\0\text{ if } c_k = 0,\end{array}\right.
\end{align*}
where $e$ $=$ $(1,\ldots, 1)$ $\in$ $\RR^n$.
See also 
\cite{GV03-Indiana}
for a similar homotopy.

Consider
\begin{equation}
\left\{\begin{array}{l}
f_t(\lambda(A_{u^{\frac{4}{n-2}}g})) = 1 \text{ in }M,\\
\lambda(A_{u^{\frac{4}{n-2}}g}) \in \Gamma_t \text{ and } u > 0 \text{ in } M,\\
\frac{\partial u}{\partial\nu} + \frac{n-2}{2}h_g\,u = c_{k,t}\,u^{\frac{n}{n-2}} \text{ on }\partial N_k, k = 1 \ldots m,
\end{array}\right.
\label{Htpy-Eqn}
\end{equation}

By \eqref{C4-Est}, there exists some $C$ $>$ $0$ independent of $t$ such that for all solutions of \eqref{Htpy-Eqn}, there holds
\[
\|u\|_{C^{4,\alpha}(M)} + \|u^{-1}\|_{C^{4,\alpha}(M)} \leq C.
\]
Therefore, as $f|_{\partial\Gamma}$ $=$ $0$, there exists $\epsilon$ $>$ $0$ independent of $t$ such that
\[
\dist(\lambda(A_{u^{\frac{4}{n-2}}g}),\partial\Gamma_t) \geq 2\epsilon.
\]

Define
\begin{align*}
D_t
	&= \Big\{u \in C^{4,\alpha}(M): u > 0, \lambda(A_{u^{\frac{4}{n-2}}g}) \in \Gamma_t,\\
		&\qquad\qquad \|u\|_{C^{4,\alpha}(M)} + \|u^{-1}\|_{C^{4,\alpha}(M)} \leq 2C, \dist(\lambda(A_{u^{\frac{4}{n-2}}g}),\partial\Gamma_t) > \epsilon\Big\},\\
G_t[u]
	&= f_t(\lambda(A_{u^{\frac{4}{n-2}}g})) - 1,\\
B_t[u](x)
	&= \partial_\nu u + \frac{n-2}{2}h_g(x) u - c_{k,t}u^{\frac{n}{n-2}} \text{ for } x \in N_k.
\end{align*}
Note that $(G_1,B_1)$ $=$ $(G,B)$. By the properties of the degree (see Appendix \ref{App-Degree}), $\deg((G_t,B_t),D_t,0)$ is defined for $0$ $\leq$ $t$ $\leq$ $1$ and 
\begin{equation}
\deg((G_0,B_0),D_0,0) = \deg((G_1,B_1),D_1,0).
\label{Exist::Eqn01}
\end{equation}

Set
\begin{align*}
\tilde G_t[u]
	&= G_0[u]\,u^{(1-t)\,\frac{n+2}{n-2}},\\
\tilde B_t(x,u)
	&= \partial_\nu u + \frac{n-2}{2}h_g(x) u - t\,c_{k,0}u^{\frac{n}{n-2}} \text{ for } x \in N_k.
\end{align*}
Note that $(\tilde G_t, \tilde B_t)(u)$ $=$ $0$ amounts to the Yamabe problem with boundary. (In particular, higher derivative estimates follows directly from $C^0$ estimates due to its semi-linear structure.) Thus, $\deg((\tilde G_t,\tilde B_t),D_0,0)$ is defined for $0$ $\leq$ $t$ $\leq$ $1$ and
\begin{equation}
\deg((\tilde G_0,\tilde B_0),D_0,0) = \deg((\tilde G_1,\tilde B_1),D_0,0) = \deg((G_0,B_0),D_0,0).
\label{Exist::Eqn02}
\end{equation}

Let $L_g$ be the conformal Laplacian of $(M,g)$, i.e.
\[
L_g = \Delta_g - \frac{n-2}{4(n-1)}R_g.
\]
For a function $v$ $\in$ $C^{2,\alpha}(M)$, let $(-L_g)^{-1} v$ be the solution to
\[
\left\{\begin{array}{l}
-L_g \phi = v \text{ in } M^\circ,\\
\partial_\nu \phi + \frac{n-2}{2}h_g(x) \phi = 0 \text{ on } \partial M.
\end{array}\right.
\]
Here we have used $\lambda_1(M,g)$ $>$ $0$.

Define $S:$ $C^{2,\alpha}(M)$ $\rightarrow$ $C^{2,\alpha}(M)$ by $Su$ $=$ $u - (-L_g)^{-1}(u^{\frac{n+2}{n-2}})$. By Proposition \ref{LinearLeading} in Appendix \ref{App-Degree},
\begin{equation}
\deg((\tilde G_0,\tilde B_0),D_0,0) = \deg_{L.S.}(S\big|_{C^{4,\alpha}(M)},D_0,0),
\label{Exist::Eqn03}
\end{equation}
where $\deg_{L.S.}$ denotes the Leray-Schauder degree. 

Let
\begin{align*}
\hat D_0
	&= \Big\{u \in C^{2,\alpha}(M): u > 0, \lambda(A_{u^{\frac{4}{n-2}g}}) \in \Gamma_1,\\
		&\qquad\qquad \|u\|_{C^{2,\alpha}(M)} + \|u^{-1}\|_{C^{2,\alpha}(M)} \leq 2C, \dist(\lambda(A_{u^{\frac{4}{n-2}g}}),\partial\Gamma_t) > \epsilon\Big\}.
\end{align*}
In \cite[Section 5]{H-L-Duke}, it was shown that 
\begin{equation}
\deg_{L.S.}(S,\hat D_0, 0) = -1,
\label{Exist::Eqn04}
\end{equation}
On the other hand, by the reduction property of the Leray-Schauder degree (see e.g. \cite[Theorem 8.7]{Deimling}), 
\begin{equation}
\deg_{L.S.}(S,\hat D_0, 0) = \deg_{L.S.}(S|_{C^{4,\alpha}(M)},\hat D_0 \cap C^{4,\alpha}(M), 0).
\label{Exist::Eqn05}
\end{equation}
In addition by the excision property of the Leray-Schauder degree,
\begin{equation}
\deg_{L.S.}(S|_{C^{4,\alpha}(M)},\hat D_0 \cap C^{4,\alpha}(M), 0) = \deg_{L.S.}(S|_{C^{4,\alpha}(M)}, D_0, 0).
\label{Exist::Eqn06}
\end{equation}

Taking \eqref{Exist::Eqn01}-\eqref{Exist::Eqn06} altogether into account, we arrive at $\deg((G_1,B_1),D_1,0)$ $=$ $-1$. Putting $D$ $=$ $D_1$ and recalling the definition of $(G_1,B_1)$, we get \eqref{DegCount}, which completes the proof.
\eproof 
 
\medskip 
\noindent{\bf Proof of Proposition \ref{PositiveView-High-k}.} (a) If $M$ is the standard half sphere, the result is obvious (though the proof below applies to this case as well). We thus assume that $M$ is not simply connected. 

Let $\pi:$ $\SS^n_+$ $\rightarrow$ $M$ be the conformal covering map. Then
\[
\sign \lambda_1(M,g) = \sign \lambda_1(\SS^n_+,\pi^*(g)) = \sign \lambda_1(\SS^n_+,g_{\SS^n}),
\]
and so $\lambda_1(M,g)$ $>$ $0$. Here $\pi^*(g)$ denotes the pull-back metric of $g$ to $\SS^n_+$ by $\pi$.

By \cite[Theorem 0.1]{H-L-Duke}, there exists a metric $\hat g$ conformal to $g$ such that $R_{\hat g}$ $\equiv$ $1$ in $M$ and $h_{\hat g}$ $\equiv$ $c$ on $\partial M$. (In fact, in this case the proof for the existence for such metric is much simpler. One first argues as in Case 2 of the proof of Proposition \ref{UpperBound} to get an a priori $C^0$ bound. Higher derivatives estimates then follow, and so a degree theory argument can be carried out to conclude the desired existence.)

Since $\pi$ is conformal, $\pi^*(\hat g)$ satisfies $R_{\pi^*(\hat g)}$ $\equiv$ $1$ on $\SS^n_+$ and $h_{\pi^*(\hat g)}$ $\equiv$ $c$ on $\partial \SS^n_+$. By the Liouville theorem in \cite[Theorem 2.1(b)]{E3} (see also \cite{LiZhu}), $\pi^*(\hat g)$ can be obtained by $g_{\SS^n}$ by a conformal transformation on $\SS^n$. In particular, $\pi^*(\hat g)$ has sectional curvature $1$. It follows that $\hat g$ also has sectional curvature $1$, which establishes (a).

\medskip 
\noindent (b) It is readily seen that (i) $\Rightarrow$ (ii) is a consequence of (a). In addition, that (ii) $\Rightarrow$ (iii) is obvious. It remains to show (iii) $\Rightarrow$ (i). Assume that $M$ is of $\Gamma_k$-positive type for some $\frac{n}{2}$ $\leq$ $k$ $\leq$ $n$. By Theorem \ref{Existence} (see also \cite{J-L-L} and \cite{Chen}), there exists a conformal metric $\hat g$ such that $\sigma_k(\lambda(A_{\hat g}))$ $=$ $1$, $\lambda(A_{\hat g})$ $\in$ $\Gamma_k$ in $M$ and $h_{\hat g}$ $=$ $0$ on $\partial M$. Let $M_2$ be the double of $M$ and equip it with the metric induced by $\hat g$. Then $M_2$ is a $C^{2,1}$ closed locally conformally flat manifold with $\lambda(A_{\hat g})$ $\in$ $\Gamma_k$, and so $M_2$ is a quotient of the standard sphere $\SS^n$ by \cite[Corollary 1]{G-V-W}. In particular, there is a covering map $\pi:$ $\SS^n$ $\rightarrow$ $M_2$ such that $\pi^*(\hat g)$ is conformal to the round metric.

Notice that on $\SS^n$, $\pi^*(\hat g)$ satisfies $\sigma_k(\lambda(A_{\pi^*(\hat g)}))$ $=$ $1$ and $\lambda(A_{\pi^*(\hat g)})$ $\in$ $\Gamma_k$. Hence, by the Liouville theorem in \cite{V}, \cite{V-TrAMS} (see also \cite{L-L-CPAM}), $(\SS^n,\pi^*(\hat g))$ and the standard $\SS^n$ differ by a M\"{o}bius transformation. Therefore, by changing $\pi$ if necessary, we can assume without loss of generality that $\pi^*(\hat g)$ is the round metric. Then any component of $\pi^{-1}(\partial M)$ $\subset$ $\SS^n$ is an umbilic minimal hypersurface and so is an equator. This implies that $\pi^{-1}(\partial M)$ is connected and is an equator, which consequently implies that $\pi^{-1}(M)$ is conformally equivalent to the standard half-sphere $\SS^n_+$. 
\eproof

\section{Counterexamples to $C^2$ estimates}\label{Counterex}

As mentioned in the introduction, the local $C^2$ estimates of Jin, Li and Li stated in Theorem \ref{C2BdryEst} fail in general if one allows $\eta$ to attain a negative value. Moreover, such estimates also fail at global level when one has some (natural) additional ellipticity assumption on the ambient manifold. Here we discuss the setting in which we construct counterexamples and provide precise statements in Lemmas \ref{NegEx} and \ref{NegEx-Global}. More specifically, we use radial solutions on annuli. It should be noted that radial solutions of the $\sigma_k$ equations were systematically analyzed by Chang, Han and Yang in \cite{C-H-Y}.

We first state a simple fact from linear algebra.

\begin{lemma}\label{EValues}
If $M$ $=$ $\mu\,x \otimes x + \nu\,I$ then $M$ has exactly two real eigenvalues $\mu\,|x|^2 + \nu$, which is simple, and $\nu$, which is of order $n-1$.
\end{lemma}

For $R$ $>$ $1$, let $A_R$ denote the annulus
\[
A_R = \{x \in \RR^n: 1 \leq |x| \leq R\}.
\]
Let $u$ be a radial function, i.e. $u(x)$ $=$ $u(|x|)$. Following
the notations of \cite{C-H-Y}, let
\[
t = \ln |x| \text{ and } \xi(t) = -\frac{2}{n-2}\ln u(|x|) - \ln |x|.
\]
A straightforward calculation using Lemma \ref{EValues} shows that the eigenvalues of $A^u$ are
\begin{align*}
\lambda_1
    &= \frac{n-2}{2}\,e^{2\xi}[\xi_{tt} - (1 - \xi_t^2)],\\
\lambda_2
    &= \ldots = \lambda_n = \frac{n-2}{4}\,e^{2\xi}(1 - \xi_t^2).
\end{align*}
We thus have:

\begin{lemma}\label{TheODE}
If a function $\xi$ satisfies
\begin{equation}
\left\{\begin{array}{l}
e^{2k\xi}(1 - \xi_t^2)^{k-1}[\xi_{tt} + \frac{n-2k}{2k}(1-\xi_t^2)] = \frac{2^{k-1}}{C_{n-1}^{k-1}} =: \Theta, \qquad 0 \leq t \leq \ln R,\\
-1 < \xi_t < 1, \qquad 0 \leq t \leq \ln R,
\end{array}\right.
\label{ODE}
\end{equation}
where $C_{n-1}^{k-1}=\frac{(n-1)!}{(k-1)!(n-k)!}$, 
then the function $u$ defined by
\[
u(x) = u(|x|) = \exp\Big[-\frac{n-2}{2}\xi(\ln|x|) - \frac{n-2}{2}\ln |x|\Big]
\]
satisfies
\[
\left\{\begin{array}{l}
\sigma_k(\lambda(A^{u})) = 1 \text{ in } B_{R} \setminus B_1,\\
u > 0 \text{ and }\lambda(A^{u}) \in \Gamma_k \text{ in } \bar B_{R} \setminus B_1.
\end{array}\right.
\]
\end{lemma}

The following lemma produces a counterexample for boundary local $C^2$ estimate.

\begin{lemma}\label{NegEx}
For any $2$ $\leq$ $k$ $\leq$ $n$ and $c$ $<$ $0$, there exist $C_0$ $>$ $0$, $R_0$ $>$ $1$ and a family of $\{u_j\}$ $\subset$ $C^\infty(\bar B_{R_0}\setminus B_1)$ satisfying
\begin{equation}
\left\{\begin{array}{l}
\sigma_k(\lambda(A^{u_j})) = 1 \text{ in } B_{R_0} \setminus B_1,\\
u_j > 0 \text{ and }\lambda(A^{u_j}) \in \Gamma_k \text{ in } \bar B_{R_0} \setminus B_1,\\
\frac{\partial u_j}{\partial r} + \frac{n-2}{2}\,u_j = -\frac{n-2}{2}c\,u_j^{\frac{n}{n-2}} \text{ on } \partial B_1,\\
\end{array}\right.
\label{CircEqn}
\end{equation}
and
\[
|u_j| + |u_j^{-1}| + |\nabla u_j| \leq C_0 \text{ in } B_{R_0} \setminus B_1,
\]
such that
\[
\lim_{j \rightarrow \infty} \inf_{\partial B_1}|\nabla^2 u_j| = \infty.
\]
\end{lemma}

\bproof Fix $c<0$. We will consider small constants $0<\epsilon<\delta<\frac 12$ whose values will be specified later. For any such $\epsilon$ and $\delta$, there is clearly a smooth function $\xi(t)\equiv \xi(t; \epsilon, \delta)$
satisfying near $t=0$
\begin{equation}
e^{2k\xi}(1 - \xi_t^2)^{k-1}[\xi_{tt} + \frac{n-2k}{2k}(1-\xi_t^2)] =  \Theta,
\label{AA1-1}
\end{equation}
and
\begin{equation}
\xi(0)=\epsilon+\ln |c|,
\qquad
\xi_t(0)=-e^{-\epsilon}=ce^{-\xi(0)}.
\label{AA1-2}
\end{equation}
(Here we have used $1 - \xi_t(0)^2>0$.) We first require that $\delta$ is small enough so that the following argument goes through
(note that $0<\epsilon<\delta$)
\begin{eqnarray}
\xi_{tt}(0)
&=&  \Theta e^{-2k\xi(0)}
(1 - \xi_t(0)^2)^{1-k}-  \frac{n-2k}{2k}
(1 - \xi_t(0)^2)
\nonumber\\
&=&
 \Theta
 e^{-2k\ln|c|} [1+O(\epsilon)]
(2\epsilon+O(\epsilon^2))^{1-k}
-  \frac{n-2k}{2k}
(2\epsilon+O(\epsilon^2))\nonumber\\
&=&
\frac  \Theta
  { |c|^{2k} } (2\epsilon)^{1-k}
[1+O(\epsilon)]\ge \frac  \Theta
  { 2 |c|^{2k} } (2\epsilon)^{1-k}>0.\nonumber
\end{eqnarray}
Therefore $\xi_t$ is strictly increasing in $t$ for small $t$.

We will only consider those values of $\epsilon$ satisfying
$-e^{-\epsilon}<  -1+ \frac \delta 2$. It follows that for small positive $t$,
\begin{equation}
-e^{-\epsilon}< \xi_t(t) <-1+\delta,
\label{AA1-3}
\end{equation}
and
\begin{equation}
 \xi_{tt}(t)>0.
\label{AA1-4}
\end{equation}

For each such $\epsilon$, we let $(0, T(\epsilon, \delta))$, $0< T(\epsilon, \delta)\le \infty$, be the largest open interval on which \eqref{AA1-1}-\eqref{AA1-4} hold.

We note that, in $[0, T(\epsilon, \delta))$, one has
\begin{align}
&0<(1- (-e^{-\epsilon})^2)\le
(1 - \xi_t^2)\le  1-( -1+\delta)^2)\le 2\delta,
\label{AA1-4*1}\\
&\xi\le \xi(0)\le \ln|c|,
\label{AA1-4*2}
\end{align}
which implies in view of \eqref{AA1-1} that
\begin{equation}
\xi_{tt} \ge 
 \Theta e^{-2k \ln|c|}
(2\delta)^{1-k}- \frac{|n-2k|}{2k}(2\delta)
\ge  \frac  \Theta 2 e^{-2k \ln|c|}
(2\delta)^{1-k}>0
\label{AA1-9}
\end{equation}
for all $\delta$ sufficiently small. 

We now fix $\delta$. By \eqref{AA1-3}, \eqref{AA1-9} and the mean value theorem,
\[
\delta \geq \xi_t(t) - \xi_t(0) \geq \frac{\Theta}{2} e^{-2k \ln|c|} (2\delta)^{1-k}\,t \text{ for } 0 < t < T(\epsilon,\delta).
\]
It follows that $T(\epsilon,\delta)$ is finite and
\begin{equation}
T(\epsilon,\delta) \leq C(n,k,|c|,\delta),
\label{AA1-9bis}
\end{equation}
where here and below $C(n,k,|c|,\delta)$ denotes some positive constant independent of $\epsilon$.  Moreover, \eqref{AA1-4*1}-\eqref{AA1-9} hold in $[0,T(\epsilon,\delta)]$. 

We next show that
\begin{equation}
T(\epsilon,\delta) \geq \frac{1}{C(n,k,|c|,\delta)} > 0.
\label{AA1-6}
\end{equation}
Evidently, this implies
\[
\liminf_{\epsilon\to 0}T(\epsilon, \delta)>0,
\]
which proves the assertion by virtue of \eqref{AA1-3}, \eqref{AA1-4*2}, \eqref{AA1-9bis} and Lemma \ref{TheODE}.

If $ T(\epsilon, \delta)\ge 1$, \eqref{AA1-6} holds. Otherwise,
$ T(\epsilon, \delta)<1$ and, in view of \eqref{AA1-3},\eqref{AA1-9}, and the definition of $T(\epsilon,\delta)$,  
$ \xi_t(T(\epsilon, \delta))=
-1+\delta$. By  (\ref{AA1-9}) and 
$\xi_t(0)< -1+\frac \delta 2< -1+\delta
= \xi_t(T(\epsilon, \delta))
$, 
there exists some $\widehat T\in (0, T(\epsilon, \delta))$ such that
$$
 \xi_t(\widehat T)= 
-1 +\frac \delta 2, \quad -1 +\frac \delta 2\le
\xi_t\le -1 + \delta \
\mbox{on}\ 
[\widehat T, T(\epsilon, \delta)].
$$
It is then easy to see from (\ref{AA1-1}) that
$$
|\xi_{tt}|\le C(n,k,|c|,\delta)\ \mbox{on}\ 
[\widehat T, T(\epsilon, \delta)].
$$
It follows from the mean value theorem that
$$
\frac \delta 2=
| \xi_t(\widehat T)- \xi_t(T(\epsilon, \delta)|
\le C(n,k,|c|,\delta) ( T(\epsilon, \delta)- \widehat T)
\le  C(n,k,|c|,\delta)  T(\epsilon, \delta),
$$
from which we deduce (\ref{AA1-6}).
\eproof

The failure of $C^2$ estimates is more dramatic in the sense that on any fixed annulus, which is of $\Gamma_k$-positive type for any $1$ $\leq$ $k$ $<$ $\frac{n}{2}$, $C^2$ estimates also fail for sufficiently negative perspective boundary mean curvature values. We announce here the statement but defer the proof to a forthcoming paper.

\begin{lemma}[\cite{LiNg}]\label{NegEx-Global}
Let $B_1$ be the open lower half-sphere of the standard $\SS^n$ and $B_2$ be an open geodesic ball centered at the north pole of $\SS^n$ with radius $0$ $<$ $r$ $<$ $\frac{\pi}{2}$. Note that $B_1$ and $B_2$ are disjoint.

For any $2$ $\leq$ $k$ $\leq$ $n$, there exist $c_*$ $=$ $c_*(n,k,r)$ $<$ $0$, a sequence $c_j$ $>$ $c_*$, $c_j$ $\rightarrow$ $c_*$ and a family of $\{u_j\}$ $\subset$ $C^\infty(\SS^n\setminus (B_1 \cap B_2))$ satisfying
\begin{equation}
\left\{\begin{array}{l}
\sigma_k(\lambda(A_{u_j^\frac{4}{n-2}g})) = 1 \text{ in } \SS^n\setminus (B_1 \cap B_2),\\
\lambda(A_{u_j^\frac{4}{n-2}g}) \in \Gamma_k \text{ and } u_j > 0 \text{ in } \SS^n\setminus (B_1 \cap B_2),\\
h_{u_j^{\frac{4}{n-2}}g} = c_j \text{ on } \partial B_1,\\
h_{u_j^{\frac{4}{n-2}}g} = 0 \text{ on } \partial B_2,
\end{array}\right.
\label{CircEqn-Global}
\end{equation}
and
\[
|u_j| + |u_j^{-1}| + |\nabla u_j| \leq C_0 \text{ in } \SS^n\setminus (B_1 \cap B_2)
\]
such that
\[
\lim_{j \rightarrow \infty} \sup_{\SS^n\setminus (B_1 \cap B_2)}|\nabla^2 u_j| = \infty.
\]
Moreover, for fixed $n$ and $k$, $c_*$ is a continuous, strictly increasing function of $r$ and
\[
\lim_{r \rightarrow \frac{\pi}{2}} c_* = 0 \text{ and } \lim_{r \rightarrow 0} c_* = -\infty.
\]
\end{lemma}
 
\section{$\Gamma$-type of a manifold}\label{GammaType}

In this appendix, we briefly study the notion of $\Gamma$-types introduced in the introduction. Throughout this appendix, we assume that $(f,\Gamma)$ 
satisfies \eqref{ConvexCone}-\eqref{f-ellipticity} and 
\eqref{f-concavity}, unless otherwise stated. We first 
recall the notion of  viscosity solutions used in Definition \ref{Type}.
\begin{definition}\label{ViscoSol}
A positive continuous function $u$ in $M$ is a viscosity supersolution (respectively, subsolution) of
\[
\lambda(A_{u^{\frac{4}{n-2}}g}) \in \partial\Gamma \text{ in } M^\circ
\]
when the following holds: If $x_0$ $\in$ $M^\circ$, $\varphi$ $\in$ $C^2(M^\circ)$, $(u - \varphi)(x_0)$ $=$ $0$, and $u - \varphi$ $\leq$ $0$ near $x_0$, then
\[
\lambda(A_{\varphi^{\frac{4}{n-2}}g}(x_0)) \in \RR^n \setminus \Gamma
\]
(respectively, if $(u - \varphi)(x_0)$ $=$ $0$, and $u - \varphi$ $\geq$ $0$ near $x_0$, then $\lambda(A_{\varphi^{\frac{4}{n-2}}g}(x_0))$ $\in$ $\bar\Gamma$). We say that $u$ is a viscosity solution if it is both a 
viscosity supersolution and a viscosity subsolution.
\end{definition}

We next show that the
 $\Gamma$-positive type and the
 $\Gamma$-nonpositive type are mutually exclusive.
\begin{lemma}\label{MutuallyExclusive}
Let $(M,g)$ be a smooth compact Riemannian manifold with boundary
 and  $\Gamma$ satisfy \eqref{ConvexCone}
and \eqref{PositiveCone}. 
Then $(M,g)$ cannot be simultaneously of $\Gamma$-positive type and of $\Gamma$-nonpositive type.
\end{lemma}

\bproof It suffices to show that if $\lambda(A_g)$ $\in$ $\Gamma$ in $M$ and $h_g$ $\geq$ $0$ on $\partial M$ and if the equation
\[
\left\{\begin{array}{l}
\lambda(A_{u^{\frac{4}{n-2}}g}) \in \partial\Gamma \text{ in } M^\circ,\\
h_{u^{\frac{4}{n-2}}g} \leq 0 \text{ on } \partial M
\end{array}\right.
\]
has a $C^{0,1}$ viscosity solution $u$, then a contradiction must occur.

First, arguing as in the proof of Theorem \ref{Main-Thm}, we can assume that $h_g$ $>$ $0$ on $\partial M$. Next, pick $x_0$ $\in$ $M$ such that
\[
u(x_0) = \max_M u > 0.
\]
If $x_0$ $\in$ $M^\circ$, then using the definition of viscosity supersolution, we see for $\varphi$ $\equiv$ $u(x_0)$ that
$$
u(x_0)^{-\frac{4}{n-2}}\lambda(A_g(x_0))=
\lambda(A_{\varphi^{\frac{4}{n-2}}g}(x_0))\in 
\RR^n \setminus \Gamma,
$$
violating $\lambda(A_g(x_0))\in \Gamma$.

If $x_0$ $\in$ $\partial M$, then we must have 
$\frac{\partial u}{\partial \nu}(x_0)$ $\leq$ $0$
in the viscosity sense, which implies that
\[
0 < \frac{\partial u}{\partial \nu}(x_0) + \frac{n-2}{2}
\, h_g(x_0)\,u(x_0) = \frac{n-2}{2}\,h_{u^{\frac{4}{n-2}}g}(x_0) \leq 0.
\]
This is also impossible. Lemma \ref{MutuallyExclusive} is established.
\eproof

Next, we show:

\begin{lemma}\label{Exhaustive}
Let $(M,g)$ be a smooth compact locally conformally flat Riemannian manifold 
with umbilic boundary and  $(f,\Gamma)$ satisfy 
\eqref{ConvexCone}-\eqref{PositiveCone}. If $\lambda_1(M,g)$ $>$ $0$ and $\Gamma$ admits a smooth concave defining function, then $(M,g)$ is either of $\Gamma$-positive type or of $\Gamma$-nonpositive type.
\end{lemma}

Recall that a smooth concave defining function for an open set $G$ is a function $h$ $\in$ $C^\infty(G) \cap C^0(\bar G)$ which is concave and positive in $G$ and vanishes on $\partial G$.

\medskip
\bproof Assume that $(M,g)$ is not of $\Gamma$-nonpositive type. We will show that $(M,g)$ is of $\Gamma$-positive type.

First, we claim that Theorem \ref{Main-Thm} holds if we replace the hypothesis ``$(M,g)$ is of $\Gamma$-positive type'' by ``$(M,g)$ is \underline{not} of $\Gamma$-nonpositive type''. Indeed, by arguing as in the proof of Theorem \ref{Main-Thm} and using Proposition \ref{UpperBound}, it suffices to show that
\[
\max_M u \geq C(n,M,g,f,\Gamma,\beta) > 0
\]
for any positive solution $u$ of \eqref{Main-Eqn}. Assume otherwise that this is incorrect, so that there is a family $\{u_j\}$ of positive solutions of \eqref{Main-Eqn} satisfying
\[
\max_M u_j =: \epsilon_j \rightarrow 0.
\]
Define $v_j$ $=$ $\frac{u_j}{\epsilon_j}\le 1$. Then
\[
\left\{\begin{array}{l}
f(\lambda(A_{v_j^\frac{4}{n-2}g})) = \epsilon_j^{\frac{4}{n-2}} \text{ in }M^\circ,\\
\lambda(A_{v_j^\frac{4}{n-2}g}) \in \Gamma \text{ and } v_j > 0 \text{ in } M,\\
\frac{\partial v_j}{\partial\nu} + \frac{n-2}{2}h_g\,v_j = c_k\,\epsilon^{\frac{2}{n-2}}\,v_j^{\frac{n}{n-2}} \text{ on } N_k, k = 1 \ldots m.
\end{array}\right.
\]
Moreover, by Proposition \ref{UpperBound} and the gradient estimate in \cite[Theorem 1.19]{Li-CPAM}, we also have
\[
|\nabla\ln v_j| = |\nabla\ln u_j| \leq C(n,M,g,f,\Gamma,\beta).
\]
Hence, by passing to a subsequence if necessary, we get $v_j$ converges in H\"{o}lder norm to some function $v$ which is a positive $C^{0,1}$ viscosity solution of
\[
\left\{\begin{array}{l}
f(\lambda(A_{v^\frac{4}{n-2}g})) = 0 \text{ in }M^\circ,\\
\lambda(A_{v^\frac{4}{n-2}g}) \in \bar\Gamma \text{ and } v_j > 0 \text{ in } M,\\
\frac{\partial v}{\partial\nu} + \frac{n-2}{2}h_g\,v = 0 \text{ on } \partial M,
\end{array}\right.
\]
which contradicts the fact that $(M,g)$ is not of $\Gamma$-nonpositive type.

By the above claim, we see that Theorem \ref{Existence} also holds if we replace the hypothesis ``$(M,g)$ is of $\Gamma$-positive type'' by ``$(M,g)$ is not of $\Gamma$-nonpositive type''. In particular, if $(M,g)$ is not of $\Gamma$-nonpositive type, and there exists a function $f$ on $\Gamma$ which satisfies \eqref{f-smoothness}-\eqref{f-concavity}, then the problem \eqref{Main-Eqn} has a solution for $c_k$ $=$ $0$, which implies that $(M,g)$ is of $\Gamma$-positive type. The result then follows from Lemma \ref{ConeFunc} below.
\eproof

\begin{lemma}\label{ConeFunc}
Let $\Gamma$ $\subset$ $\RR^n$ be a cone satisfying \eqref{ConvexCone} and \eqref{PositiveCone}. There exists a function $f$ satisfying \eqref{f-smoothness}-\eqref{f-concavity} if and only if $\Gamma$ admits a concave defining function $h$ $\in$ $C^\infty(\Gamma) \cap C(\bar \Gamma)$.
\end{lemma}

\bproof The necessity is clear. For the sufficiency, assume that $\Gamma$ admits a defining function $h$ $\in$ $C^\infty(\Gamma) \cap C(\bar \Gamma)$. By considering 
\[
\tilde h(\lambda) = \sum_{x \text{ is a permutation of } \lambda} h(x),
\]
instead of $h$, we can assume without loss of generality that $h$ is symmetric.

Let $\Omega_\Gamma$ $=$ $\Gamma \cap \{\lambda: \lambda_1 + \ldots + \lambda_n = 1\}$, which is open, symmetric and convex. Let $\nabla_T$ denote the gradient on $\Omega_\Gamma$. Observe that for $x$ $\in$ $\Omega_\Gamma$ and $p_0$ $\in$ $\bar\Omega_\Gamma$, the concavity of $h$ implies that
\[
-\nabla_T h(x) \cdot (x - p_0) \geq - h(x) + h(p_0) \geq - h(x).
\]
Therefore, for $0$ $<$ $\alpha$ $<$ $1$ and $g$ $=$ $h^\alpha$, there holds
\begin{multline}
g(x) - \nabla_T g(x) \cdot (x - p_0) = h(x)^{\alpha-1}\big[h(x) - \alpha\,\nabla_T h(x) \cdot (x - p_0)\big]\\
 \geq (1-\alpha)h(x)^\alpha > 0 \text{ for any } x \in \Omega_\Gamma \text{ and } p_0 \in \bar\Omega_\Gamma.
\label{EllConstr}
\end{multline}

Define $f$ by
\[
f(\lambda) = (\lambda_1 + \ldots + \lambda_n)\,g\big(\frac{\lambda}{\lambda_1 + \ldots + \lambda_n}\big).
\]
We only need to show that $\partial_i f$ $>$ $0$, $\sum_i \partial_i f$ $>$ $\delta$ and $f$ is concave in $\Gamma$. 

For simplicity, we write
\[
[\lambda] = \lambda_1 + \ldots \lambda_n \text{ and }\lambda' = \frac{\lambda}{[\lambda]}.
\]
We compute
\begin{multline*}
\partial_i f(\lambda)
	= g(\lambda') + [\lambda]\,\partial_j g(\lambda')\,\frac{\delta_{ij}[\lambda] - \lambda_j}{[\lambda]^2}\\
		= g(\lambda') + \partial_i g(\lambda') - \partial_j g(\lambda')\lambda'_j
			= g(\lambda') - \nabla_T g(\lambda') \cdot(\lambda' - p^i),
\end{multline*}
where $p^i_j$ $=$ $\delta^i_j$. Since $\Gamma$ $\supset$ $\Gamma_n$, it follows that $p^i$ $\in$ $\bar\Omega_\Gamma$. Hence, by \eqref{EllConstr},
\[
\partial_i f(\lambda) \geq (1-\alpha)[\lambda]^{-1}\,f(\lambda) > 0 \text{ in }\Gamma.
\]

To prove the concavity of $f$, we calculate its Hessian. We have
\begin{align*}
[\lambda]\partial_{ij} f(\lambda)
	&= \partial_k g(\lambda')\frac{\delta_{kj}[\lambda] - \lambda_k}{[\lambda]} + \partial_{ki} g(\lambda')\frac{\delta_{kj}[\lambda] - \lambda_k}{[\lambda]}\\
		&\qquad\qquad - \partial_l g(\lambda')\frac{\delta_{lj}[\lambda] - \lambda_l}{[\lambda]} - \partial_{kl} g(\lambda')\,\lambda'_l\,\frac{\delta_{kj}[\lambda] - \lambda_k}{[\lambda]}\\
	&= \partial_{ij} g(\lambda') - \partial_{ki} g(\lambda')\lambda'_k - \partial_{lj} g(\lambda')\lambda'_l + \partial_{kl} g(\lambda')\lambda'_k\,\lambda'_l.
\end{align*}
Hence, for any $p$ $\in$ $\RR^n$, we have
\begin{align}
[\lambda]\partial_{ij} f(\lambda)\,p_i\,p_j
	&= \partial_{ij} g(\lambda')\,p_i\,p_j - \partial_{ki} g(\lambda')\lambda'_k\,p_i\,p_j - \partial_{lj} g(\lambda')\lambda'_l\,p_i\,p_j + \partial_{kl} g(\lambda')\lambda'_k\,\lambda'_l\,p_i\,p_j\nonumber\\
	&= \partial_{ij} g(\lambda')\,p_i\,p_j - 2\partial_{ki} g(\lambda')\lambda'_k\,p_i\,[p] + \partial_{kl} g(\lambda')\lambda'_k\,\lambda'_l\,[p]^2.
\label{Hessian--}
\end{align}
On the other hand, since $\nabla^2 g$ $\leq$ $0$ in $\Gamma$,
\begin{align*}
&\partial_{kl} g(\lambda')\lambda'_k\,\lambda'_l \leq 0,\\
&\big(\partial_{ki} g(\lambda')\lambda'_k\,p_i\big)^2 - \big(\partial_{ij} g(\lambda')p_i\,p_j\big)\big(\partial_{kl} g(\lambda')\lambda'_k\,\lambda'_l\big) \leq 0.
\end{align*}
Also, as $\Gamma$ $\subset$ $\Gamma_1$, $[\lambda]$ $>$ $0$ in $\Gamma$. Therefore, \eqref{Hessian--} implies that $\nabla^2 f$ $\leq$ $0$ in $\Gamma$, i.e. $f$ is concave in $\Gamma$.

Finally, we show \eqref{delta}. For $\mu$ $>$ $0$, since $\Gamma$ $\subset$ $\Gamma_n$, the concavity and the homogeneity of $f$ implies that
\begin{align*}
\mu\,f\big(\frac{1}{n} + \frac{\lambda_1}{\mu}, \ldots, \frac{1}{n} + \frac{\lambda_n}{\mu}\big) 
	&= f\big(\frac{\mu}{n} + \lambda_1, \ldots, \frac{\mu}{n} + \lambda_n\big) \\
	&\leq f(\lambda_1, \ldots, \lambda_n) + \frac{\mu}{n}\sum_{i=1}^n \partial_i f(\lambda).
\end{align*}
Dividing both side by $\mu$ and let $\mu$ $\rightarrow$ $\infty$, we get
\[
\sum_{i=1}^n \partial_i f(\lambda) \geq nf\big(\frac{1}{n},\ldots,\frac{1}{n}\big) = n\,h\big(\frac{1}{n},\ldots,\frac{1}{n}\big)^{\frac{1}{\alpha}} =:\delta.
\]
The proof is complete.
\eproof

\bigskip
\noindent{\bf Proof of Example \ref{ExampleOfGammaType}.} First observe that $M$ $=$ $\SS^n \setminus (B_1 \cup B_2)$ is conformally equivalent to some annulus $A_R$ $=$ $B_R \setminus B_1$. Consider the cylindrical metric $g_{cyl}$ $=$ $|x|^{-2}g_{flat}$ on $A_R$. We have
\[
|h_{g_cyl}| = (r^{-\frac{n-2}{2}})^{\frac{n}{n-2}}\big|\partial_r r^{-\frac{n-2}{2}} + \frac{n-2}{2r}\,r^{-\frac{n-2}{2}}\big| = 0 \text{ for } r = 1 \text{ or } r = R.
\]
Hence $\partial A_R$ is mean curvature free with respect to $g_{cyl}$. By a direct calculation, the eigenvalues of $A_{g_{cyl}}$ are found to be $-\frac{n-2}{2}$ with multiplicity $1$ and $\frac{n-2}{4}$ with multiplicity $n-1$. Hence
\begin{align*}
\sigma_k(\lambda(A_{g_{cyl}}))
	&= \frac{(n-2)^k}{4^k}\frac{(n-1)!}{(k-1)!(n-k)!}\big[-2 + \frac{n}{k}\big] > 0 \text{ if } n > 2k.
\end{align*}
By well-known properties of $\Gamma_k$, this shows that $\lambda(A_{g_{cyl}})$ $\in$ $\Gamma_k$ in $M$ for any $1$ $\leq$ $k$ $<$ $\frac{n}{2}$. The first assertion follows.

For the second part, observe that $M$ is not covered by the standard half-sphere and invoke Proposition \ref{PositiveView-High-k}(b).
\eproof

\appendix
\section{Conformally invariant boundary differential operators on Euclidean domains}\label{BdryOpClassProof}

\noindent{\bf Proof of Theorem \ref{BdryOpClass}.} First, applying \eqref{BdryDiffOp-ConfInv} for $\psi$ being a translation of $\RR^n$, we infer that
\begin{equation}
B(0,s,p,\nu,H) = B(x,s,p,\nu,H).
\label{ConfInvBdryDiffOp::x-indep}
\end{equation}
Thus, for simplicity, we will write $B(s,p,\nu,H)$ in place of $B(x,s,p,\nu,H)$.

\medskip
\noindent (a) Fix some $s$, $p$, $\nu$ and $H$. Let $u$ be a smooth function such that $u(x)$ $=$ $s$ and $\nabla u(x)$ $=$ $p$.

Assume for the moment that $p \cdot \nu$ $\neq$ $0$. Set
\[
\lambda = - \frac{(n-2)s}{p \cdot \nu} \text{ and }
    x = \lambda\nu.
\]
Define a conformal transformation $\psi$ by $\psi(z)$ $=$ $\lambda^2\,\frac{z}{|z|^2}$. Writing $r$ for $|z|$, we calculate
\begin{align*}
u_\psi
    &= \frac{|\lambda|^{n-2}}{r^{n-2}}\,u \circ \psi,\\
\nabla_i u_\psi
    &= \frac{|\lambda|^n}{r^{n+2}}(\delta_{ij} r^2 - 2 z_i\,z_j)\,\nabla_j u\circ \psi - (n-2)\frac{|\lambda|^{n-2}}{r^n}\,z_i\,u \circ \psi.
\end{align*}
As $\psi(x)$ $=$ $x$ and $|x|$ $=$ $|\lambda|$, it follows that
\begin{align*}
u_\psi(x)
	&= s,\\
\nabla_i u_\psi(x)
	&= (\delta_{ij} - 2 \nu_i\,\nu_j)\,p_j - (n-2)\frac{1}{\lambda}\,\nu_i\,s\\
	&= p_i - \nu_i\Big(2p_j\nu_j + \frac{(n-2)s}{\lambda}\Big) = p_i - (p \cdot \nu)\nu_i.
\end{align*}
Also,
\[
\nu_\psi(x)
    = -\nu \text{ and } H_\psi(x) = H - \frac{2}{\lambda}I = H + \frac{2 p \cdot \nu}{(n-2)s}I.
\]
Hence, by \eqref{BdryDiffOp-ConfInv},
\[
B(s,p- (p \cdot \nu)\nu,\nu,H) = B\Big(s,p,-\nu,H + \frac{2 p \cdot \nu}{(n-2)s}I\Big).
\]
Therefore,
\[
B(s,p,\nu,H) = B\Big(s,p - (p \cdot \nu)\nu, - \nu, H + \frac{2 p \cdot \nu}{(n-2)s}I\Big) \text{ for any } s, p, \nu \text{ with } p \cdot \nu \neq 0.
\]
By the continuity of $B$, the restriction that $p \cdot \nu$ $\neq$ $0$ can be dropped. We thus get
\begin{equation}
B(s,p,\nu,H) = B\Big(s,p - (p \cdot \nu)\nu, - \nu, H + \frac{2 p \cdot \nu}{(n-2)s}I\Big) \text{ for any } s, p, \nu.
\label{ConfInvBdryDiffOp::ReductionToPerp}
\end{equation}

Assume next that $p \cdot \nu$ $=$ $0$. This time we define
\[
\lambda = - \frac{(n-2)s}{|p|} \text{ and } x = \lambda\frac{p}{|p|}.
\]
For $\psi(z)$ $=$ $\lambda^2\,\frac{z}{|z|^2}$, a direct calculation using $p \cdot \nu$ $=$ $0$ shows that
\[
u_\psi(x) = s,  \nabla u_\psi(x) = 0,  \nu_\psi(x) = \nu,  \text{ and } H_\psi(x) = H.
\]
Therefore, by \eqref{BdryDiffOp-ConfInv},
\[
B(s,0,\nu,H) = B(s,p,\nu,H) \text{ provided } p \cdot \nu = 0.
\]
Combining with \eqref{ConfInvBdryDiffOp::ReductionToPerp}, we get
\begin{equation}
B(s,p,\nu,H) = B\Big(s,0, - \nu, H + \frac{2 p \cdot \nu}{(n-2)s}I\Big) \text{ for any } s, p, \nu.
\label{ConfInvBdryDiffOp::PRemoval}
\end{equation}

\medskip
\noindent (b) Fix $s$, $\nu$ and $H$. Pick a smooth function $u$ such that $u(0)$ $=$ $s$ and $\nabla u(0)$ $=$ $0$.

Applying \eqref{BdryDiffOp-ConfInv} for $\psi(z)$ $=$ $Rz$, we get
\[
B(R^{\frac{n-2}{2}}s,0,\nu,H) = B(s,0,\nu,R^{-1}H).
\]
Hence
\begin{equation}
B(s,0,\nu,H) = B(1,0,\nu, s^{-\frac{2}{n-2}}H). 
\label{ConfInvBdryDiffOp::SRemoval}
\end{equation}

Next, let $e$ be any unit vector. Pick an orthonormal matrix $O$ such that $O\nu$ $=$ $e$. Applying \eqref{BdryDiffOp-ConfInv} for $\psi(z)$ $=$ $Oz$ we get
\begin{equation}
B(s,0,\nu,H) = B(s,0,e,H).
\label{ConfInvBdryDiffOp::NuRemoval}
\end{equation}

\medskip
\noindent (c) Finally, combining \eqref{ConfInvBdryDiffOp::x-indep}, \eqref{ConfInvBdryDiffOp::PRemoval}, \eqref{ConfInvBdryDiffOp::SRemoval}, and \eqref{ConfInvBdryDiffOp::NuRemoval}, we conclude that
\[
B(x,s,p,\nu,H) = B\Big(0,1,0,e,s^{-\frac{2}{n-2}}\big(H + \frac{2}{n-2}\frac{p \cdot \nu}{s}I\big)\Big).
\]
The proof is complete.
\eproof

\section{Degree theory for second order elliptic operators with Neumann boundary conditions}\label{App-Degree}

In this appendix, we give a modification of \cite{L-CPDE} to define a degree theory for (nonlinear) second order elliptic equations with (nonlinear) oblique derivative boundary conditions. Let $\Omega$ $\subset$ $\RR^n$ be a domain with smooth boundary. For a fixed $f$ $\in$ $C^{3,\alpha}(\bar\Omega \times \RR \times \RR^n \times \textrm{Sym}^{n \times n})$, define a differential operator $F:$ $C^{4,\alpha}(\bar\Omega)$ $\rightarrow$ $C^{2,\alpha}(\bar\Omega)$ by
\[
F[u] = f(\cdot, u, \nabla u, \nabla^2 u).
\]
We say that $F$ is elliptic on some bounded open subset $\mathcal{O}$ of $C^{4,\alpha}(\bar\Omega)$ if for any $u$ $\in$ $\mathcal{O}$, $x$ $\in$ $\bar\Omega$ and $\xi$ $\in$ $\RR^n$ there holds,
\begin{equation}
\frac{\partial f}{\partial u_{ij}}(x,u,\nabla u,\nabla^2 u)\,\xi_i\,\xi_j > \rho|\xi|^2 \text{ for some } \rho > 0.
\label{FNL::Ellipticity}
\end{equation}

For $\beta_1$, \ldots, $\beta_n$, $\gamma$ $\in$ $C^{4,\alpha}(\partial\Omega \times \RR)$, define a boundary condition operator $B:$ $C^{4,\alpha}(\bar\Omega)$ $\rightarrow$ $C^{3,\alpha}(\partial\Omega)$ by
\[
B[u] = \Big(\beta_i(x, u)\,u_i + \gamma(x, u)\Big)\Big|_{\partial\Omega},
\]
where $\nu$ denotes the outer unit normal to $\partial\Omega$. We say that $B$ is oblique on some bounded open subset $\mathcal{O}$ of $C^{4,\alpha}(\bar\Omega)$ if for any $u$ $\in$ $\mathcal{O}$ and $x$ $\in$ $\partial\Omega$ there holds
\begin{equation}
\beta_i(x,u)\,\nu_i > \mu \text{ for some } \mu > 0.
\label{BC::Obliqueness}
\end{equation}

Let $\mathcal{O}$ $\subset$ $C^{4,\alpha}(\bar\Omega)$ be a bounded open set such that $\partial\mathcal{O} \cap (F,B)^{-1}(0)$ $=$ $\emptyset$, $F$ is elliptic on $\mathcal{O}$ and $B$ is oblique on $\mathcal{O}$.
We will define an integer-valued degree for $(F,B)$ on $\mathcal{O}$ at $0$ along the line of \cite{L-CPDE}.

Define
\begin{equation}
\left.\begin{array}{rl}
S: C^{2,\alpha}(\bar\Omega)
	&\rightarrow C^\alpha(\bar\Omega) \times C^{1,\alpha}(\partial\Omega)\\
u
	&\mapsto \big(S^{(1)}[u], S^{(2)}[u]\big) := \Big(\Delta u, \big(\frac{\partial u}{\partial\nu} + u\big)\big|_{\partial\Omega}\Big).
\end{array}\right.
\label{SDef}
\end{equation}
It is well-known that $S$ is an isomorphism.

Consider
\begin{align*}
T: C^{4,\alpha}(\bar\Omega)
	&\rightarrow C^\alpha(\bar\Omega) \times C^{1,\alpha}(\partial\Omega) \times C^{3,\alpha}(\partial\Omega)\\
u
	&\mapsto \big(T^{(1)}[u], T^{(2)}[u], T^{(3)}[u]\big) := \Big(S\circ F[u], B[u]\Big).
\end{align*}
As in \cite{L-CPDE}, we write
\begin{align*}
T^{(1)}[u]
	&= a_{st}(x, u, \nabla u, \nabla^2 u)\,u_{iist} + C_*(x,u,\nabla u, \nabla^2 u, \nabla^3 u),\\
T^{(2)}[u]
	&= \Big(a_{st}(x, u, \nabla u, \nabla^2 u)\,u_{sti}\,\nu_i  + \tilde f(x, u, \nabla u, \nabla^2 u)\Big)\Big|_{\partial\Omega},\\
T^{(3)}[u]
	&= \Big(\beta_i(x,u)\,u_i + \gamma(x,u)\Big)\Big|_{\partial\Omega},
\end{align*}
where
\begin{align*}
a_{st}(x, u, \nabla u, \nabla^2 u)
	&= \frac{\partial f}{\partial u_{st}}(x, u, \nabla u, \nabla^2 u).
\end{align*}
We split
\[
T[u] = L_{u,N}\,u + R_{u,N}[u]
\]
where $N$ is to be determined and
\begin{align*}
L_{u,N}^{(1)}\,w
	&= a_{st}(x, u, \nabla u, \nabla^2 u)\,w_{iist} - N\,a_{st}\,w_{st},\\
L_{u,N}^{(2)}\,w
	&= \Big(a_{st}(x, u, \nabla u, \nabla^2 u)\,w_{sti}\,\nu_i\Big)\Big|_{\partial\Omega},\\
L_{u,N}^{(3)}\,w
	&= \Big(\beta_i(x,u)\,w_i + w\Big)\Big|_{\partial\Omega},\\
R_{u,N}^{(1)}[w]
	&= N\,a_{st}\,w_{st} + C_*(x, u, \nabla u, \nabla^2 u, \nabla^3 u),\\
R_{u,N}^{(2)}[w]
	&= \tilde f(x, u, \nabla u, \nabla^2 u)\big|_{\partial\Omega},\\
R_{u,N}^{(3)}[w]
	&= \Big(-w + \gamma(x,u)\Big)\Big|_{\partial\Omega}.
\end{align*}

We claim that, there exists $N_0$ $\geq$ $0$ depending only on $\|a_{st}\|_{C^{1,\alpha}}$, $\|\beta_i\|_{C^{3,\alpha}}$, $\|\gamma\|_{C^{3,\alpha}}$, the ellipticity constant $\rho$ and the obliqueness constant $\mu$ such that for all $N$ $\geq$ $N_0$
\begin{equation}
L_{u,N}: C^{4,\alpha}(\bar\Omega) \rightarrow C^\alpha(\bar\Omega) \times C^{1,\alpha}(\partial\Omega) \times C^{3,\alpha}(\partial\Omega) \text{ is an isomorphism}.
\label{LInvertability}
\end{equation}
To see this, consider $L_0:$ $C^{4,\alpha}(\bar\Omega)$ $\rightarrow$ $C^\alpha(\bar\Omega) \times C^{1,\alpha}(\partial\Omega) \times C^{3,\alpha}(\partial\Omega)$ defined by
\[
Lw = \Big(\Delta^2 w, \big(\frac{\partial(\Delta w)}{\partial \nu} + \Delta w\big)\big|_{\partial\Omega}, \big(\frac{\partial w}{\partial \nu} + w\big)\big|_{\partial\Omega}\Big).
\]
By \cite[Theorem 12.1]{ADN-1}, $t\,L_0 + (1-t)L_{u,N}$ has finite dimensional kernel and so is Fredholm for $t$ $\in$ $(0,1)$. (Here we have used $a_{st} \in C^{1,\alpha}$, $\beta_i \in C^{3,\alpha}$, $\gamma \in C^{3,\alpha}$.) Moreover, by the stability of the Fredholm index (see e.g. \cite[Theorem 5.17]{Kato}), the Fredholm index of $L_{u,N}$ is the same as that of $L_0$, which is zero. Thus, to show that $L_{u,N}$ is isomorphic, it suffices to show that it has a trivial kernel. To this end, assume that $w$ $\in$ $\textrm{Ker}\,L_{u,N}$. First, using $H^2$ and $H^3$ estimates for linear elliptic equation of second order with oblique derivative boundary condition (see e.g. \cite[Theorem 15.2]{ADN-1}) and $L_{u,N}^{(3)}\,w$ $=$ $0$, we have
\begin{align}
\|w\|_{H^2}
	&\leq C\Big[\|a_{st}w_{st}\|_{L^2(\Omega)} + \|w\|_{L^2(\Omega)}\Big],
\label{PreH2Est}\\
\|w\|_{H^3}
	&\leq C\Big[\|a_{st}w_{st}\|_{H^1(\Omega)} + \|w\|_{L^2(\Omega)}\Big],
\label{PreH3Est}
\end{align}
where $C$ depends on $\|a_{st}\|_{C^{1,\alpha}}$, $\|\beta_i\|_{C^{2,\alpha}}$, $\|\gamma\|_{C^{2,\alpha}}$, $\rho$ and $\mu$. On the other hand, as the problem $(a_{st}\,\phi_{st},L_{u,N}^{(3)}\phi)$ $=$ $0$ has a unique solution for $N$ sufficiently large (due to \eqref{BC::Obliqueness} and the maximum principle),  the $\|w\|_{L^2}$ terms on the right hand sides of \eqref{PreH2Est} and \eqref{PreH3Est} can be dropped so that
\begin{align}
\|w\|_{H^2} \leq C\,\|a_{st}w_{st}\|_{L^2(\Omega)} \text{ and }
\|w\|_{H^3} \leq C\,\|a_{st}w_{st}\|_{H^1(\Omega)}.
\label{H2H3Est}
\end{align}
Next, using $L_{u,N}\,w$ $=$ $0$, we compute
\begin{align*}
0
	&=\int_{\Omega} L_{u,N}^{(1)}\,w\big[-a_{kl}\,w_{kl}\big]\,dx\\
	&\geq c_1\int_{\Omega} \big|\nabla(a_{st}w_{st})\big|^2 + N\int_{\Omega}\big|a_{kl}\,w_{kl}\big|^2\,dx - c_2\,\|w\|_{H^3(\Omega)}\,\|w\|_{H^2(\Omega)}\\
	&\geq \frac{c_1}{2}\int_{\Omega} \big|\nabla(a_{st}w_{st})\big|^2 + (N-c_3)\int_{\Omega}\big|a_{kl}\,w_{kl}\big|^2\,dx.
\end{align*}
In the above, the first inequality follows from integration by parts, $L_{u,N}^{(2)}w$ $=$ $0$, and H\"{o}lder's inequality, while the second inequality follows from Cauchy-Schwarz's inequality and \eqref{H2H3Est}. It is evident that for $N$ sufficiently large, we must have $a_{kl}\,w_{kl}$ $\equiv$ $0$. As $L_{u,N}^{(3)}\,w$ $=$ $0$, we must have $w$ $\equiv$ $0$, i.e. $L_{u,N}$ is injective. The claim follows.

Note that, by \cite[Theorem 7.3]{ADN-1}, $(L_{u,N})^{-1}$ maps $C^{1,\alpha}(\bar\Omega) \times C^{2,\alpha}(\partial\Omega) \times C^{4,\alpha}(\partial\Omega)$ into $C^{5,\alpha}(\bar\Omega)$, and its norm as a linear map between these spaces is bounded by a constant that depends only on $\|a_{st}\|_{C^{2,\alpha}}$, $\|\beta_i\|_{C^{4,\alpha}}$, $\|\gamma\|_{C^{4,\alpha}}$, $\rho$ and $\mu$. It follows that $u$ $\mapsto$ $V_N[u]$ $:=$ $(L_{u,N})^{-1}\circ R_{u,N}[u]$ is a compact operator from $\mathcal{O}$ to $C^{4,\alpha}(\bar\Omega)$. Moreover, $\partial\mathcal{O} \cap (Id + V_N)^{-1}(0)$ $=$ $\partial\mathcal{O} \cap (F,B)^{-1}(0)$ $=$ $\emptyset$. Therefore, we can define
\begin{equation}
\deg((F,B),\mathcal{O},0) = \deg_{L.S.}(Id + V_N, \mathcal{O},0), \qquad N \geq N_0,
\label{DegDef}
\end{equation}
where $\deg_{L.S.}$ denotes the Leray-Schauder degree. As in \cite{L-CPDE}, it can be shown that this definition of the degree is independent of $N$.

It is also standard to check that the degree constructed above satisfies the following standard requirements.
\begin{enumerate}[(a)]
\item If $\deg((F,B),\mathcal{O},0)$ $\neq$ $0$, then there exists $u$ $\in$ $\mathcal{O}$ such that $(F[u],B[u])$ $=$ $0$.
\item The excision property: If $~\mathcal{\bar U}$ $\subset$ $\mathcal{O}$ and $~\mathcal{\bar U}\, \cap\, (F,B)^{-1}(0)$ $=$ $\emptyset$, then $\deg((F,B),\mathcal{O},0)$ $=$ $\deg((F,B),\mathcal{O} \setminus \mathcal{\bar U},0)$.
\item The homotopy invariance property: If $t$ $\mapsto$ $(f_t,b_t)$ is continuous from $[0,1]$ to $C^{3,\alpha}(\bar\Omega \times \RR \times \RR^n \times \textrm{Sym}^{n \times n}) \times C^{4,\alpha}(\partial\Omega \times \RR \times \RR^n)$, $F_t$ is elliptic on $~\mathcal{O}$, $B_t$ is oblique on $~\mathcal{O}$, and $\partial\mathcal{O} \cap (F_t,B_t)^{-1}(0)$ $=$ $\emptyset$, then $\deg((F_t,B_t),\mathcal{O},0)$ is independent of $t$.
\end{enumerate}

We next consider the ``compatability'' of the degree defined by \eqref{DegDef} with the Leray-Schauder degree in case $(F,B)$ has linear leading terms. We begin with a result about the semi-finiteness of a linear operator, which is known to the experts but relatively hard to find in the literature.

\begin{lemma}\label{SemiFiniteness}
Assume $a_{ij}$ $\in$ $C^1(\bar\Omega)$, $b_i$, $c$ $\in$ $C^0(\bar\Omega)$, $\beta_i$ $\in$ $C^1(\partial\Omega)$ and $\gamma$ $\in$ $C^0(\partial\Omega)$. Assume furthermore that $(a_{ij})$ $>$ $\rho\,I$ in $\bar\Omega$ and $\beta \cdot \nu$ $>$ $\mu$ on $\partial\Omega$ for some $\rho$ $>$ $0$ and $\mu$ $>$ $0$. There exists $\lambda_*$ depending on  $\|a_{ij}\|_{C^1(\bar\Omega)}$, $\|b_i\|_{C^0(\bar\Omega)}$, $\|c^+\|_{C^0(\bar\Omega)}$, $\|\beta_i\|_{C^1(\partial\Omega)}$, $\|\gamma^-\|_{C^0(\partial\Omega)}$, $\rho$ and $\mu$ such that for any $\lambda$ $>$ $\lambda_*$, the problem 
\begin{equation}
\left\{\begin{array}{l}
a_{ij}(x)\,u_{ij} + b_i(x)\,u_i + c(x)\,u = \lambda u \text{ in } \Omega,\\
\beta_i(x)\,u_i + \gamma(x)\,u = 0 \text{ on } \partial\Omega
\end{array}\right.
\label{LinearEVP}
\end{equation}
has no non-trivial solution in $C^2(\bar\Omega)$.
\end{lemma}

\begin{remark}
If $\gamma$ $\geq$ $0$ on $\partial\Omega$, the result follows directly from the maximum principle. In fact, $\lambda_*$ can then be taken to be $\|c^+\|_{C^0(\bar\Omega)}$.
\end{remark}

\bproof Assume that $u$ $\in$ $C^2(\bar\Omega)$ is a solution to \eqref{LinearEVP}. Let $\nu$ be the outer normal to $\partial \Omega$. Split $\beta$ $=$ $(\beta \cdot \nu)\,\nu + \beta_T$. Then by the second equation in \eqref{LinearEVP},
\begin{equation}
\frac{\partial u}{\partial \nu} = -\frac{1}{\beta \cdot \nu}\big(\beta_T \cdot \nabla u + \gamma\,u) \text{ on } \partial\Omega.
\label{LSF::NormalDer}
\end{equation}

Multiplying the first equation in \eqref{LinearEVP} by $u$ then integrating over $\Omega$, we get
\begin{align*}
\lambda\int_\Omega u^2\,dx
	&\leq -c_1\int_\Omega |\nabla u|^2\,dx + c_2\int_\Omega u^2\,dx + \int_{\partial\Omega} u\,a_{ij}\,u_i\,\nu_j\,d\sigma(x),
\end{align*}
where $c_1$ and $c_2$ are some positive constant that depends only on $\|a_{ij}\|_{C^{1}}$, $\|b_i\|_{C^0}$, $\|c^+\|_{C^0}$ and $\rho$. To proceed, we write $a_{ij}\nu_j$ $=$ $p\nu_i + q_i$ with $q \cdot \nu$ $=$ $0$. It is easily seen that $p$ $=$ $a_{ij}\nu_i\,\nu_j$ $>$ $0$ on $\partial\Omega$. Then, by \eqref{LSF::NormalDer}, 
\begin{align}
c_1\int_\Omega |\nabla u|^2\,dx + (\lambda - c_2)\int_\Omega u^2\,dx
	&\leq \int_{\partial\Omega} u\,\Big[p\,\frac{\partial u}{\partial \nu} + q_i\,u_i\Big]\,d\sigma(x)\nonumber\\
	&= -\int_{\partial \Omega} \frac{p}{\beta \cdot \nu}\,\gamma\,u^2\,d\sigma(x) + \int_{\partial\Omega} u\,X(u)\,d\sigma(x),\nonumber\\
	&\leq c_3\int_{\partial \Omega} u^2\,d\sigma(x) + \int_{\partial\Omega} u\,X(u)\,d\sigma(x),\label{LSF::Bound1}
\end{align}
where $c_3$ is a constant that depends only on $\|a_{ij}\|_{C^0}$, $\|\gamma^-\|_{C^0}$ and $\mu$, and $X$ $=$ $X_i\,\partial_i$ is some vector field along $\partial\Omega$ which is tangential to $\partial\Omega$ and $\|X\|_{C^1(\partial\Omega)}$ is bounded from above by a constant depending only on $\|a_{ij}\|_{C^1}$, $\|\beta\|_{C^1}$ and $\mu$.

On the other hand, by the smoothness of $\partial\Omega$, there exists a finite partition of unity $\{(U_k,\eta_k)\}$ of $\partial\Omega$, i.e. $\partial\Omega$ $=$ $\cup U_i$, $\mathbf{1}_{\partial\Omega}$ $=$ $\sum \eta_k$, $\eta_k$ $\in$ $C^\infty_0(U_k)$, such that each $U_k$ is diffeomorphic to the unit ball $B_1'$ of $\RR^{n-1}$ by a diffeomorphism $\psi_k:$ $B_1'$ $\rightarrow$ $U_k$. For $\tilde u_k$ $=$ $u\circ \psi_k$ and $X_k$ $=$ $\eta_k\,X$, we compute
\[
\int_{\partial\Omega} u\,X(u)\,d\sigma(x)
	= \frac{1}{2}\sum_k \int_{U_k} X_k\,(u^2)\,d\sigma(x) = \frac{1}{2}\sum_k \int_{B_1'} \sum_{i=1}^{n-1}\tilde X_i^{(k)}\,(\tilde u_k^2)_i\,dy,
\]
where $\|\tilde X_i^{(k)}\|_{C^1}$ is bounded from above by a constant depending only on $\|a_{ij}\|_{C^1}$, $\|\beta\|_{C^1}$ and $\mu$. Integrating by parts and noting that $X_i^{(k)}$ has compact support, we obtain
\[
\int_{\partial\Omega} u\,X(u)\,d\sigma(x)
	\leq C\sum_k\int_{B_1'} |\tilde u_k|^2\,dy \leq c_4\int_{\partial \Omega} |u|^2\,d\sigma(x),
\]
where $c_4$ depends on $\|a_{ij}\|_{C^1}$, $\|\beta\|_{C^1}$ and $\mu$. Returning to \eqref{LSF::Bound1} we hence get
\begin{align*}
c_1\int_\Omega |\nabla u|^2\,dx + (\lambda - c_2)\int_\Omega u^2\,dx
	&\leq (c_3 + c_4)\int_{\partial \Omega} u^2\,d\sigma(x)\\
	&\leq \frac{1}{2}c_1 \int_\Omega |\nabla u|^2\,dx + c_5\int_\Omega u^2\,dx
\end{align*}
where in the second estimate we have used the compactness of the embedding $H^1(\Omega)$ $\hookrightarrow$ $L^2(\partial\Omega)$. It is readily seen that, for $\lambda$ $>$ $c_2 + c_5$, the above inequality forces $\|u\|_{L^2(\Omega)}$ $=$ $0$, which implies $u$ $=$ $0$.
\eproof

\begin{proposition}\label{LinearLeading}
Assume that $(F,B)$ $=$ $(F_1, B_1) + (F_2, B_2)$ where
\begin{align*}
(F_1[u],B_1[u])
	&= \Big(a_{ij}(x)\,u_{ij} + b_i(x)\,u_i + c(x)\,u,\big(\beta_i(x)\,u_i + \gamma(x)\,u\big)\big|_{\partial\Omega}\Big),\\
(F_2[u],B_2[u])
	&= \Big(f_*(x,u,\nabla u), b_*(\cdot,u)\big|_{\partial\Omega}\Big),
\end{align*}
$a_{ij}$, $b_{i}$, $c$ $\in$ $C^{3,\alpha}(\bar\Omega)$, $\beta$, $\gamma$ $\in$ $C^{4,\alpha}(\partial\Omega)$, $f_2$ $\in$ $C^{3,\alpha}(\bar\Omega \times \RR \times \RR^n)$, and $b_2$ $\in$ $C^{4,\alpha}(\partial\Omega \times \RR)$. Assume that $F$ is elliptic, i.e. $(a_{ij})$ $>$ $0$ in $\bar\Omega$, and $B$ is oblique, i.e. $\beta \cdot \nu$ $>$ $0$ on $\partial\Omega$.

If $(F_1, B_1)$ is invertible, then for any open bounded set $\mathcal{O}$ $\subset$ $C^{4,\alpha}(\bar\Omega)$ such that $\partial\mathcal{O} \cap (F,B)^{-1}(0)$ $=$ $\emptyset$, we have
\[
\deg((F,B),\mathcal{O},0) = (-1)^{\dim E^-(F_1,B_1)}\deg_{L.S.}(Id + (F_1,B_1)^{-1}\circ(F_2,B_2),\mathcal{O},0),
\]
where
\[
E^-(F_1,B_1) = \bigoplus_{\lambda_i < 0} \Big\{u \in C^{4,\alpha}(\bar\Omega): -(F_1[u],B_1[u]) = (\lambda_i\,u,0)\Big\}.
\]
\end{proposition}

\begin{remark}
(a) $\dim E^-(F_1,B_1)$ is finite due to Lemma \ref{SemiFiniteness} and standard elliptic estimates.

(b) If $f_*$ $\equiv$ $0$ and $b_*$ $\equiv$ $0$, then the conclusion simplifies to
\[
\deg((F,B),\mathcal{O},0) = (-1)^{\dim E^-(F_1,B_1)}.
\]
\end{remark}

\bproof As before, set $T$ $=$ $(S\circ F,B):$ $C^{4,\alpha}(\bar\Omega)$ $\rightarrow$ $C^\alpha(\bar\Omega) \times C^{1,\alpha}(\partial\Omega) \times C^{3,\alpha}(\partial\Omega)$, where $S$ is given by \eqref{SDef}. Define $L:$ $C^{4,\alpha}(\bar\Omega)$ $\rightarrow$ $C^\alpha(\bar\Omega) \times C^{1,\alpha}(\partial\Omega) \times C^{3,\alpha}(\partial\Omega)$ by
\[
L\,w = \Big(a_{st}(x)\,w_{iist} - N\,a_{st}\,w_{st}, \Big(a_{st}(x)\,w_{sti}\,\nu_i\Big)\Big|_{\partial\Omega}, \Big(\beta_i(x,u)\,w_i + w\Big)\Big|_{\partial\Omega}\Big).
\]
By \eqref{DegDef}, we can pick $N$ large enough such that $L$ is invertible, $L^{-1} \circ T:$ $C^{4,\alpha}(\bar\Omega)$ $\rightarrow$ $C^{4,\alpha}(\bar\Omega)$ is of the form $\textrm{Id} + \textrm{Compact}$ and
\[
\deg((F,B),\mathcal{O},0) = \deg_{L.S.}(L^{-1}\circ T,\mathcal{O},0).
\]

Set $G$ $=$ $(S\circ F_1, B_1)$. By our hypotheses, $G$ is invertible. Thus, by the product rule of the Leray-Schauder degree, 
\[
\deg((F,B),\mathcal{O},0) = \sum_{\mathcal{U}}\deg_{L.S.}(L^{-1}\circ G,\mathcal{U},0)\,\deg_{L.S.}(G^{-1}\circ T,\mathcal{O},\mathcal{U}),
\]
where the summation is made over the connected components of $C^{4,\alpha}(\bar\Omega) \setminus (F_1,B_1)^{-1}\circ(F,B)(\partial\mathcal{O})$. It is evident that $\deg_{L.S.}(L^{-1}\circ G,\mathcal{U},0)$ $=$ $0$ if $0$ $\notin$ $\mathcal{U}$. Hence
\[
\deg((F,B),\mathcal{O},0) = \deg_{L.S.}(L^{-1}\circ G,\mathcal{\tilde O},0)\,\deg_{L.S.}(G^{-1}\circ T,\mathcal{O},0),
\]
where $\mathcal{\tilde O}$ is the connected component of $C^{4,\alpha}(\bar\Omega) \setminus (F_1,B_1)^{-1}\circ(F,B)(\partial\mathcal{O})$ containing $0$. As $G^{-1} \circ T$ $=$ $(F_1,B_1)^{-1} \circ (F,B)$ $=$ $Id + (F_1,B_1)^{-1}\circ(F_2,B_2)$, it remains to show that
\begin{equation}
d := \deg_{L.S.}(L^{-1}\circ G,\mathcal{\tilde O},0) = (-1)^{\dim E^-(F_1,B_1)}.
\label{LL::Eqn01}
\end{equation}

Define $(\tilde F, \tilde B):$ $C^{4,\alpha}(\bar\Omega)$ $\rightarrow$ $C^{2,\alpha}(\bar\Omega) \times C^{3,\alpha}(\partial\Omega)$ by
\[
(\tilde F[w],\tilde B[w]) = \Big(a_{ij}(\cdot)\,w_{ij},\big(\beta_i(\cdot)\,w_i + w\big)\big|_{\partial\Omega}\Big).
\]
For $0$ $\leq$ $t$ $\leq$ $1$, define $L_t:$ $C^{4,\alpha}(\bar\Omega)$ $\rightarrow$ $C^\alpha(\bar\Omega) \times C^{1,\alpha}(\partial\Omega) \times C^{3,\alpha}(\partial\Omega)$ by
\[
L_t\,w 	= \Big((1-t)\,a_{st}(x)\,w_{iist} + t\Delta \tilde F[w] - N\,\tilde F[w], \Big((1-t)a_{st}(x)\,w_{sti}\,\nu_i + t\,\frac{\partial \tilde F[w]}{\partial \nu}\Big)\Big|_{\partial\Omega}, \tilde B[w]\Big).
\]
As in the proof of \eqref{LInvertability}, we can select $N$ large enough so that each $L_t$ is an isomorphism for $t$ $\in$ $[0,1]$. Furthermore, as $L_t - G:$ $C^{4,\alpha}(\bar\Omega)$ $\rightarrow$ $C^{1,\alpha}(\bar\Omega) \times C^{2,\alpha}(\partial\Omega) \times C^{4,\alpha}(\partial\Omega)$, $L_t^{-1}\circ G$ is a legitimate homotopy for the Leray-Schauder degree. It follows that
\begin{equation}
d = \deg_{L.S.}(L^{-1}\circ G,\mathcal{\tilde O},0) = \deg_{L.S.}(L_1^{-1}\circ G,\mathcal{\tilde O},0).
\label{LL::Eqn02}
\end{equation}
Next, set
\[
\tilde L_t\,w = \Big(\Delta(\tilde F[w]) - (1-t)\,N\,\tilde F[w], \Big(\frac{\partial \tilde F[w]}{\partial \nu} + t\,\tilde F[w]\Big)\Big|_{\partial\Omega}, \tilde B[w]\Big).
\]
Arguing as before, we have $\tilde L_t$ is invertible and 
\begin{equation}
\deg_{L.S.}(L_1^{-1}\circ G,\mathcal{\tilde O},0) = \deg_{L.S.}(\tilde L_0^{-1}\circ G,\mathcal{\tilde O},0) = \deg_{L.S.}(\tilde L_1^{-1}\circ G,\mathcal{\tilde O},0).
\label{LL::Eqn03}
\end{equation}

Note that $\tilde L_1^{-1}$ $=$ $(S\circ \tilde F, \tilde B)$ and so $\tilde L_1^{-1} \circ G$ $=$ $(\tilde F,\tilde B)^{-1} \circ (F_1, B_1)$. Hence, by \eqref{LL::Eqn02} and \eqref{LL::Eqn03}
\[
d = \deg_{L.S.}((\tilde F, \tilde B)^{-1} \circ (F_1, B_1), \mathcal{O},0)
	= \deg_{L.S.}((F_1,B_1)^{-1} \circ (\tilde F, \tilde B), \mathcal{O},0).
\]
Set
\[
A_t = (F_1,B_1)^{-1} \circ \big[(1-t)(\tilde F, \tilde B) - t(Id,0)\big],
\]
where $(Id,0)$ is considered as an operator from $C^{4,\alpha}(\bar\Omega)$ into $C^{2,\alpha}(\bar\Omega) \times C^{3,\alpha}(\partial\Omega)$. By the maximum principle and obliqueness, $A_t$ is a continuous family of invertible linear operators acting on $C^{4,\alpha}(\bar\Omega)$. Moreover, for $t$ $\in$ $[0,1)$, $(1-t)^{-1}\,A_t$ is of the form $\textrm{Id} + \textrm{Compact}$. Hence, by the homotopy invariance property of the Leray-Schauder degree,
\[
d = \deg((1-t)^{-1}A_t,\mathcal{\tilde O},0) \text{ for any } t \in [0,1),
\]
which implies
\[
d = (-1)^{\dim E^-(A_t)} \text{ for any } t \in [0,1),
\]
where
\[
E^-(A_t) = \bigoplus_{\lambda_i(t) < 0} \Big\{u \in C^{4,\alpha}(\bar\Omega): A_t\,u = \lambda_i(t)\,u\Big\}.
\]

To proceed, we claim that there exists some $C$ $>$ $0$ and $\delta$ $\in$ $(0,1)$ such that, for any $t$ $\in$ $(1-\delta,1]$
\begin{equation}
-C < \lambda < -\frac{1}{C} \text{ for any negative eigenvalue $\lambda$ of $A_t$}.
\label{EigenBnd}
\end{equation}
Indeed, let $\lambda$ be an eigenvalue of some $A_t$ and $u$ be a corresponding eigenfunction. Since $A_t$ is invertible, Then
\[
\left\{\begin{array}{l}
a_{ij}\,u_{ij} + b_i\,u_i + c\,u = \frac{1}{\lambda}\big[(1-t)a_{ij}\,u_{ij} - tu\big] \text{ in } \Omega,\\
\beta_i\,u_i + \gamma\,u = \frac{1}{\lambda}\,(1-t)(\beta_i\,u_i + u) \text{ on } \partial\Omega,
\end{array}\right.
\]
which is equivalent to
\[
\left\{\begin{array}{l}
a_{ij}\,u_{ij} + \frac{\lambda}{\lambda - (1-t)}\,b_i\,u_i + \frac{\lambda}{\lambda - (1-t)}\,c\,u + \frac{t}{\lambda - (1-t)}u = 0 \text{ in } \Omega,\\
\beta_i\,u_i + \frac{\lambda}{\lambda - (1-t)}\gamma\,u - \frac{1-t}{\lambda - (1-t)}\,u = \text{ on } \partial\Omega,
\end{array}\right.
\]
It is readily seen that the first inequality in \eqref{EigenBnd} follows from the invertability of $(F,B)$ while the second follows from Lemma \ref{SemiFiniteness} for $\delta$ sufficiently small.

By \eqref{EigenBnd} and the compactness of $A_1$, we can pick a (simply connected) neighborhood $\mathcal{N}$ of $[-C,-\frac{1}{C}]$ in the complex plane such that in the set of eigenvalues of $A_1$ lying in $\mathcal{N}$ consists of all negative real eigenvalues of $A_1$. Furthermore, we can assume that $\mathcal{N}$ is symmetric about the real axis. Set
\begin{align*}
E(A_t,\mathcal{N})
	&= \bigoplus_{\lambda_i(t) \in \mathcal{N}}\Big\{u \in C^{4,\alpha}(\bar\Omega): A_t\,u = \lambda_i(t)u\Big\},\\
E^*(A_t,\mathcal{N})
	&= \bigoplus_{\lambda_i(t) \in \mathcal{N}\setminus \RR}\Big\{u \in C^{4,\alpha}(\bar\Omega): A_t\,u = \lambda_i(t)u\Big\}
\end{align*}
By the continuity of a finite system of eigenvalues (see e.g. \cite[pp. 213-214]{Kato}), for $\delta$ $>$ $0$ sufficiently small,
\[
\dim E(A_t,\mathcal{N}) \text{ is independent of } t \in (1-\delta,1].
\]
Also, since $A_t$ has real coefficients,
\[
\dim E^*(A_t,\mathcal{N}) \text{ is even }.
\]
Therefore, by \eqref{EigenBnd},
\[
d = \lim_{t \rightarrow 1}(-1)^{\dim E^-(A_t)} = \lim_{t \rightarrow 1}(-1)^{\dim E(A_t,\mathcal{N})} = (-1)^{\dim E(A_1,\mathcal{N})} = (-1)^{\dim E^-(A_1)}.
\]
As $A_1$ $=$ $-(F_1,B_1)^{-1}\circ(Id,0)$, \eqref{LL::Eqn01} follows. The proof is complete.
\eproof


\bibliography{paris}{}
\bibliographystyle{amsplain}
\end{document}